\documentclass[11pt,leqno]{amsart}

\usepackage{graphicx}
\usepackage[margin=1.15in]{geometry}
\usepackage{amsmath}
\usepackage{amssymb}
\usepackage{amsthm}
\usepackage{amsfonts}
\usepackage[export]{adjustbox}
\usepackage{graphicx}
\usepackage{enumitem}
\usepackage[scriptsize,hang,raggedright]{subfigure}
\usepackage[cmyk,usenames,dvipsnames]{xcolor}
\usepackage{mathtools}

\usepackage{hyperref}
\usepackage{pdfsync}
\usepackage{dsfont}
\usepackage{color}

\usepackage{cite}
\usepackage{bbm}
\usepackage{multirow}
\usepackage[font=footnotesize,labelfont=bf]{caption}

\numberwithin{equation}{section}

\newtheorem{theorem}{Theorem}[section]
\newtheorem{proposition}[theorem]{Proposition}
\newtheorem{lemma}[theorem]{Lemma}
\newtheorem{corollary}[theorem]{Corollary}

\newtheorem{definition}[theorem]{Definition}
\newtheorem{as}[theorem]{Assumption}
\theoremstyle{remark}
\newtheorem{remark}[theorem]{Remark}

\newcommand{\e}{\varepsilon}
\newcommand{\R}{\mathbb{R}}

\newcommand{\ba}{\begin{array}}
\newcommand{\ea}{\end{array}}

\newcommand{\bthm}{\begin{theorem}}
\newcommand{\ethm}{\end{theorem}}
\newcommand{\bprop}{\begin{proposition}}
\newcommand{\eprop}{\end{proposition}}
\newcommand{\blemma}{\begin{lemma}}
\newcommand{\elemma}{\end{lemma}}

\newcommand{\beqn}{\begin{equation}}
\newcommand{\eeqn}{\end{equation}}
\newcommand{\beqns}{\begin{equation*}}
\newcommand{\eeqns}{\end{equation*}}

\newcommand{\supp}{\operatorname{supp}}
\newcommand{\diam}{\operatorname{diam}}

\newcommand{\pr}{\prime}
\newcommand{\pt}{\partial}
\newcommand{\arrow}{\rightarrow}

\renewcommand{\leq}{\leqslant}
\renewcommand{\geq}{\geqslant}

\definecolor{mygreen}{rgb}{0.1,0.75,0.2}

\newcommand{\N}{\mathbb{N}}

\newcommand{\eps}{\epsilon}

\newcommand{\E}{\mathsf{E}}

\newcommand{\F}{\mathsf{F}}
\newcommand{\Se}{\mathsf{S}}
\newcommand{\Om}{\Omega}

\DeclareMathOperator{\id}{Id}

\newcommand{\Rd}{{\mathord{\mathbb R}^d}}
\newcommand{\ird}{\int_\Rd}
\newcommand{\loc}{{\rm loc}}
\newcommand{\grad}{\nabla}
\newcommand{\wto}{\rightharpoonup}
\newcommand{\wsto}{\stackrel{*}{\rightharpoonup}}
\newcommand{\la}{\langle}
\newcommand{\ra}{\rangle}
\newcommand{\C}{\mathcal{C}}
\newcommand{\K}{\mathsf{K}}
\def\P{{\mathcal P}}
\def\S{{\mathcal S}}
\newcommand{\bt}{\mathbf{t}}

\title[Aggregation-diffusion to constrained interaction]{Aggregation-diffusion to constrained interaction:  \\ minimizers \& gradient flows in the slow diffusion limit}

\author{Katy Craig}
\address{Department of Mathematics, University of California, Santa Barbara, CA}
\email{kcraig@math.ucsb.edu}
\author{Ihsan Topaloglu}
\address{Department of Mathematics and Applied Mathematics, Virginia Commonwealth University, Richmond, VA}
\email{iatopaloglu@vcu.edu}
\thanks{KC's research is supported by U.S. National Science Foundation grants DMS 1401867 and DMS 1811012.}
\date{\today}                                        

\begin{document}

\begin{abstract}
Inspired by recent work on minimizers and gradient flows of constrained interaction energies, we prove that these energies arise as the slow diffusion limit of well-known aggregation-diffusion energies. We show that minimizers of aggregation-diffusion energies converge to a  minimizer of the constrained interaction energy and gradient flows converge to a gradient flow. Our results apply to a  range of  interaction potentials, including singular attractive and repulsive-attractive power-law potentials. In the process of obtaining the slow diffusion limit, we also extend the well-posedness theory for aggregation-diffusion equations and Wasserstein gradient flows to admit a wide range of  nonconvex interaction potentials. We conclude by applying our results to develop a numerical method for constrained interaction energies, which we use to investigate open questions on  set valued minimizers.
\end{abstract}

\maketitle

\tableofcontents


\section{Introduction}\label{sec:intro}

Nonlocal interactions arise  throughout the natural world, from collective dynamics in biological swarms to vortex motion in superconductors and gravitational interactions among stars. In each case, agents experience pairwise attractive or repulsive forces, and these pairwise interactions  are often coupled with additional repulsive effects, such as diffusion or a height constraint, which penalize accumulations. 
The simplest mathematical model for nonlocal interactions and diffusion is the \emph{aggregation-diffusion equation},
\begin{align*}
\partial_t \rho - \grad \cdot ( (\grad K*\rho) \rho) = \Delta \rho^m ,   \quad K: \Rd \to \R , \quad m \geq 1,
\end{align*}
where the \emph{interaction potential} $K$ governs the pairwise interactions and the \emph{diffusion exponent} $m \geq 1$ controls the strength at which diffusion is felt  at different heights of the density $\rho$.
 Likewise, nonlocal interactions coupled with a height constraint can be heuristically modeled by the \emph{constrained aggregation equation}
\begin{align*} \label{constrainedagg1}
\begin{cases} \partial_t \rho - \grad \cdot ( (\grad K*\rho) \rho) = 0 & \text{ if } \rho <1 ,\\ \rho \leq 1  \text{ always}, &\end{cases}
\end{align*}
where, again, $K: \Rd \to \R$ is the interaction potential.
 We note that this equation is merely a heuristic partial differential equation, as we do not specify the sense in which the height constraint $\rho \leq 1$ is enforced. We provide a rigorous formulation below.

Both the aggregation-diffusion equation and constrained aggregation equation have  gradient flow structures with respect to the \emph{2-Wasserstein metric}. The aggregation-diffusion equation is formally the gradient flow of the sum of an interaction energy and R\'enyi entropy
\beqn \label{eqn:reg_energy}
		\E_m(\rho) \, = \, \begin{dcases*} 
											\frac{1}{2}\iint K(x-y)\,\rho(x) \rho(y) \,dx \,dy+ \frac{1}{m-1}\int \rho(x)^m\,dx & \quad if $\rho \in     L^m(\Rd) $, \\
											 + \infty & \quad otherwise,
									 \end{dcases*}
	\eeqn
and the constrained aggregation equation can be rigorously posed as the gradient flow of the constrained interaction energy
	\beqn\label{eqn:const_energy}
		\E_\infty(\rho) \, = \,  \begin{dcases*} 
											\frac{1}{2}\iint K(x-y)\,\rho(x) \rho(y) \,dx \,dy& \quad  if $\rho\in L^\infty(\Rd) $ and $\|\rho\|_{\infty} \leq 1$, \\
										    + \infty & \quad otherwise.
									  \end{dcases*}
	\eeqn

Over the past fifteen years, there has been significant work on  aggregation-diffusion equations, analyzing dynamics of solutions, asymptotic behavior, and minimizers of the energy $\E_m$ \cite{CaHiVoYa16,Bed2010,Bed2011, BertozziSlepcev10,BlanchetCalvezCarrillo, BlanchetDolbeaultPerthame, CalvezCarrillo, CarrilloMcCannVillani2, CarrilloMcCannVillani,KimYao,Carrillo:2015ke, carrillo2017ground,YaoBertozzi,yao2014asymptotic,chayes2013aggregation}. The vast majority of the literature has considered one of two choices of interaction potential: either purely attractive power-laws or repulsive-attractive power-laws,
\begin{align} \label{powerlaweqn0} K(x)=  |x|^{p}/p  \quad \text{ or } \quad K(x) = |x|^q/q - |x|^p/p \quad \text{for} \quad 2-d\leq  p<q\leq 2 , \quad q >0 , 
\end{align}
with the convention that $|x|^0/0 = \log(|x|)$.
For the purely attractive case, the literature has largely studied the competition between the attraction parameter $p$ and the diffusion exponent $m$,   along with the effects this competition has on properties such as global existence of  solutions or finite time blowup; see \cite{BlCaLa09, BCC,sugiyama2006global,CaSu16,calvez2017geometry,calvez2017equilibria,BlanchetCarrilloMasmoudi,BedrossianRodriguezBertozzi,BL13,BlanchetDolbeaultPerthame}. For the repulsive-attractive case, the requirement    $p < q$ ensures that the nonlocal interactions are repulsive at short length scales and attractive at long length scales. This competition between short-range and long-range effects  leads to rich pattern formation in both the steady states of solutions and the minimizers of the corresponding energy $\E_m$; see \cite{ChFeTo14,CraigTopaloglu,SiSlTo2014,CrBe14,BaSh15,CarrilloCraigPatacchini,carrillo2017uniform,Bertozzietal_RingPatterns,BalagueCarrilloYao}.

More recently, several works have also considered the constrained aggregation equation and minimizers of the constrained interaction energy $\E_\infty$. Minimizers of $\E_\infty$  are directly related to a shape optimization problem introduced by Burchard, Choksi, and the second author  \cite{BuChTo2016}: given a repulsive-attractive power-law interaction potential $K$, as in equation (\ref{powerlaweqn0}),
	\begin{align} \label{setvaluedproblem}		\text{minimize} \qquad \E(\Omega) \,=\, \frac{1}{2}\int_{\Omega} \! \int_{\Omega} K(x-y)\,dxdy \quad \text{ over sets } \Omega \subseteq \Rd \text{ of  volume } M .
\end{align}
Competition between the attraction parameter $q$  and the repulsion parameter $p$ in the definition of $K$ determines existence, nonexistence, and qualitative properties of minimizers, providing a counterpoint to the well-studied nonlocal isoperimetric problem. (See \cite{ChMuTo2017} for a survey.) 
Burchard, Choksi, and the second author showed that the shape optimization problem admits a  solution if and only if  the constrained interaction energy $\E_\infty$ admits a \emph{set valued minimizer}, i.e., a minimizer $\rho$ that is a characteristic function of a set $\Omega$, $\rho=\chi_\Omega$. 
Furthermore, they proved that for attraction $q=2$ and repulsion $-d<p<0$,  there are critical values of the mass $M_1 \leq M_2$  so that set valued minimizers  with mass $M$  exist for $M \geq M_2$ and do not exist for $M \leq M_1$. Subsequently, Frank and Lieb extended this result to $q >0$ and $p = 2-d$ and  proved that there are also critical values of the mass that separate the \emph{liquid} and \emph{solid} phases of minimizers of $\E_\infty$:  if $M< M_1^*$, then minimizers of $\E_\infty$   satisfy $|\{ \rho = 1\}| = 0$ (liquid), and if $M> M_2^*$, then $|\{ \rho = 1\}| = M$ (solid) \cite{FrLi2016} . On one hand, it is known that  
\begin{align} \label{criticalmassinequality}
 M_1^* \leq M_1 \leq M_2 \leq M_2^*,
\end{align}
and for Newtonian repulsion and quadratic attraction ($p=2-d$, $q =2$, $d>2$), all four values  equal 1. 
On the other hand, Lopes provided an explicit example for which $M_1^* < M_2^*$ \cite{Lopes}. In general, it remains unknown   for which values of   $p$ and $q$ strict inequality holds in any of the three inequalities in (\ref{criticalmassinequality}), as well as how the values of the  critical masses depend on $p$ and $q$.

Concurrently with this work on minimizers of the constrained interaction energy $\E_\infty$, Kim, Yao, and the first author studied gradient flows of $\E_\infty$, which formally solve the constrained aggregation equation \cite{CrKiYa16}. This work was inspired by the vast literature on height constrained problems, which arise in both models of crowd motion and tumor growth (see, e.g., \cite{MaRoSa10,MaRoSaVe11,chizat2016scaling,chizat2017tumor,perthame2014hele, mellet2017hele,kim2018hele}).  In the case of a purely attractive Newtonian interaction potential (equation (\ref{powerlaweqn0}) with $p = 2-d$), they characterized gradient flows of $\E_\infty$ with set valued initial data   in terms of a Hele-Shaw type free boundary problem. A key element of their proof was that, formally, gradient flows of $\E_m$ converge to gradient flows of $\E_\infty$ as $m \to +\infty$. Indeed, Alexander, Kim, and Yao had proved the analogous results for \emph{drift} diffusion equations  in previous work \cite{AKY}. However, in the case of \emph{aggregation} diffusion equations, rigorous analysis of this limit was not considered, due to the lack of convexity of the interaction potential $K$.

The objective of the present work is to prove that, indeed, minimizers and gradient flows of $\E_m$ do converge to minimizers and gradient flows $\E_\infty$ in the \emph{slow diffusion limit }as $m \to +\infty$. We consider measures with a fixed mass $M>0$, and without loss of generality, we   rescale  so that $M =1$. For a general class of interaction potentials $K$,  including  both attractive  and repulsive-attractive power-law potentials (\ref{powerlaweqn0}), we prove that minimizers of $\E_m$ converge to a minimizer of $\E_\infty$ (up to a subsequence and translations) and  gradient flows of $\E_m$ converge to a gradient flow of $\E_\infty$ (up to a subsequence)  in the weak-* topology of probability measures. The latter result extends the famous Mesa Problem for the porous medium equation to include a singular nonlocal interaction term (see e.g. \cite{CaffarelliFriedman}).

In the process of proving these results, we also rigorously prove the equivalence between solutions of  aggregation-diffusion equations and gradient flows of the energies $\E_m$. Likewise, we extend the well-posedness theory for aggregation-diffusion equations  to include singular repulsive-attractive power-law potentials, thereby filling a gap in the existing theory. Finally, we  succeed in characterizing the minimal subdifferential of $\E_\infty$ along the gradient flow, a key quantity in the study of gradient flows, which was  identified formally in previous works on constrained energies \cite{MaRoSa10,MaRoSaVe11}.
We believe that one of our main contributions  is the extension of the theory of Wasserstein gradient flows  to energies, such as $\E_m$ and $\E_\infty$, that   satisfy neither the classical $\lambda$-convexity assumption of Ambrosio, Gigli, and Savar\'e \cite{AGS} nor the more recent $\omega$-convexity assumption  \cite{CarrilloLisiniMainini, craig2017nonconvex, AmbrosioSerfaty, AmbrosioMaininiSerfaty, CarrilloMcCannVillani2}.

Finally, we apply these theoretical results to develop a numerical method for   gradient flows and minimizers of the constrained interaction energy $\E_\infty$. We use Carrillo, Patacchini, and the first author's \emph{blob method for diffusion} (see \cite{CarrilloCraigPatacchini}) to simulate gradient flows and minimizers of $\E_m$ for $m$  large,  thereby approximating the corresponding gradient flows and minimizers of $\E_\infty$. While there exist other  numerical methods for constrained problems---such as  Liu, Wang, and Zhou's method for purely attractive Newtonian interactions  \cite{liu2018positivity} and several Eulerian methods based on the JKO scheme \cite{chizat2017tumor, carlier2017convergence, peyre2015entropic, gallouet2017unbalanced}---our particle method is unique in its ability to resolve the nonlocal interaction term for a range of interaction potentials $K$. 
As the primary goal of the present work is theoretical analysis of the slow diffusion limit, we restrict our numerical study to one dimension, though our  method  naturally extends to all dimensions $d \geq 1$. 

 We conclude with  several numerical simulations that shed light on  open questions for minimizers of the constrained interaction energy. These numerical results indicate that the critical values of the mass $M_1$ and $M_2$ that separate   nonexistence and existence of set valued minimizers of $\E_\infty$ are in fact equal, and we explore how $M_1 = M_2$ depends on the attraction and repulsion parameters $q$ and $p$. We also observe that, for $p=1$, the  critical masses $M_1^*$ and $M_2^*$ that separate the liquid and solid phases are in general not equal, except for  $q=2$, so that the existence of an intermediate phase is indeed the generic behavior for minimizers of the constrained interaction energy.  

We now describe the  assumptions we impose on the interaction potentials $K$ and then provide  a precise statement of our main results. We conclude the introduction with an outline of our approach and a brief summary of our notation.

\subsection{Assumptions on Interaction Potentials}

We  impose the following assumptions on the interaction potential $K$ and diffusion exponent $m$. To ensure lower semicontinuity of the energies $\E_m$ and $\E_\infty$ with respect to weak-* convergence of  measures, we suppose that $m \geq m_0$, where $m_0$ and the interaction potential satisfy the following condition.

	 \medskip
	
		\begin{enumerate}[label= (LSC)] \addtolength{\itemsep}{6pt}

\item $K:\Rd \to [-\infty, \infty]$ is even, locally integrable,  and $K = K_{a} + K_{b}$ for two lower semicontinuous functions $K_a$ and $K_b$, where $K_a$ is bounded below and $K_b \in  L^{r, \infty}(\Rd)$, for  $r \in (1,+\infty)$; the  lower bound on the diffusion exponent satisfies $m_0>1 + 1/r$. \label{GFcont}
		\end{enumerate}

\medskip		
		
\begin{remark}[Diffusion Dominated Regime]
In the case of an attractive power-law interaction potential, $K(x) = |x|^p/p$ for $-d< p < 0$, hypothesis \ref{GFcont} is equivalent to the requirement that we are in the diffusion dominated regime, $m_0 > 1- p/d$ (cf. \cite{BedrossianRodriguezBertozzi,carrillo2017ground,sugiyama2006global}). More generally, for repulsive-attractive power-law interaction potentials $K(x)=|x|^q/q-|x|^p/p$ with $-d < p < q $, hypothesis \ref{GFcont} merely requires that $m_0 >1$. (See Proposition \ref{GFcontprop}.)  \label{ddremark}
\end{remark}

\medskip

In order to establish the existence of compactly supported minimizers and prove that minimizers of $\E_m$ converge to a minimizer of $\E_\infty$ as $m \to +\infty$, we impose the following assumptions on the regularity and growth of the interaction potential.

	\medskip
	
\begin{enumerate}[label = (ATT)] 
	\item Either $K$ is purely attractive and  approaches some constant $\ell \in \mathbb{R}$ at infinity or $K$ grows to infinity at infinity. Namely, either
			\medskip
							\begin{itemize} \addtolength{\itemsep}{6pt}
								\item[(i)] $\lim_{|x| \to +\infty} K(x) = \ell$, and $\ell-K \in L^p(\Rd\setminus B)$ for some $1\leq p < \infty$;\\   $K \in C^1(\Rd\setminus\{0\})$, $\pt_{|x|} K \geq 0$, and  for all $|x|>1$, $\pt_{|x|} K \leq C$; or,  \label{weakatt}  
								\item[(ii)]   $\lim_{|x| \to +\infty} K(x) = +\infty$.						\label{strongatt}
							\end{itemize} 
							\label{GFgrowth}  
					\end{enumerate}
					
	\medskip

To prove that gradient flows of $\E_m$ converge to a gradient flow of $\E_\infty$, we impose the following assumptions on the growth, regularity, and stability of $K*\rho$ and $\grad K *\rho$ for all $\rho \in \P_2(\Rd) \cap L^{m}(\Rd)$, where $\P_2(\Rd)$ denotes the set of probability measures with finite second moment, $M_2(\rho) = \int |x|^2 d \rho(x) < +\infty$.						
We assume that there exists  a constant $C\geq1$ and a continuous, nondecreasing, concave function $\psi:[0,+\infty)\to[0,+\infty)$ with  $\psi(0) = 0$ so that for all $m \geq m_0$ and $\rho,\, \mu, \, \nu \in\P_2(\Rd) \cap L^{m}(\Rd)$, 
\pagebreak
	\medskip
							\begin{enumerate}[label = (GF\arabic*)]
	\addtolength{\itemsep}{6pt}

						\item    $\|\nabla K * \rho\|_{L^2(\nu)} \leq C(1+\|\rho\|_{{m}} +M_2(\rho)^{1/2} + M_2(\nu)^{1/2} )$; \label{GFbounds}

						\item   $K*\rho \in C^1(\Rd)$ and $	\big| \nabla K * \rho(x) - \nabla K * \rho (y) \big|^2 \leq C(1+ {\|\rho\|^2_{{m}}} ) \psi \big(|x-y|^2\big) ; $  \label{GFcty}

						\item $\| \nabla K*(\rho -  \nu) \|_{L^2(\mu)} \leq C\Big(1+\|\rho\|_{{m}} +   \|\nu\|_{{m}} +\|\mu\|_{{m}}\Big) \psi( d_{W_{2- \epsilon}}(\rho,\nu))$ for some $\epsilon \in (0,1)$. \label{GFdualsobolev}
\end{enumerate}

\begin{remark}[differentiability vs convexity]
Hypothesis  \ref{GFcty} is weaker than the analogous hypotheses in previous work \cite{CraigTopaloglu, craig2017nonconvex}, since in the present context we merely require differentiability of $\E_m$ instead of convexity (or $\omega$-convexity) of $\E_m$.
\end{remark}

Hypotheses \ref{GFcont}, \ref{GFgrowth}, and \ref{GFbounds}-\ref{GFdualsobolev} are satisfied by the attractive and repulsive-attractive power law   potentials described in the introduction (\ref{powerlaweqn0}); see Theorem \ref{powerlawhypothesesthm}.

\medskip

\subsection{Main Results}

Our first main result establishes the convergence of energy minimizers.

\bthm[minimizers weak-* converge to minimizer] \label{thm:gamma_conv} Suppose $K$ satisfies  hypotheses \ref{GFcont} and \ref{GFgrowth}.  Then for any sequence $\rho_{m} \in \P_2(\Rd)$ of minimizers of $\E_m$, there  exists $\rho \in \P_2(\Rd)$ so that, up to a subsequence and translations, $\rho_{m} \wsto \rho$ in $\P(\Rd)$   and $\rho$ minimizes $\E_\infty$. 
\ethm

\begin{remark}[existence of minimizers of $\E_m$]
We use hypothesis \ref{GFgrowth}  to conclude existence of minimizers of $\E_m$. If one were able to obtain existence by other means, hypothesis \ref{GFcont} is sufficient to conclude the $m \to +\infty$ limit.
\end{remark}

For attractive or repulsive-attractive power-law potentials, we adapt the arguments by Rein \cite{Re2001} and Frank and Lieb \cite{FrLi2016}, respectively, to prove that minimizers of the energies $\E_m$ are compactly supported \emph{uniformly} in $m$. 

\bthm[uniform bound on support] \label{thm:unif_comp_supp}
Let $\rho_m$ be a minimizer of the energy $\E_m$. For 
	\beqn \label{eqn:kernels_powerlaw}
			K(x) = |x|^p/p \quad \text{ with } -d<p<0 \quad \text{ or } \quad K(x)=\frac{1}{q}|x|^q - \frac{1}{p}|x|^p \quad \text{ with }-d<p<0<q
	\eeqn
there exists $R>0$ so that $\supp \rho_m \subset B_R(0)$ for all $m>1$ sufficiently large. 
\ethm

As a consequence of the previous two theorems, we obtain the convergence of minimizers of $\E_m$ to a minimizer of $\E_\infty$  in the stronger 2-Wasserstein distance.

\begin{corollary}[minimizers converge to minimizer] \label{cor:P2_conv}
For interaction potentials $K$ of the form \eqref{eqn:kernels_powerlaw}, any sequence of minimizers of $\E_m$ converges, up to a subsequence and translations, to a minimizer of $\E_\infty$ in the 2-Wasserstein metric.
\end{corollary}

We next turn our attention to  gradient flows of the energies $\E_m$ and $\E_\infty$. We begin by showing that, for $m$ sufficiently large, gradient flows of $\E_m$ exist and solve the aggregation-diffusion equation, for all  initial data in the domain of the energy $D(\E_m) = \{ \rho \in \P_2(\Rd) : \E_m(\rho) < +\infty \}$.

\begin{theorem}[well-posedness of gradient flows] \label{GFexistthm}
Suppose $K$ satisfies hypotheses \ref{GFcont} and \ref{GFbounds}--\ref{GFdualsobolev} and $m \in [ m_0, +\infty]$. 
\begin{enumerate}[label = (\roman*)]
\item For all $\rho_m^{(0)} \in D(\E_m)$, the gradient flow of $\E_m$ with initial data $\rho_m^{(0)}$ exists. 
\item If the modulus $\psi(s)$ in \ref{GFcty}--\ref{GFdualsobolev} satisfies $\psi(s) \geq s$ and $\int_0^1 (s\psi(s))^{-1/2} ds = +\infty$, then the gradient flow is unique.
\item For $m < +\infty$ and $\rho_m^{(0)} \in D(\E_m)$,  $\rho_m(t)$ is a gradient flows of $\E_m$ if and only if it solves the aggregation-diffusion equation in the duality with $C^\infty_c(\Rd\times [0,T]) $,
\begin{align*} \partial_t \rho_m  + \grad \cdot ((\grad K*\rho_m)\rho_m) = \Delta \rho_m^m  , \quad \rho_m(0) = \rho_m^{(0)}  
\end{align*}
\end{enumerate}
\end{theorem}

\begin{remark}[existence] \label{existenceremark}
In the particular case that $K$ is a singular attractive power-law potential, $K(x) = |x|^p/p$ for $2-d \leq p \leq 0$, and the diffusion exponent is sufficiently large, $m \geq m_0>   \max \{   d/(d+p-1), 1\}$, the previous theorem extends the range of initial data $\rho_0$ for which it is known that  solutions to the aggregation-diffusion equation  exist from $\rho_0 \in L^\infty(\Rd)\cap \P(\Rd)$ to $\rho_0 \in D(\E_m)$; see \cite{BedrossianRodriguezBertozzi, BedrossianRodriguez, sugiyama2006global}. In Proposition \ref{GFpropprop}, we also strengthen the energy dissipation \emph{inequality} from these previous works to an energy dissipation \emph{identity}.

To our knowledge, the previous theorem provides the first existence  results for agg\-re\-ga\-tion-diffusion equations with singular repulsive-attractive power-law potentials of the form $K(x) = |x|^q/q - |x|^p/p$, $2-d \leq p < q \leq 2$, provided that $m \geq m_0> \max \{   d/(d+p-1), 1\}$. This complements recent work by Carrillo and Wang, which studied global boundedness of solutions, under the assumption that solutions exist locally in time  \cite{carrillo2017uniform}.
\end{remark}

\begin{remark}[uniqueness] \label{uniquenessremark}
If   $K(x) = |x|^p/p$ or $K(x) = |x|^q/q - |x|^p/p$ for $2-d \leq p < q \leq 2$ and $m_0 \geq d/(p+d-2)$, then we may take $\psi(s) = s |\log(s)|$ for $s$ near zero in hypotheses \ref{GFcty}--\ref{GFdualsobolev}; see \cite[Proposition 4.4]{craig2017nonconvex}. Consequently, the gradient flow of $\E_m$ is unique for $m \geq m_0 \geq d/(p+d-2)$ and the gradient flow of $\E_\infty$ is unique.
\end{remark}

We  apply this result to show that, up to a subsequence,  gradient flows of $\E_m$ with well-prepared initial data converge to a gradient flow of $\E_\infty$ as $m \to +\infty$.

\bthm[subsequence of gradient flows converges to gradient flow]\label{thm:grad_flow_diff}
Suppose $K$ satisfies hypotheses \ref{GFcont} and \ref{GFbounds}--\ref{GFdualsobolev}. Let $\rho_{m}(t) $ be a gradient flow of $\E_m$. Suppose  that the initial data $\rho_{m}^{(0)}$ is well-prepared: $\sup_m M_2(\rho_{m}^{(0)}) < +\infty$  and for some $\rho^{(0)}\in D(\E_\infty)$
	\[
		\rho_{m}^{(0)}\wsto\rho^{(0)} \text{ weak-* in $\P(\Rd)$ and } \lim_{m\to\infty  } \E_m(\rho_{m}^{(0)}) = \E_\infty(\rho^{(0)}).
	\]
Then $\rho_{m} $ has a weak-* convergent subsequence so that $\rho_m(t) \wsto \rho(t)$ for almost every $t \geq 0$, and $\rho(t)$ is  a gradient flow of $\E_\infty$ with initial data $\rho^{(0)}$. Furthermore, as $m\to\infty$,
	\begin{gather*}
 \E_m(\rho_{m}(t)) \to \E_\infty(\rho(t)) \text{ for all }  \ t \geq 0;\\
|\pt \E_m|(\rho_{m}) \to |\pt \E_\infty|(\rho)  \text{ and } \ |\rho^\pr_{m}|_{d_W} \to |\rho^\pr|_{d_W} \text{ in } L^2_\loc(0,+\infty).
\end{gather*}

Finally, for almost every $t\geq 0$, there exists
\begin{align} \label{sigma}
\sigma(t) \in H^1(\Rd) \text{ satisfying }\sigma(t) \geq 0 \text{ and }\sigma(t) =0 \text{ almost everywhere on } \{ \rho(t) < 1\} ,
\end{align} 
so that, up to a subsequence,  $\rho_m^m(t) \wto \sigma(t)$ in $L^2( \Rd)$ and $\grad K*\rho(t)+ \grad \sigma(t)/\rho(t)$  is the element of $\partial E_\infty(\rho(t))$ with minimal $L^2(\rho(t))$ norm.
\ethm

\begin{remark}[minimal subdifferential of $\E_\infty$]
A byproduct of our result on the convergence of gradient flows is that we are able to characterize the minimal element of the subdifferential of $\E_\infty$ along the gradient flow. This result can be easily extended to allow  a $\lambda$-convex drift potential $V(x)$  in the energies $\E_m$ and $\E_\infty$ by adding term of the form $\partial V$ to the subdifferential. This makes rigorous formal characterizations of the subdifferential from previous works \cite{MaRoSa10,MaRoSaVe11}.
\end{remark}

For attractive or repulsive-attractive power-law potentials, we may use uniqueness of the gradient flow of $\E_\infty$ (see Remark \ref{uniquenessremark}) to immediately obtain a stronger convergence result. 

\begin{corollary}[gradient flows converge to gradient flow] \label{gammaconvGFstrong}
Given an interaction potential 
\begin{align*} K(x)=  |x|^{p}/p  \quad \text{ or } \quad K(x) = |x|^q/q - |x|^p/p \quad \text{for} \quad 2-d\leq  p<q\leq 2 ,  
\end{align*}
consider  gradient flows $\rho_m(t)$ of $\E_m$ with well-prepared initial data: $\sup_m M_2(\rho_{m}^{(0)}) < +\infty$  and for some $\rho^{(0)}\in D(\E_\infty)$
	\[
		\rho_{m}^{(0)}\wsto\rho^{(0)} \text{ weak-* in $\P(\Rd)$ and } \lim_{m\to\infty  } \E_m(\rho_{m}^{(0)}) = \E_\infty(\rho^{(0)}).
	\]
Then $\rho_m(t) \wsto \rho(t)$ for almost every $t\geq0$, where $\rho(t)$ is the unique gradient flow of $\E_\infty$ with initial data $\rho^{(0)}$.
\end{corollary}

\subsection{Outline and Notation}
The remainder of the paper is organized as follows. In section \ref{sec:background}, we recall fundamental results on the Wasserstein metric, gradient flows, and $\Gamma$-convergence. In section \ref{examplessection}, we prove that attractive and repulsive-attractive power-law  potentials satisfy our main hypotheses \ref{GFcont}, \ref{GFgrowth}, and \ref{GFbounds}-\ref{GFdualsobolev} (Theorem \ref{powerlawhypothesesthm}). In section \ref{energiessection}, we  prove that the energies $\E_m$ and $\E_\infty$ are lower semicontinuous with respect to weak-* convergence of measures (Proposition \ref{LSCproposition}) and bounded below, uniformly in $m$ (Proposition \ref{firstLmbound}). We then characterize the minimal element of the subdifferential of $\E_m$ (Proposition \ref{subdifftheorem}), identify an element of subdifferential of $\E_\infty$ (Proposition \ref{Einftysubdiff}), and prove well-posedness of gradient flows (Theorem \ref{GFexistthm}). In section \ref{sec:proof1}, we prove our main results on convergence of minimizers (Theorem \ref{thm:gamma_conv}) and the uniform bound on the support of minimizers (Theorem \ref{thm:unif_comp_supp}). In section \ref{sec:proof2}, we prove our main result on convergence of gradient flows (Theorem \ref{thm:grad_flow_diff}). Finally, in section \ref{sec:numerics}, we apply these theoretical results to develop a numerical method for simulating gradient flows and minimizers of $\E_\infty$, which we use to explore the open questions about minimizers of $\E_\infty$ described in the introduction.

We conclude by briefly reviewing our notation. When a probability measure $\rho \in \P(\Rd)$ is absolutely continuous with respect to Lebesgue measure, $\rho \ll \mathcal{L}^d$,  we commit a mild abuse of notation and  denote both the measure and its density by $\rho$,  $d \rho = \rho(x) dx$. Differentials in integrals will likewise be written either as $d\rho(x)$ or $\rho(x)dx$, depending on the context. Norms with respect to Lebesgue measure will be denoted by single subscripts (e.g., $\|\cdot\|_p$) whereas $L^p$-norms with respect to a measure $\mu \in \P(\Rd)$ will be explicitly marked (e.g., $\|\cdot\|_{L^p(\mu)}$). We denote  convergence with respect to the weak-* topology   by $\wsto$.  For measures that depend on time $\mu(t) \in \P(\Rd)$, we commit a mild abuse of notation and identify 
\[ (\mu_t)_{t \in (0,1)} \sim \int_0^1 \delta_t \otimes \mu_t dt \text{ so that } \iint f(t,x) d \mu = \int_0^1 \int_\Rd f(t,x) d \mu_t(x) dt . \]
We let $\chi_\Om$ denote the characteristic function on a set $\Om \subseteq \Rd$ and $\Om^c$ denote the complement of $\Omega$.
We allow all constants $C>0$ to change from line to line.

\section{Preliminaries}\label{sec:background}

\subsection{The Wasserstein Metric} \label{sec:wasserstein}
For $b \in [1,2]$, we consider measures belonging to the space
	\[
		\P_b(\Rd) := \left\{ \mu\in\P(\Rd) \colon \int |x|^b\,d\mu(x) < +\infty \right\}
	\] 
of probability measures with finite $b$th moments. We endow this space with the $b$-Wasserstein metric, which we recall briefly now. For further background, we refer the reader to the books by Ambrosio, Gigli and Savar\'e \cite{AGS} and Villani \cite{Villani}.

The $b$-Wasserstein distance between $\mu, \nu \in \P_b (\Rd)$ is given by
	\beqn\label{eqn:Wass_dist}
		d_{W_b}(\mu,\nu):=\left(\min\left\{\iint |x-y|^b\,d\gamma(x,y)\colon \gamma\in\C(\mu,\nu)\right\}\right)^{1/b},
	\eeqn
where $\C(\mu,\nu)$ is the set of transport plans between $\mu$ and $\nu$,	
	\[
		\C(\mu,\nu):=\left\{\gamma\in\P(\Rd\times\Rd)\colon (\pi_1)_\#\gamma=\rho\quad\text{and}\quad(\pi_2)_\#\gamma=\nu\right\}.
	\]
Here $\pi_1$, $\pi_2$ denote the projections $\pi_1(x,y)=x$ and $\pi_2(x,y)=y$. For $i=1,2$, $(\pi_i)_{\#}\gamma$ denotes the pushforward of $\gamma$ defined by $(\pi_i)_{\#}\gamma (U):=\gamma(\pi_i^{-1}(U))$ for any measurable set $U\subset\Rd$. By H\"older's inequality for the probability measure $\gamma \in \P(\Rd \times \Rd)$, we have
\begin{align} \label{Wasserstein order}
 d_{W_b}(\mu,\nu) \leq d_{W_a}(\mu,\nu) , \quad \text{ for all } b \leq a . 
 \end{align}

For any $\mu, \nu \in \P_2(\Rd)$, the minimization problem \eqref{eqn:Wass_dist} admits a solution: there exists an optimal transport plan $\gamma_0\in \C_0(\mu,\nu)$ so that
	\[
		d_{W_b}(\mu,\nu)= \left( \iint |x-y|^b\,d\gamma_0(x,y) \right)^{1/b} .
	\]
Furthermore, if $b \geq1$ and  $\nu$ is absolutely continuous with respect to Lebesgue measure,  
	\[ 
		\nu \in \P_{b,ac}(\Rd) := \left\{ \rho \in \P_b(\Rd) \colon \rho \ll \mathcal{L}^d \right\},
	\]
then there exists an optimal transport plan $\gamma_0$ that is given by the product of the identity map $\id(x) = x$ and a Borel measurable function $\bt_\nu^\mu:\Rd \to \Rd$, i.e.,  $\gamma_0 = (\id \times \bt_\nu^\rho)_\# \nu$  (cf.  \cite[Theorem 6.2.10]{AGS},  \cite[Theorem 7.1]{ambrosio2003existence}). The function $\bt_\nu^\mu$ is an \emph{optimal transport map} from $\nu$ to $\mu$.

Along with these characterizations of optimal transport plans, for all $b \in [1,2]$, $(\P_b(\Rd),d_{W_b})$ is a complete and separable metric space  \cite[Proposition 7.1.5]{AGS}. We now suppose $b>1$. While bounded subsets of $(\P_b(\Rd),d_{W_b})$ are not generally relatively compact in the $b$-Wasserstein metric \cite[Remark 7.1.9]{AGS}, they are relatively compact with respect to $d_{W_{a}}$ for $a< b$.
Likewise, convergence in $d_{W_b}$ can be characterized as follows (c.f \cite[Remark 7.1.11]{AGS}):
\smallskip
\begin{center}
	\begin{tabular}{lcl}
		$d_{W_b}(\mu_n,\mu)\arrow 0$ & $\Longleftrightarrow$ & $\mu_n\arrow\mu$ weak-$*$ in $\P(\Rd)$ and $\int |x|^b\,d\mu_n(x)\arrow\int |x|^b\,d\mu(x)$,
		\\
		\\
		& $\Longleftrightarrow$  & $\int f(x)\,d\mu_n(x) \to \int f(x) \,d\mu(x)$,\\[0.25cm]
		&   & for all $f\in C(\Rd)$ such that $|f(x)| \leq C(1+ |x-x_0|^b)$.

	\end{tabular}
\end{center}
When $b = 2$, we abbreviate $d_{W} = d_{W_2}$.

\subsection{Gradient Flows and their $\Gamma$-convergence}

We now briefly recall the notion of a \emph{curve of maximal slope} in a compete metric space $(\S, d)$, which generalizes the concept of gradient flows outside the Riemannian context. We refer again to the book by Ambrosio, Gigli, and Savar\'e \cite{AGS} for further details. A curve $u(t): (a,b) \to \S$ is \emph{2-absolutely continuous} if there exists $m \in L^2(a,b)$ so that 
	\beqn \label{absctsdef} 
		d(u(t),u(s)) \leq \int_s^t m(r)\,dr \text{ for all } a <s \leq t < b.
	\eeqn
We denote the space of 2-absolutely continuous curves by $AC^2([a,b],\S)$.

For any 2-absolutely continuous curve, the limit
\[ |u^\pr(t)| = \lim_{s \to t} \frac{d(u(s),u(t))}{|s-t|} \]
exists for a.e. $t \in (a,b)$. Furthermore $m(t):= |u^\pr(t)| \in L^2(a,b)$ satisfies \eqref{absctsdef} and for any $m \in L^2(a,b)$ satisfying \eqref{absctsdef}, we have $|u^\pr(t)| \leq m(t)$ for a.e. $t \in (a,b)$.

Given a functional $\F: \S \to (-\infty, + \infty]$ that is \emph{proper}, i.e., $D(\F) = \{ u \in \S : \F(u) < +\infty \} \neq \emptyset$, its upper gradient is a generalization of the modulus of the gradient from Euclidean space. Specifically, $g: \S \to [0,+\infty]$ is a \emph{strong upper gradient} for $\F$ if for every $u \in AC^2([a,b],\S)$ the function $g \circ u$ is measurable and
	\beqn \label{stronguppergradienteqn}
		|\F(u(t))-\F(u(s))| \leq \int_s^t g(u(r))|u^\pr|(r)\,dr \text{ for all } a<s\leq t<b.
	\eeqn
When $\F$ is  convex and lower semicontinuous, one example of a strong upper gradient is given by the metric local slope \cite[Corollary 2.4.10]{AGS} ,
	\begin{align} \label{metricslopedef}
		|\partial \F|(u) := \limsup_{ v \to u} \frac{(\F(u)-\F(v))_+}{d(u,v)}.
	\end{align}

Next, we recall the definition of a curve a maximal slope. A locally 2-absolutely continuous curve $u: (a,b) \to \S$ is a \emph{curve of maximal slope} for $\F$ with respect to the strong upper gradient $g$ if there exists a non-increasing function $\phi$ so that $\phi(t) = \F \circ u(t)$ for a.e. $t \in (a,b)$ and
\beqn \label{curvemaxslopeeqn}
 \phi^\pr(t) \leq -\frac{1}{2} |u^\pr|^2(t) - \frac{1}{2} g^2(u(t)) \text{ for a.e. } t \in (a,b) .
\eeqn

Suppose $\F:\P_2(\Rd) \to \R \cup \{+\infty\}$  is proper, lower semicontinuous, and $D(|\partial \F|) \subseteq \P_{2, ac}(\Rd)$. For any  $\mu \in D(|\partial \F|)$, a map $\xi \in L^2(\mu)$  belongs to the \emph{subdifferential} of $\F$ at $\mu$ if 
	\beqn \label{subdiffdef}
\F(\nu) - \F(\mu) \geq \int \big\la \xi, \bt_\mu^\nu - \id \big\ra \, d\mu + o(d_W(\mu,\nu)) \text{ for all } \nu \xrightarrow{d_W} \mu .
	\eeqn
We denote this by $\xi \in \partial \F(\mu)$. 

\begin{remark}[subdifferential and metric slope] \label{subdiffmetricslope}
For any $\xi \in \partial \F(\mu)$, we have $\|\xi\|_{L^2(\mu)} \geq |\partial \F|(\mu)$. 
\end{remark}

When $\S=\P_2(\Rd)$ is endowed with the 2-Wasserstein metric $d_W$, a locally 2-absolutely continuous curve $\mu:(0,\infty)\to\P_2(\Rd)$ with $|\mu^\pr|\in L^2_\loc(0,\infty)$ is called a \emph{gradient flow} relative to the functional $\F$ if its velocity vector $v(t)$ satisfies
	\begin{align} \label{gradflowdef}
		-v(t) \in \pt \F(\mu(t)) , \quad v(t)\in \text{Tan}_{\mu(t)} \P_2(\Rd) \qquad \text{ for a.e. } t\in(0,\infty).
	\end{align}
The velocity vector field $v$  is associated to $\mu$ through the continuity equation $\pt_t \mu + \nabla \cdot (v \mu)=0$. 
If $\F$ is proper, lower semicontinuous, bounded below, and regular (see Definition \ref{regulardefinition}) and its metric slope $|\partial \F|$ is a strong upper gradient, then $\mu(t)$ is a gradient flow of $\F$ if and only if $\mu(t)$ is a curve of maximal slope for $|\partial \F|$  \cite[Theorem 11.1.3]{AGS}. In particular, if $\mu(t)$ is a gradient flow of $\F$, then $|\mu^\pr|(t) = |\partial \F|(\mu(t))$ for almost every $t$. 

\medskip

With these definitions in hand, we now recall a general result of Serfaty on the $\Gamma$-convergence of gradient flows on a metric space. We state a mild variant of this result, similar to that used in \cite[Theorem 5.6]{CarrilloCraigPatacchini}, which is a direct consequence of Serfaty's original proof. 

\medskip

\begin{theorem}[cf.{\cite[Theorem 2]{Serfaty} }]\label{thm:Serfatygammagradflow}
	Let $\F_n$ and $\F$ be functionals defined on $(\P_2(\Rd), d_W)$ with strong upper gradients $|\partial \F_n|$ and $|\partial \F|$ . Suppose that $\mu_n$ is a curve of maximal slope of $\F_n$ with well-prepared initial data $\mu_n(0)$, i.e., there exists $\mu(0) \in D(\F)$ so that
	\[ \mu_n(0) \wsto \mu(0)  \quad \text{and} \quad \lim_{n\to\infty} \F_n(\mu_n(0)) = \F(\mu(0)). \]
If there exists some $\mu \in AC^2([0,T], \P_2(\Rd))$ so that $\mu_n(t) \wsto \mu(t)$  for all $t \in [0,T]$ and
	 	\begin{enumerate}[label = (\roman*)]
		\item 
		$ \liminf_{n\to\infty} \F_n(\mu_n(t)) \geq \F(\mu(t)) $,
		\item 
		$\liminf_{n\to\infty} \int_0^t |\mu_n^\pr|^2 (s)\, ds \geq \int_0^t |\mu^\pr|^2(s)\, ds$,
		\item 
		$\liminf_{n\to\infty} \int_0^t  |\partial \F_n|^2(\mu_n(s)) \, ds  \geq   \int_0^t |\partial \F|^2(\mu(s)) \,ds $, \label{Serfatyslopes}
	\end{enumerate}
then $\mu$ is a curve of maximal slope of $\F$ and
	\begin{gather*}
\lim_{n\to\infty} \F_n(\mu_n(t)) = \F(\mu(t)) \quad\text{for all }  \ t \in [0, T], \\
|\partial \F_n| (\mu_n(t)) \to |\partial \F|(\mu(t)) \text{ and } |\mu^\pr_n|(t) \to |\mu^\pr|(t) \text{ in } L^2([0,T]).
	\end{gather*}
\end{theorem}

\medskip

For the 2-Wasserstein metric $d_W$ the second criterion in Theorem \ref{thm:Serfatygammagradflow} above holds independent of the choice of the energy functionals, and follows from the properties of the metric only, as the following elementary lemma shows. For lack of a reference, we include a proof.

\blemma[Lower bound on metric derivatives] \label{lem:lower_bd_met_der}
Suppose $\mu_n $ and $\mu \in AC^2([0,T], \P_2(\Rd))$ for all $n \in \mathbb{N}$. If  $\mu_n (t) \wsto \mu(t)$ in $\P_2(\Rd)$ for all $t \in [0,T]$, then 
		\[ \liminf_{n\to\infty} \int_0^s |\mu_n^\pr|^2(t)\, dt \geq \int_0^s |\mu^\pr|^2(t)\, dt. \text{ for all } s \in [0,T] . \]
\elemma

\begin{proof}
We may assume, without loss of generality, that there exists $0\leq C<+\infty$ so that 
	\[ 
		C = \liminf_{n \to +\infty} \int_0^s |\mu_n^\pr|^2(t)\, dt.
	\]
Choose a subsequence $|{\mu}_n^\pr|(t)$ so that $\lim_{n\to+\infty} \int_0^s |{\mu}_n^\pr|^2(t)\, dt = C$. Then $|{\mu}_n^\pr|(t)$ is bounded in $L^2(0,s)$ as a sequence in $n\in\N$, so, up to a further subsequence, it is weakly convergent to some $\nu(t) \in L^2(0,s)$. Consequently,  for any $0 \leq s_0 \leq s_1 \leq s$,
	\[ 
		\lim_{n\to+\infty} \int_{s_0}^{s_1} |{\mu}_{n}^\pr|(t) \, dt= \int_{s_0}^{s_1} \nu(t)\, dt.
	\]
By taking limits in the definition of the metric derivative and using the lower semicontinuity of $d_W$ with respect to weak-* convergence,
	\[
		d_W(\mu_{n}(s_0), \mu_{n}(s_1)) \leq \int_{s_0}^{s_1} |\mu_{n}^\pr|(t)\, dt \quad \text{ yields } \quad d_W(\mu(s_0),\mu(s_1)) \leq \int_{s_0}^{s_1} \nu(t)\, dt. 
	\]
By \cite[Theorem 1.1.2]{AGS}, this implies that $|\mu^\pr|(t) \leq \nu(t)$ for a.e. $t \in (0,s)$. Thus, by the lower semicontinuity of the $L^2(0,s)$-norm with respect to weak convergence,
	\[ 
		\liminf_{n\to +\infty} \int_0^s |\mu_{n}^\pr|^2(t)\, dt  \geq \int_0^s |\nu(t)|^2\, dt \geq \int_0^s |\mu^\pr|^2(t)\, dt,
	\]
and we obtain the result.
\end{proof}

\section{Power-Law Interaction potentials} \label{examplessection}

In the present section, we prove that the interaction potentials described in the introduction satisfy our main hypotheses.
\begin{theorem} \label{powerlawhypothesesthm}
Suppose $K$ is a power-law interaction potential of the form
\begin{align} \label{powerlaweqn} K(x)=  |x|^{p}/p  \quad \text{ or } \quad K(x) = |x|^q/q - |x|^p/p \quad \text{with} \quad 2-d\leq p<q\leq 2,
\end{align}
where  we adopt the convention  $|x|^0/0 = \log(|x|)$. Then for all $m \geq m_0>  \max\{ d/(d+p-1),1 \}$, $K$ satisfies hypotheses   \ref{GFcont}, \ref{GFbounds}--\ref{GFdualsobolev} for all $2-d \leq p < q \leq 2$, and \ref{GFgrowth} for $2-d \leq p < q \leq 2$ and $q>0$.  
\end{theorem}
 
\begin{remark}
Note that nonnegative combinations of the above potentials continue to satisfy the hypotheses, where the constraint on $m$ depends on the most singular part of the potential.
\end{remark}

We start by showing that power-law interaction potentials $K$ satisfy hypothesis \ref{GFcont}.

\begin{proposition} \label{GFcontprop}
If $K(x) = |x|^p/p$ with $-d < p <0$, then  \ref{GFcont} holds for $r = -d/p$. If $K(x) = |x|^p/p$ with $p \geq 0$ or $K(x)=|x|^q/q-|x|^p/p$ with $-d < p < q $, then $K$ satisfies   \ref{GFcont} for all $r \in (1, +\infty)$. 
\end{proposition}

\begin{proof}
By definition, $K$ is even and locally integrable.
Suppose $K(x) = |x|^p/p$ for $p>0$ or $K(x)=|x|^q/q-|x|^p/p$ with $-d < p < q $. Then $K$ is lower semicontinuous and bounded from below; hence, \ref{GFcont} is satisfied with $K_a=K$ and $K_b=0 \in L^r(\Rd)$ for all $r \in [1,+\infty)$.

Now, suppose $K(x) = |x|^0/0 = \log(|x|)$. Let $B = B_1(0)$ and define $K_a = K \chi_{\Rd \setminus B}$ and $K_b = K \chi_B$. Then $K_a$ is continuous and bounded below and $K_b \in L^r(\Rd)$ for all $r \in [1,+\infty)$. 
Finally, suppose $K(x) = |x|^p/p$ with $-d\leq p<0$. Then we have $K \in L^{-d/p, \infty}(\Rd)$, and \ref{GFcont} is satisfied with $K_a=0$ and $K_b=K$. 
\end{proof}

We now show that power-law interaction potentials from Theorem  \ref{powerlawhypothesesthm} satisfy   \ref{GFgrowth}.

\begin{proposition} \label{GFgrowthprop}
Suppose $K(x) = |x|^p/p$ with $-d < p \leq 2$ or $K(x)=|x|^q/q-|x|^p/p$ with $-d < p < q \leq 2$ and $q>0$. Then $K$ satisfies hypothesis  \ref{GFgrowth}.
\end{proposition}

\begin{proof}
Suppose $K(x) = |x|^p/p$ with $p>0$ or $K(x) = |x|^q/q-|x|^p/p$ with $-d < p < q \leq 2$ and $q>0$. Then $\lim_{|x|\to\infty} K(x)=+\infty$, and \ref{GFgrowth}(ii) is satisfied. Similarly, for $K(x) = |x|^p/p$ with $p=0$, $K$ grows to infinity and satisfies hypothesis \ref{GFgrowth}(i). On the other hand, when $K(x) = |x|^p/p$ with $p<0$, $K$ is strictly increasing with $\lim_{|x|\to\infty} K(x)=0$. Moreover, for $B= B_1(0)$, we have that  $K \in L^{a}(\Rd\setminus B)$ for some fixed $1<a<\infty$, $K\in C^1(\Rd\setminus\{0\})$, and $|\nabla K| \leq 1$ for $|x|> 1$. Therefore, $K$ satisfies the hypothesis \ref{GFgrowth}(i).
\end{proof}

\medskip

Next we verify the hypotheses \ref{GFbounds}--\ref{GFdualsobolev}. As these are preserved under finite linear combinations of interaction potentials $K$, it suffices to show them for potentials of the form
\begin{align} \label{simplepowerlaw} K(x) = |x|^p/p, \quad \ 2-d \leq p \leq 2 . 
\end{align}
We begin with a few results concerning potentials of this form and conclude with the proof of Theorem  \ref{powerlawhypothesesthm} at the end of the section.

\begin{lemma} \label{uniformboundlem}
Suppose $K(x) = |x|^p/p$ for $2-d \leq p \leq1$. Then for all $m_0 >  d/(d+p-1)$ and  $\rho\in \P_2(\Rd) \cap L^{m_0}(\Rd)$ there exists $C_{d,p}>0$ so that $\|\grad K*\rho\|_\infty \leq C_{d,p}(1+ \|\rho\|_{{m_0}})$.
\end{lemma}

\begin{proof}
Since ${m_0} > d/(d+p-1)$, its H\"older exponent satisfies ${{m_0}}^\pr < d/(1-p)$ and $\grad K \in L^{{m_0}^\pr}(B) \cap L^\infty( \Rd \setminus B)$, where $B=B_1(0)\subset\Rd$ is the unit ball.  Thus,
\begin{align*} \label{uniformboundK2} \|\grad K*\rho\|_\infty &\leq \|(\grad K \chi_{B})*\rho \|_\infty + \|(\grad K \chi_{\Rd \setminus B})*\rho\|_\infty \\
& \leq  \|\grad K\|_{L^{{m_0}^\pr}(B)} \|\rho \|_{{m_0}} + \| \grad K\|_{L^\infty(\Rd \setminus B)} \|\rho\|_1, \nonumber
\end{align*}
where $\chi_B$ denotes the characteristic function of $B$.
\end{proof}

We now consider the hypothesis \ref{GFbounds}.

\begin{proposition} \label{GFboundsprop}
Suppose $K(x) = |x|^p/p$ for $2-d \leq p \leq 2$. Then $K$ satisfies hypothesis  \ref{GFbounds} for all ${m_0} > \max \{  d/(d+p-1),1 \}$. 
\end{proposition}

\begin{proof}
When $p \leq 1$, we may use the uniform bound from Lemma \ref{uniformboundlem} to conclude that for any $\nu \in \P_2(\Rd)$, we have $\|\grad K*\rho\|_{L^2(\nu)} \leq \|\grad K *\rho\|_\infty \|\nu\|_1^{1/2} \leq C(\|\rho\|_{{m_0}} + 1)  \leq C'(\|\rho\|_{{m}} + 1)$.
Now, assume $1< p \leq 2$. Then there exists $c>0$ so that $|\grad K(x)| \leq c(|x|+1)$ for all $x \in \Rd$. Thus, by Minkowski's integral inequality and Jensen's inequality for the concave function $s \mapsto s^{1/2}$,
\begin{align*}
\|\grad K*\rho\|_{L^2(\nu)} &\leq \int \left( \int |\grad K(x-y)|^2 \, d \rho(y) \right)^{1/2} \, d \nu(x) \\
&\leq c  \left( \iint (2|x|^2 + 2|y|^2 +1) \, d \rho(y)  d \nu(x) \right)^{1/2} \leq \sqrt{2}c ( M_2(\rho)^{1/2} + M_2(\nu)^{1/2} + 1) ,
\end{align*}
which gives the result.
\end{proof}

We next turn to hypothesis \ref{GFcty}.

\begin{proposition} \label{GFctyprop}
Suppose $K(x) = |x|^p/p$ for $2-d \leq p \leq 2$. Then $K$ satisfies hypothesis  \ref{GFcty} for all ${m_0} >  \max\{ d/(d+p-1),1 \}$.
\end{proposition}

\begin{proof}
It suffices to estimate each component of the gradient $\grad K*\rho = [ \partial_i K*\rho ]$ separately. Our approach is classical (cf. \cite[Theorem 2.2]{Mizuta}), extending known results about continuity properties of singular integrals to interaction potentials with at most quadratic growth at infinity.

Let $R = |x-y|$. Then,
\begin{align*}
\left|\partial_i K *\rho(x) - \partial_i K*\rho(y) \right| &= \left| \int \big(\partial_i K(x-z) - \partial_i K(y-z)\big)\, d \rho(z) \right| \\
&= \left| \int_{B_{2R}(x)} + \int_{\Rd \setminus B_{2R}(x)} \right| \leq  \left| \int_{B_{2R}(x)} \right| + \left| \int_{\Rd \setminus B_{2R}(x)} \right| =: { \rm I} + {\rm II}   \end{align*}

We begin by estimating I. If $p=2$, let $\beta = 1$. Otherwise, choose $\beta \in (0,1]$ so that ${m_0} >d/(d+p-1-\beta) >1$. Define $r:=d/(d+p-1-\beta) \in [1, {m_0})$, and let $r^\pr$ be the conjugate index of $r$. Since   $p-1 +(d/r^\pr) = \beta \in (0,1]$,
\begin{align*}
{ \rm I } &\leq \int_{B_{2R}(x)} |\partial_i K(x-z)|\, d \rho(z) + \int_{B_{2R}(x)} |\partial_i K(y-z)|\, d \rho(z) \\
&\leq \int_{B_{2R}(x)} | x-z|^{p-1}\, d \rho(z) + \int_{B_{3R}(y)} | y-z|^{p-1}\, d \rho(z) \\
&\leq \|\rho\|_{r} \left\| \left| x- \ \cdot \ \right|^{p-1} \right\|_{L^{r^\pr}({B_{2R}(x)})} +\|\rho\|_{r}     \left\| \left| y- \ \cdot \ \right|^{p-1} \right\|_{L^{r^\pr}({B_{3R}(y)})} \\
&\leq C \|\rho\|_{r} |x-y|^\beta,
\end{align*}
for $C = C_{d, p,r}$. Now, we estimate II.   For $x_\alpha:= \alpha x + (1-\alpha) y$, $\alpha\in [0,1]$, we have
\begin{align*}
{\rm II} &= \left| \int_{\Rd \setminus B_{2R}(x)}  \int_0^1 \frac{d}{d \alpha} \partial_i K (x_\alpha-z) \, d \alpha d \rho(z) \right| \\
& = \left| \int_{\Rd \setminus B_{2R}(x)}   \int_0^1 \big\la \grad \partial_i K(x_\alpha-z) , y-x \big\ra \, d \alpha  d \rho(z) \right| \\
&\leq C \, |x-y| \max_{\alpha \in [0,1]} \int_{\Rd \setminus B_{2R}(x)}  | x_\alpha -z |^{p-2}\, d\rho(z).
\end{align*}
If $p =2$, the integral is bounded by 1 and $\beta = r = 1$. Thus, ${\rm II} \leq C|x-y|^{p-1 + (d/r^\pr)} = C|x-y|^{\beta}$. If $p < 2$,  then $p-2 +(d/r^\pr) \leq 0$, and we may bound the integral as follows,
\begin{align*}
\int_{\Rd \setminus B_{2R}(x)} | x_\alpha -z|^{p-2} \, d\rho(z) & \leq \|\rho\|_r \left( \int_{\Rd \setminus B_{2R}(x)} | x_\alpha -z |^{(p-2)r^\pr} \, dz \right)^{1/r^\pr} \\
& \leq \|\rho\|_{r} \left( \int_{\Rd \setminus B_{R}(x_\alpha)} |x_\alpha -z|^{(p-2)r^\pr}\, dz \right)^{1/r^\pr} \\
& \leq   C \|\rho\|_r |x-y|^{p-2 +(d/r^\pr)} ,
\end{align*}
Therefore,
\begin{align*}
{\rm I } + {\rm II} \leq C \|\rho\|_r |x-y|^{p-1 + (d/r^\pr)} \leq C (\|\rho\|_{{m}} +1)\psi\big(|x-y|\big),
\end{align*}
for $\psi(s) = s^{p-1 + (d/r^\pr)} = s^{\beta}$, $\beta \in (0,1]$, which completes the proof of \ref{GFcty}.
\end{proof}

In order to show that power-law interaction potentials satisfy property \ref{GFdualsobolev}, we begin with the following estimate, quantifying the stability of $\grad K *\rho$ in the 1-Wasserstein metric when $K$ is sufficiently regular.
\begin{lemma} \label{pgeq1W1}
Suppose $K(x) = |x|^p/p$ for $1 < p \leq 2$. Then there exists $C>0$ so that, for all $\rho, \nu \in \mathcal{P}_1(\Rd)$,
\[ \| \grad K*\rho - \grad K*\nu \|_\infty \leq C (d_{W_1}(\rho,\nu))^{(p-1)} . \]
\end{lemma}
\begin{proof}
Let $\gamma_0 \in \mathcal{C}_0(\rho,\nu)$ be the optimal transport plan from $\rho$ to $\nu$.
Using   $|\grad K(v) - \grad K(w)| \leq C |v-w|^{p-1}$ and the concavity of the right hand side, we obtain
\begin{align*}
\| \grad K*\rho - \grad K*\nu \|_\infty
&= \sup_{x \in \Rd} \left| \iint \left( \grad K(x-y) - \grad K(x-z) \right)  \, d\gamma_0(y,z) \right| \\
&\leq C  \iint |  z - y|^{p-1}  \, d\gamma_0(y,z)  \leq C \left( \int |z - y|   \, d\gamma_0(y,z) \right)^{(p-1)} ,
\end{align*}
which completes the proof.
\end{proof}

We now use the previous estimate to quantify the stability of $\grad K*\rho$ in the $b$-Wasserstein metric for all $b \in [1,2]$ and for general attractive power-law potentials. This generalizes a result of Loeper \cite[Theorem 4.4]{Loeper} to general power-law   potentials, $L^p$ spaces, and $b$-Wasserstein metrics for $b\leq2$. This generalization plays a key role in our proof of the $\Gamma$-convergence of gradient flows, since the 2-Wasserstein gradient flow structure merely provides compactness in $b$-Wasserstein metrics for $b<2$. To obtain convergence of the subdifferentials of the interaction energies, we require continuity of $\grad K*\rho$ with respect to weaker Wasserstein metrics. 

\begin{proposition} \label{LoeperGen}
Suppose $K(x) = |x|^p/p$ for $2-d \leq p \leq 2$ and fix $\beta$ so that
\[ (p-1)_+ < \beta \leq1 \text{ for }p <2 \quad \text{ or } \quad  \beta = 1 \text{ for } p=2 . \]
 Then for all $0\leq \eps <  1$, there exists $C>0$, depending on $d, p, \epsilon,$ and $\beta$, so that
\[ \|\grad K*\rho- \grad K*\nu \|_{p_*} \leq C  \max\big\{ \|\rho\|_{m_*}, \|\nu\|_{m_*} \big\}^{(1-\eps)/(2-\eps)} d_{W_{2-\e}}(\rho,\nu) \]
for all $\rho$, $\nu \in \P_2(\Rd)\cap   L^{m_*}(\Rd)$,
where 
\[m_*=\frac{d}{d+p-1-\beta}  \quad \text{ and } \quad  p_* = \frac{d(2-\eps)}{3-p-\beta -\eps(1-\beta)} . \]
\end{proposition}

\begin{proof}
Let $b = 2-\eps \in (1,2]$, and let  $\bt_\rho^\nu$ be an optimal transport map from $\rho$ to $\nu$ with respect to the $b$-Wasserstein metric.
If $p =2$, then $m_* = 1$, $p_* = +\infty$, and by Lemma \ref{pgeq1W1},
\begin{align*}
\|\grad K*\rho- \grad K*\nu \|_{p_*}  \leq C d_{W_1}(\rho,\nu) \leq C d_{W_b}(\rho,\nu) .
\end{align*}

Now, suppose $p< 2$, so $p_* \in (1,+\infty)$. Let $\rho_\alpha := ((1-\alpha) \bt_\rho^\nu + \alpha \id)_\# \rho$,  $\alpha \in [0,1]$, be the constant speed geodesic in the $b$-Wasserstein metric, and let  $\bt_{\rho_\alpha}^{\rho_1}$ and $\bt_{\rho_\alpha}^{\rho_0}$ the optimal transport maps from intermediate points along the geodesics to the endpoints, $\alpha \in (0,1)$ \cite[Lemma 7.2.1]{AGS}.
 Then, by Minkowski's integral inequality,
\begin{align*}
& \|\grad K*\rho - \grad K*\nu\|_{p_*}   = \left\| \int_0^1 \frac{d}{d \alpha} \grad K* \rho_\alpha \, d \alpha \right\|_{p_*}  \leq \int_0^1 \left\| \frac{d}{d \alpha} \grad K* \rho_\alpha\right\|_{p_*} \, d \alpha \\
& = \int_0^1 \left\| \int D^2 K(\cdot- (1-\alpha) y - \alpha \bt_{\rho}^{\nu}(y) ) (\bt_{\rho}^{\nu}(y) - y) \, d \rho_0(y) \right\|_{p_*} \, d \alpha \\
& = \int_0^1 \left\| \int D^2 K(\cdot- y ) (\bt_{\rho_\alpha}^{\rho_1}(y) - \bt_{\rho_\alpha}^{\rho_0}(y)) \, d \rho_\alpha(y) \right\|_{p_*} \, d \alpha   = \int_0^1 \left\| D^2 K* \big[ (\bt_{\rho_\alpha}^{\rho_1} - \bt_{\rho_\alpha}^{\rho_0})\rho_\alpha \big] \right\|_{p_*} \, d \alpha.
\end{align*}
 By Sobolev's inequality for Riesz potentials ($p>2-d$; cf.  \cite[Section 4.2]{Mizuta}) and the Calder\'on Zygmund inequality ($p = 2-d$; cf. \cite[Theorem V.1]{Stein}) and H\"older's inequality,
\begin{align*} 
\|D^2 K*[(\bt_{\rho_\alpha}^{\rho_1} - \bt_{\rho_\alpha}^{\rho_0})\rho_\alpha] \|_{p_*} &\leq C_{d, p,p_*} \|(\bt_{\rho_\alpha}^{\rho_1} - \bt_{\rho_\alpha}^{\rho_0})\rho_\alpha\|_{l} \\
&\leq   C_{d, p,p_*} \left\| (\bt_{\rho_\alpha}^{\rho_1} - \bt_{\rho_\alpha}^{\rho_0}) {\rho_\alpha}^{1/b} \right\|_b \| {\rho_\alpha}^{(b-1)/b}\|_r \\
&=   C_{d, p,p_*} d_{W_b}(\rho,\nu) \|\rho_\alpha\|_{m_*}^{(b-1)/b}
\end{align*}
for 
\[ l = \frac{(2-\eps)d}{2d+p-1-\beta +\eps(\beta-d-p+1)} \in (1,2] \quad \text{ and } \quad r = \frac{d(2-\eps)}{(d+p-1-\beta)(1-\eps)} \in [2,+\infty]  . \]


Finally, by convexity of $L^p$-norms along $b$-Wasserstein geodesics  \cite[Proposition 9.3.9]{AGS},  
\[ \|\rho_\alpha\|_{m_*} \leq (1-\alpha) \|\rho\|_{m_*} + \alpha \|\nu\|_{m_*} \leq \max\{ \|\rho\|_{m_*}, \|\nu\|_{m_*}\}, \]
 which, combined with the previous inequalities gives the result.
\end{proof}

Finally, we use the result of Proposition \ref{LoeperGen} to show that power-law interaction potentials satisfy hypothesis \ref{GFdualsobolev}. 

\begin{proposition} \label{GFdualsobolevprop}
Suppose $K(x) = |x|^p/p$ for $2-d \leq p \leq 2$. Then $K$ satisfies hypothesis  \ref{GFdualsobolev} for all ${m_0} >  \max\{ d/(d+p-1),1 \}$. In particular, for all $0 \leq \eps < 1$, there exists $C_\eps = C(p,d,m_0, \eps)>0$ and $\alpha_\eps = \alpha (p,d,m_0,\epsilon) \in (0,1]$ so that 
\begin{align*}
\|\grad K*\rho - \grad K *\nu \|_{L^2(\mu)} 
\leq C_\eps  \Big(1+ \|\rho\|_{{m_0}} + \|\nu\|_{{m_0}} + \|\mu\|_{{m_0}} \Big) d_{W_{2-\eps}}^{\alpha_\eps}(\rho,\nu).
\end{align*}
\end{proposition}

\begin{remark}[range of ${m_0}$ and $\eps$]
For ${m_0}$ sufficiently large and $\epsilon = 0$,  we obtain $\alpha = 1$, and this proposition reduces to \cite[Assumption 4.1(v)]{craig2017nonconvex}. Consequently, the main contribution of this new estimate is that it lowers the range of admissible values of ${m_0}$ and strengthens the Wasserstein metric in the estimate. This extension is crucial is the case $p = 2-d$, $d>2$, since previous works required ${m_0} = +\infty$, while  this new estimate allows all ${m_0}  > d$. This extension is also crucial in that it provides stability with respect to weaker Wasserstein metrics when  $\eps >0$. Such a stability result is needed in our proof of $\Gamma$-convergence of the gradient flows, where we merely obtain compactness of the gradient flows in $b$-Wasserstein metrics for $1 \leq b < 2$.
\end{remark}

 \begin{proof}
When $1 \leq p \leq 2$,  Lemma \ref{pgeq1W1} and the fact that $\mu \in \P(\Rd)$ ensures 
\begin{align*}
\|\grad K*\rho - \grad K *\nu \|_{L^2(\mu)}  \leq \| \grad K*(\rho - \nu) \| _\infty \leq d_{W_1}(\rho,\nu)^{p-1} \leq d_{W_{2-\eps}}(\rho,\nu)^{p-1}
\end{align*}
 which gives the result.

Now, suppose that $p<1$. Let
\[ \beta =  \min\big\{d+ p-1-d/{m_0}, 1\big\}  \in ((p-1)_+,1] . \]
 Note that  $\beta = 1$ if and only if ${m_0} \geq \frac{d}{d +p-2}$.
Fix $0 \leq \eps <1$ and define $m_*$ and $p_*$ as in Theorem \ref{LoeperGen}, so that ${m_0} \geq m_*>1$.
Since $2m_*/(m_*-1) \geq p_* \geq 1$, there exists $\alpha_\eps = \alpha(p,d,m_0,\epsilon) \in (0,1]$ so that
\[  \alpha \frac{2m_*}{m_*-1} = p_*. \]
By H\"older's inequality and interpolation of $L^p$-norms, 
\begin{align*} \label{firstprop410}
\|\grad K*\rho - \grad K *\nu \|_{L^2(\mu)} &\leq \|\grad K*(\rho-\nu)\|_{2 m_*/(m_*-1)}\|\mu\|_{m_*}^{1/2} \\
&\leq \|\grad K*(\rho-\nu)\|_{p_*}^{\alpha} \|\grad K*(\rho-\nu)\|_{\infty}^{1-\alpha}\|\mu\|_{m_*}^{1/2} 
\end{align*}
Lemma \ref{uniformboundlem} ensures  $\|\grad K*(\rho-\nu)\|_{\infty} \leq C_{d,p}(1+ \|\rho\|_{m_*}+ \|\nu\|_{m_*}) $. Therefore, applying this and Proposition \ref{LoeperGen} gives
\begin{align*}
&\|\grad K*\rho - \grad K *\nu \|_{L^2(\mu)} \\
&\quad \leq C \|\grad K*(\rho-\nu)\|_{p_*}^{\alpha} \left(1+ \|\rho\|_{m_*} + \|\nu\|_{m_*}  \right)^{1-\alpha}  \|\mu\|_{m_*}^{1/2} \\
&\quad \leq C d_{W_{2-\eps}}(\rho,\nu)^\alpha \max\{ \|\rho\|_{m_*}, \|\nu\|_{m_*} \}^{ \frac{\alpha(1-\eps)}{2-\eps}} (1+ \|\rho\|_{m_*} + \|\nu\|_{m_*})^{1-\alpha} \|\mu\|_{m_*}^{1/2} .
\end{align*}
Simplifying and using that ${m_0} \geq m_*$, we conclude the result for $\psi(s) = s^\alpha$.
\end{proof}

We conclude this section with the proof of Theorem \ref{powerlawhypothesesthm}.

\begin{proof}[Proof of Theorem \ref{powerlawhypothesesthm}]
Hypothesis \ref{GFcont} follows from Proposition \ref{GFcontprop}.
Hypothesis  \ref{GFbounds}--\ref{GFdualsobolev} are preserved under finite linear combinations of interaction potentials $K$, so it suffices to show them for power-law potentials of the form $K(x) = |x|^p/p$ for  $m_0 > \max\{d/(d+p-1), 1\}$. Hypothesis \ref{GFbounds} follows from Proposition \ref{GFboundsprop}. 
To see hypotheses \ref{GFcty}--\ref{GFdualsobolev}, let $\psi_1(s)= s^\beta$ be the modulus of continuity from Proposition \ref{GFctyprop} and let $\psi_2(s) = s^\alpha$ be the modulus of continuity from Proposition \ref{GFdualsobolevprop}. Then $K$ satisfies hypotheses  \ref{GFcty}--\ref{GFdualsobolev} with $\psi(s) :=  \psi_1(s)+ \psi_2(s) $. 
\end{proof}

\begin{remark}(restriction to power-law potentials with $p \geq 2-d$)
As shown in Propositions \ref{GFcontprop} and \ref{GFgrowthprop}, we merely require $p>-d$ for the corresponding energies $\E_m$ to  satisfy hypotheses \ref{GFcont} and \ref{GFgrowth}. (Furthermore, by Proposition \ref{firstLmbound}, hypothesis \ref{GFcont} is sufficient to ensure that $\E_m$ is lower semicontinuous and bounded below; hence its gradient flows exist.) Still, in order to characterize the subdifferentials of the gradient flows of $\E_m$ and study their limits as $m \to +\infty$, we need control over derivatives of $K*\rho$ when $\rho \in L^m(\Rd)$. Consequently, our results on $\Gamma$-convergence of the gradient flows require that our energies $\E_m$ satisfy   hypotheses \ref{GFbounds}--\ref{GFdualsobolev}, which hold merely for $p \geq2-d $.
\end{remark}

\section{Energies, Gradient Flows, and Aggregation-Diffusion Equations} \label{energiessection}
In this section, we develop several fundamental properties of the energies and $\E_m$ and $\E_\infty$, prove existence and (in some cases) uniqueness of their gradient flows. We also rigorously connect these gradient flows  to aggregation-diffusion equations. In the process, we extend the well-posedness theory for such equations and, in some cases, obtain sharper estimates on solutions  than has been previously obtained by pure PDE methods; see Remarks \ref{existenceremark} and \ref{uniquenessremark}.

The key obstacle in analysis of the energies $\E_m$ and $\E_\infty$ is to quantify the competing effects of the interaction, diffusion, and height constraint terms. We denote the diffusion and height constraint parts of the energies by
	\beqn \label{eqn:Renyi}
		\Se_m(\rho) := \begin{dcases*} 
											\frac{1}{m-1}\int \rho^m\,dx &  if $\rho\in L^1(\Rd) \cap L^m(\Rd)$, \\
											+ \infty &  otherwise,
									 \end{dcases*}
		\
		\Se_\infty(\rho) := \begin{dcases*} 
											0 & if $\|\rho\|_{\infty}\leq 1$, \\
											+ \infty &  otherwise,
									 \end{dcases*}
	\eeqn
and the interaction part by 
	\beqn \label{eqn:interac_energy}
		\K(\rho) := \frac{1}{2} \iint K(x-y)\,d\rho(x)d\rho(y)
	\eeqn
so that $\E_m = \Se_m + \K$ and $\E_\infty = \Se_\infty + \K$.

\subsection{Basic Properties of $\E_m$ and $\E_\infty$}
We now develop some basic properties of the energy functionals $\E_m$ and $\E_\infty$. (For ease of notation, we often consider both energies at the same time by proving properties for $\E_m$, allowing $m =+\infty$.) The results in this section merely rely on hypothesis \ref{GFcont}. In the case of attractive power-law interaction potentials, this is equivalent to requiring that we are in the diffusion dominated regime. (See Remark \ref{ddremark}.)

We first show that the energies are bounded below and that an upper bound on $\E_m(\rho)$ implies an upper bound on $\|\rho\|_m$. 

\begin{proposition} \label{firstLmbound}
Suppose $K$ satisfies hypothesis \ref{GFcont} and  $m \in [m_0 ,+\infty]$. Then $\E_m$ is  bounded below, uniformly in $m$, and there exists $ C_{r} >0$ s.t. 
\begin{align}  \|\rho \|_m^{1+\frac{1}{r}} \leq \E_m(\rho) + {C}_{r}.  \label{Lmbound1}
\end{align} 
\end{proposition}

\begin{proof}
First, we show  there exists $C_r >0$ so that, for all $m \in[ m_0, +\infty]$ and  $\rho \in L^m(\Rd) \cap \P(\Rd)$, 
\begin{align} \label{HLS1} \frac12 \int (K*\rho)\,d\rho \geq -C_r(1+\|\rho\|_m^{{1 + \frac{1}{r}}}) . 
\end{align}
By hypothesis \ref{GFcont} and the weak Young inequality (cf. \cite[Section 4.3]{LiLo}, \cite[Lemma 4]{BedrossianRodriguezBertozzi})
\begin{align*}  \frac{1}{2} \int (K*\rho) \, d \rho=  \frac{1}{2} \int (K_a*\rho)\, d \rho + \frac{1}{2} \int  (K_b*\rho)\, d \rho &\geq - \frac12 \|K_a\|_\infty  - \|K_b\|_{L^{r, \infty}} \|\rho \|_{1 + \frac{1}{r}}^{1 + \frac{1}{r}} . 
\end{align*}
Since $m \geq m_0 \geq 1+ \frac{1}{r}$ and $\rho \in \P(\Rd)$, interpolating $L^p(\Rd)$ norms,
\begin{align*}
\|\rho \|_{{1 + \frac{1}{r}}}^{1 + \frac{1}{r}} \leq \|\rho\|_m^{\frac{m}{r(m-1)}} \leq (1+\|\rho\|_m)^{\frac{m}{r(m-1)}} \leq (1+\|\rho\|_m)^{{1 + \frac{1}{r}}} \leq C_{r} \left(1+\|\rho\|_m^{{1 + \frac{1}{r}}} \right).
\end{align*}
Combining the two previous inequalities shows (\ref{HLS1}).

We now show inequality (\ref{Lmbound1}).  For $m =+\infty$, this follows from inequality (\ref{HLS1}) and the definition of $\E_{\infty}$. Suppose $m \in [m_0, +\infty)$.
By definition of $\E_m$ and inequality (\ref{HLS1}), for all  $\rho \in D(\E_m)$,
\begin{align} \label{HLSeqn0}
		 \E_m(\rho) \geq \frac{1}{m-1} \|\rho\|_m^m - C  (1+ \|\rho\|_m^{ 1+\frac{1}{r}}) = \|\rho\|_m^{ 1+\frac{1}{r}} \left( \frac{1}{m-1} \|\rho\|_m^{m-1 -\frac{1}{r}}  -C \right)-C .
\end{align}
  First, assume 
\begin{align} \frac{1}{m-1} \|\rho \|_m^{m-1 -\frac{1}{r}} -C \geq 1. \label{case10} \end{align}
Then, $\|\rho\|_m^{1+\frac{1}{r}} \leq \E_m(\rho)+ C $ and (\ref{Lmbound1}) holds. Alternatively, suppose (\ref{case10}) does not hold. Then, 
\[ \|\rho\|_m  \leq \left((1+ C)(m-1)\right)^{1/(m-1 -\frac{1}{r})} . \]
Since $m \geq m_0 > 1+\frac{1}{r}$, there exists $C_1 = C_1(C,r) >0$ so that $\|\rho \|_m^{1+\frac{1}{r}} \leq  C_1$.

To conclude (\ref{Lmbound1}), it suffices to show that $\E_m$ is bounded below, uniformly in $m \geq m_0$. We will show it is  bounded below on the set $\{ \rho \in D(\E_m): \E_m(\rho) \leq 1 \}$. By the previous inequalities, on this set, we have $\|\rho_m\|_m^{ 1+\frac{1}{r}} \leq 1 + C + C_1$.  Combining this with inequality (\ref{HLSeqn0}), we obtain that the energy is uniformly bounded below.
\end{proof}

We now turn to further properties of the energies. For all $m \in [1,+\infty]$, the energies $\Se_m$ are proper, lower semicontinuous with respect to weak-* convergence, and convex in the 2-Wasserstein metric  (cf. \cite[Proposition 9.3.9]{AGS}, \cite[Proposition 4.5]{craig2017nonconvex}).  Likewise, for all $p \in [1, +\infty)$ the $L^p(\Rd)$ norms
\begin{align} \label{Lpnormdef}
\|\rho\|_p := \begin{cases} \|\rho\|_p & \quad \text{ if }\rho\in L^1(\Rd) \cap L^p(\Rd) , \\ +\infty & \quad \text{ otherwise,} \end{cases}
\end{align}
are also proper, lower semicontinuous with respect to weak-* convergence, and convex in the 2-Wasserstein metric.

We now show that the interaction energy is also lower semicontinuous with respect to  weak-* convergence, on $L^{m_0}(\Rd)$-bounded sets.
\begin{proposition}[lower semicontinuity of interaction energy] \label{lcsinteraction}
Suppose $K$ satisfies  hypothesis \ref{GFcont}. If $\rho_n \wsto \rho$ and $\sup_{n \in \mathbb{N}} \|\rho_n\|_{m_0} < +\infty$, then $\liminf_{n \to +\infty} \K(\rho_n) \geq \K(\rho)$.
\end{proposition}
\begin{proof}
First, note that since $\rho_n \wsto \rho$, we have also have $\rho_n \otimes \rho_n \wsto \rho \otimes \rho$. Since $K_a$ is  lower semicontinuous and bounded below,  by Portmanteau theorem (see \cite[Theorem 1.3.4]{van1996weak}),
\begin{align*}
\liminf_{n \to +\infty} \frac12 \iint K_a(x-y) \, d \rho_n(x) d \rho_n(y) \geq \frac12 \iint K_a(x-y) \, d \rho (x) d \rho (y)  .
\end{align*}

Now, we consider  $K_b$.  For any $k>0$, define $K_b \wedge (-k) = \min \{ K_b, -k \}$. Since $K_b: \Rd \to [-\infty,+\infty]$  is lower-semicontinuous, $K_b \wedge (-k)$ is lower semicontinuous and bounded below for all $k>0$. Furthermore, since $K_b \in L^{r, \infty}(\Rd)$, if we define $S_k :=\{ x: K_b(x) \leq -k\}$, by the weak Young inequality  (cf. \cite[Section 4.3]{LiLo}, \cite[Lemma 4]{BedrossianRodriguezBertozzi}), for any $m_* \in ( 1+ \frac{1}{r},2)$ with $m_* \leq m_0$,
\begin{align*} \left| \int ((K_b \wedge (-k))*\rho_n) \,d\rho_n - \int (K_b*\rho_n) \,d\rho_n \right|  \leq \|K_b \chi_{S_k}\|_{L^{r_*, \infty}(\Rd)} \|\rho_n\|_{m_*}^{m_*} , \quad r_* = (m_*-1)^{-1} \nonumber
\end{align*}
Since we assume  $\sup_{n \in \mathbb{N}} \|\rho_n\|_{m_0} < +\infty$ and $m_* \leq m_0$, the second term is  bounded uniformly in $n \in \mathbb{N}$. Likewise, since $r_* < r$, by definition of $L^{r, \infty}(\Rd)$, 
\begin{align*}
\|K_b \chi_{S_k}\|_{L^{r_*, \infty}(\Rd)}  &= \sup_{ \lambda >0} \lambda | \{  |K_b  \chi_{S_k} | >\lambda \}|^{1/r_*} \leq \sup_{ \lambda \geq k} \lambda | \{   |K_b | >\lambda \}|^{1/r_*}  \\
& =  \sup_{ \lambda \geq k} \lambda^{(r_*-r)/r_*} \Big( \lambda | \{ |K_b  | >\lambda \}|^{1/r} \Big)^{r/r_*}   \leq k^{(r_*-r)/r_*}   \|K_b\|_{L^{r, \infty}(\Rd)}^{r/r_*} \xrightarrow{k \to +\infty} 0.
\end{align*}
Therefore, for all $\epsilon >0$, we may choose $k >0$ so that
\begin{align*}
\liminf_{n \to +\infty} \frac12  \int (K_b*\rho_n) \,d\rho_n  &\geq    \liminf_{n \to +\infty} \frac12  \iint (K_b(x-y) \wedge (-k))\, d \rho_n(x) d \rho_n(y) - \epsilon \\
& \geq  \frac12  \int ((K_b  \wedge (-k))*\rho)\,d \rho - \epsilon \geq \frac12 \int (K_b*\rho) \,d\rho  -\epsilon.
\end{align*}
Since $\epsilon >0$ was arbitrary and  $K = K_a + K_b$, we obtain the result.
\end{proof}

We conclude this section by applying the previous proposition to show that $\E_m$ is proper, lower semicontinuous, and bounded below.

\begin{proposition} \label{LSCproposition} Suppose $K$ satisfies hypothesis \ref{GFcont} and $m \in [m_0, +\infty]$. Then $\E_m$  is proper, lower semicontinuous with respect to weak-* convergence of probability measures, and bounded below.
\end{proposition}
\begin{proof}
The fact that $\E_m$ is bounded below is an immediate consequence of Proposition \ref{firstLmbound} and the fact that $\Se_m \geq 0$. We now show $\E_m$ is proper.
By hypothesis \ref{GFcont}, $K \in L^1_\loc(\Rd)$. Therefore, evaluating the energy $\E_m$ on there characteristic function of a ball $B$ of volume 1, we have 
 \[ \E_m(\chi_B) \leq \frac{1}{m_0-1} +  \frac12 \iint_{B \times B} K(x-y)\, dx dy  < +\infty . \]
Hence, $\E_m$ is proper.  

We conclude by proving that $\E_m$ is lower semicontinuous. Suppose that $\rho_n \wsto \rho$. Without loss of generality, we may assume that $\liminf_{n \to +\infty} \E_m(\rho_n) < +\infty$. Taking a subsequence $\rho_n$ so that the $\liminf$ is attained, we may also assume that $\sup_{n \in \mathbb{N}} \E_m(\rho_n) < +\infty$. Applying Proposition \ref{firstLmbound} and interpolation of $L^p$ norms for $1 < m_0 \leq m$, we obtain that $\sup_{n \in \mathbb{N}} \|\rho_n\|_{m_0} < +\infty$. Therefore, by Proposition \ref{lcsinteraction}, $\liminf_{n \to +\infty} \K(\rho_n) \geq \K(\rho)$. Since \cite[Corollary 3.5]{Mc1} ensures $\liminf_{n \to +\infty} \Se_m(\rho_n) \geq \Se_m(\rho)$, this gives the result. 
\end{proof}

\begin{remark}[sharpness of condition $m \geq m_0$]
The condition $m \geq m_0$ in Proposition \ref{LSCproposition} is sharp to ensure that the energy $\E_m$ is lower semicontinuous with respect to weak-* convergence. In particular, for all $\epsilon >0$, there exists $K(x)$ and $m_0$ satisfying \ref{GFcont}  and $m \in (m_0 - \epsilon, m_0)$ so that the energy $\E_m$ is not lower semicontinuous. For example,  we may take $K(x) = |x|^{p}/p$ for $-d<p<0$, $m_0= 1- (p/d)+ \epsilon/4$, and $m = 1- (p/d)-\epsilon/4$. (We assume, without loss of generality, that $\epsilon>0$ is sufficiently small so that $m>1$.) By Proposition \ref{GFcontprop}, $K$ and $m_0$ satisfy hypothesis \ref{GFcont}. For any $\rho_* \in \P_2(\Rd)$ we may consider its sequence of dilations $\rho_\lambda = \lambda^d \rho_*(\lambda x)$, which converges in the weak-* topology as $\lambda \to +\infty$ to a Dirac mass at the origin $\delta_0$. Along this sequence the energies $\E_m$ satisfy $\lim_{\lambda \to +\infty} \E_m(\rho_\lambda)= -\infty$ and $\E_m(\delta_0) = +\infty$ (see, e.g. \cite[equations (12)-(13)]{carrillo2018aggregation}). Therefore, the energy $\E_m$ is not lower semicontinuous in the weak-* topology. 
\end{remark}

\subsection{Subdifferentials and Gradient Flows} \label{sec:subdiff}
We now  characterize   the minimal elements of the subdifferential  of $\E_m$, $m \in [m_0, +\infty)$, and identify elements belonging to the subdifferential of $\E_\infty$. We defer our full characterization of minimal elements of the subdifferential of $\E_\infty$ to the proof of Theorem \ref{thm:grad_flow_diff} in section \ref{sec:proof2}.
Following these results on the subdifferentials, we prove that   gradient flows of $\E_m$ and $\E_\infty$ exist and provide conditions under which they are unique. Throughout this section, we use hypotheses \ref{GFcont}, \ref{GFbounds}--\ref{GFdualsobolev} on the interaction potential and suppose $m \geq m_0$.

In order to analyze the subdifferentials of $\E_m$ and $\E_\infty$, we begin with the following  lemma, which bounds the variation of the nonlocal interaction energy $\K$  along measures in $L^m(\Rd)$. This lemma extends   \cite[Proposition 4.6]{craig2017nonconvex}, where 
 \ref{GFbounds}, \ref{GFcty}, and  \ref{GFdualsobolev} generalize \cite[Assumption 4.1]{craig2017nonconvex}. We defer its proof to appendix section \ref{powerlawappendix}.

\begin{proposition} \label{Kderivative}
Suppose $K$ satisfies \ref{GFbounds}--\ref{GFdualsobolev} and $\rho_0, \rho_1 \in L^m(\Rd)$ for $m \geq m_0$. Then  
\begin{align*}
\left| \K(\rho_1) - \K(\rho_0) -  \int  \big\la \grad K * \rho_{0} , \bt_{\rho_0}^{\rho_1} - \id \big\ra \, d \rho_0   \right| \leq  f(\rho_0,\rho_1) \, \psi(d_W(\rho_0,\rho_1))  \, d_W(\rho_0,\rho_1) 
\end{align*}
for $f(\rho_0,\rho_1) =  C'(1+ \|\rho_0\|_m +\|\rho_1\|_m )$ and $C' = C'(d_W(\rho_0,\rho_1))>0$ is an increasing function of the distance from $\rho_0$ to $\rho_1$.
\end{proposition}
We now apply the previous proposition to obtain the following generalization of Ambrosio, Gigli, and Savar\'{e}'s characterization of the subdifferential for $\lambda$-convex energies to the energies  $\E_m$ and  $\E_\infty$.

\begin{proposition} \label{subdiffalmostconvex}
Suppose $K$ satisfies  \ref{GFcont},  \ref{GFbounds}--\ref{GFdualsobolev}, $m \in [m_0, +\infty]$, and $\rho  \in D(\E_m)$. Then
$\xi \in L^2(\rho)$ belongs to $\partial \E_m(\rho)$ if and only if
\begin{align} \label{Emsubdiffchar1} \E_m(\nu) -\E_m(\rho) \geq \int \la \xi, \bt_\rho^\nu - \id \ra \, d \rho - f(\rho,\nu)  \psi(d_W(\rho,\nu)) d_W(\rho,\nu)  ,  \ \forall \nu \in D(\E_m) ,\end{align} 
for $f(\rho,\nu) $   as in Proposition \ref{Kderivative}.
\end{proposition}

\begin{proof}
First, suppose $\xi \in L^2(\rho)$ satisfies inequality (\ref{Emsubdiffchar1}). We show that it satisfies the subdifferential inequality (\ref{subdiffdef}). For any sequence $\nu_n \to \rho$, $\nu_n \in D(\E_m)$,  we may assume without loss of generality that $\sup_{n} \E_m(\nu_n) <+\infty$, or else the subdifferential inequality (\ref{subdiffdef}) is satisfied trivially.  
Hence, by Proposition \ref{firstLmbound},  $\sup_n \|\nu_n\|_m < +\infty$.  Therefore,  $f(\rho,\nu_n)$ is uniformly bounded as $\nu_n \to \rho$. Consequently, by the definition of the subdifferential (\ref{subdiffdef}),   $\xi \in \partial \E_m(\rho)$.

Now, suppose $\xi \in L^2(\rho)$ belongs to $\partial \E_m(\rho)$ and $\nu \in D(\E_m)$. We show that inequality (\ref{Emsubdiffchar1})  holds. Let $\rho_\alpha= ((1-\alpha) \id + \alpha \bt_\rho^\nu) \# \rho $ be the Wasserstein geodesic from $\rho$ to $\nu$. Then by definition of the subdifferential,  inequality (\ref{subdiffdef}),
\[ \left. \frac{d}{d\alpha} \E_m(\rho_\alpha) \right|_{\alpha = 0}= \lim_{\alpha \to 0} \frac{\E_m(\rho_\alpha) - \E_m(\rho)}{\alpha} \geq \int \la \xi, \bt_\rho^\nu - \id \ra \, d \rho . \]

By Proposition \ref{Kderivative},
\begin{align*}
 \left. \frac{d}{d\alpha} \K(\rho_\alpha) \right|_{\alpha =0} =  \int  \big\la \grad K * \rho_{0} , \bt_{\rho_0}^{\rho_1} - \id \big\ra \, d \rho_0  &\leq  \K(\nu) - \K(\rho) + f(\rho,\nu) \, \psi(d_W(\rho,\nu))  \, d_W(\rho,\nu) .
 \end{align*}
 Likewise, by the convexity of $\Se_m$ for all $m \in [m_0, +\infty]$,
 \begin{align*}
  \left. \frac{d}{d\alpha} \Se_m(\rho_\alpha) \right|_{\alpha =0}  &\leq  \Se_m(\nu) - \Se_m(\rho)  ,
\end{align*}
Adding the three previous inequalities gives the result.

\end{proof}

Next, we apply the previous proposition to characterize elements belonging to the subdifferential of $\E_m$ for $m_0 \leq m < +\infty$.  Our proof generalizes Ambrosio, Gigli, and Savar\'e's characterization of the subdifferential of aggregation-diffusion energies to the case of nonconvex, singular interaction potentials satisfying hypotheses \ref{GFcont}, \ref{GFbounds}--\ref{GFdualsobolev}   (cf. \cite[Theorem 10.4.13]{AGS}), and we defer its proof to appendix section \ref{powerlawappendix}.
\begin{proposition} \label{subdifftheorem}
Suppose $K$ satisfies hypotheses \ref{GFcont}, \ref{GFbounds}--\ref{GFdualsobolev} and $m \in [m_0, +\infty)$. Then, 
	\beqn  \label{subdiffeqnEm}
 |\partial \E_m|(\rho) < +\infty  
	\iff \begin{cases} \rho^m \in W^{1,1}(\Rd), \\
		 (\grad K*\rho)  + \frac{\grad \rho^m}{\rho}  \in \partial \E_m(\rho), \\
		    |\partial \E_m|(\rho) =  \left\|   (\grad K*\rho)  + \frac{\grad \rho^m}{\rho} \right\|_{L^2(\rho)}.
		\end{cases}  
	\eeqn
In particular, for all $\rho \in D(|\partial \E_m|)$, $ (\grad K*\rho)  + \frac{\grad \rho^m}{\rho} $ is the unique element of the subdifferential of $\E_m$ at $\rho$ with minimal $L^2(\rho)$-norm. 
\end{proposition}

\medskip

\begin{remark}[division by $\rho$]
For simplicity, we commit a small notational abuse in the above expression of the minimal element of the subdifferential:  we divide by $\rho$, even though $\rho$ may not be strictly positive. More precisely,  let $w = (\grad K*\rho)  + \frac{\grad \rho^m}{\rho}$ represent a function satisfying  $ (\grad K*\rho) \rho + \grad \rho^m = w \rho$ almost everywhere. This function is unique $\rho$-almost everywhere.
   \end{remark}

We now identify elements belonging to the subdifferential of $\E_\infty$,  providing an upper bound on the metric slope of $\E_\infty$.  (See Theorem \ref{thm:grad_flow_diff} for the  characterization of minimal elements of $\partial \E_\infty$ along gradient flows.)

\begin{proposition} \label{Einftysubdiff}
If $K$ satisfies hypotheses \ref{GFcont}, \ref{GFbounds}--\ref{GFdualsobolev} and $\rho \in D(\E_\infty)$, 
\[ \grad K*\rho \in \partial \E_\infty(\rho)  \text{ and } \|\grad K*\rho \|_{L^2(\rho)} \geq |\partial \E_\infty|(\rho) .\]
\end{proposition}

\begin{proof}
It suffices to show $\grad K*\rho \in \partial \E_\infty(\rho)$, as the second inequality then follows from Remark \ref{subdiffmetricslope}.
By Proposition \ref{subdiffalmostconvex}, it suffices to show for all $\nu \in D(\E_\infty)$,
\begin{multline} \label{Emsubdiffchar2}
	 \E_\infty(\nu) -\E_\infty(\rho)  \\ \geq \int \big\la \grad K*\rho , \bt_\rho^\nu - \id \big\ra \, d \rho  - f(\rho,\nu) \psi(d_W(\rho,\nu)) d_W(\rho,\nu)  ,  \ \forall \nu \in D(\E_m),
\end{multline} 
By definition of $\E_\infty$ and the fact that $\nu, \rho \in D(\E_\infty)$, this is equivalent to 
\begin{align} \label{Emsubdiffchar3} 
	\K(\nu) -\K(\rho) \geq \int \la \grad K*\rho , \bt_\rho^\nu - \id \ra d \rho  - f(\rho,\nu)  \psi(d_W(\rho,\nu)) d_W(\rho,\nu) ,
\end{align} 
which is an immediate consequence of Proposition \ref{Kderivative}.
\end{proof}

 We   apply the previous results to prove Theorem \ref{GFexistthm}, which ensures that, for any initial data $\rho_0 \in D(\E_m)$, the gradient flow exists. It also provide sufficient conditions for the gradient flow to be unique.

\begin{proof}[Proof of Theorem \ref{GFexistthm}]
Existence follows from Proposition \ref{lcsinteraction}, \cite[Corollary 11.1.8]{AGS}, and \cite[Example 11.1.9]{AGS}.   Uniqueness follows from Proposition \ref{Kderivative}, the convexity of $\Se_m$ for all $m \in (1,+\infty]$, and \cite[Proposition 2.8, Theorem 3.12]{craig2017nonconvex}, with $\omega(s) = \sqrt{s\psi(s)}$.  The correspondence between solutions of the aggregation-diffusion equation and gradient flows of $\E_m$ when $m< +\infty$ follows from the characterization of the minimal element of the subdifferential from Proposition \ref{subdifftheorem} and \cite[Corollary 11.1.8]{AGS}. 
\end{proof}

We conclude by proving that gradient flows of $\E_m$ satisfy an energy dissipation identity.

\begin{proposition} \label{GFpropprop}
Suppose $K$ satisfies hypotheses \ref{GFcont} and \ref{GFbounds}--\ref{GFdualsobolev}, $m \in [ m_0, +\infty]$, and  $\rho_m^{(0)} \in D(\E_m)$. Then if $\rho_m(t)$ is the gradient flow of $\E_m$ with initial data $\rho_m^{(0)}$, for all $T>0$, 
\begin{enumerate}[label = (\roman*)]
\item $ |\rho_m^\pr|(t) =|\partial \E_m|(\rho_m(t))< +\infty$  for a.e. $t \in [0,T]$; \label{slopesequal}
\item $\E_m(\rho_m(T)) + \int_0^T |\partial \E_m|^2(\rho_m(t)) \, dt = \E_m(\rho_m^{(0)})$. \label{EDE}
\end{enumerate}
\end{proposition}

\begin{proof}
The result follow from the fact that  the gradient flow is a curve of maximal slope for the strong upper gradient $|\partial \E_m|$; see Corollary \ref{regular corollary}, \cite[Theorem 11.1.3]{AGS}, and \cite[Remark 1.3.3]{AGS}.
\end{proof}

\section{Convergence of Minimizers}\label{sec:proof1}

In this section, we prove our first main results: up to a sequence, minimizers of the energies $\E_m$ converge to minimizers of the energies $\E_\infty$   and these minimizers have uniformly bounded support. We begin by proving the $\Gamma$-convergence of $\E_m$ to $\E_\infty$.

\begin{theorem}[$\Gamma$-convergence of $\E_m$ to $\E_\infty$] \label{liminfgamma}
Suppose $K$ satisfies \ref{GFcont}. If $\rho_m \wsto \rho$, then 
	\[
		\liminf_{m \to +\infty} \E_m(\rho_m) \geq \E_\infty(\rho).
	\]
Furthermore,   for any $\rho \in \P_2(\Rd)$, we have $ \limsup_{m \to +\infty} \E_m(\rho) \leq \E_\infty(\rho)$.
\end{theorem}
\begin{proof}
 Without loss of generality, we may assume that $\liminf_{n \to +\infty} \E_m(\rho_m) < +\infty$. Taking a subsequence $\rho_m$ so that the $\liminf$ is attained, we may also assume that $\sup_{m} \E_m(\rho_m) < +\infty$. Applying Proposition \ref{firstLmbound} and interpolation of $L^p(\Rd)$ norms for $1< m_0 \leq m$, we obtain that $\sup_{m} \|\rho_m\|_{m_0} < +\infty$. Therefore, by Proposition \ref{lcsinteraction}, $\liminf_{m \to +\infty} \K(\rho_m) \geq \K(\rho)$. In particular,  $\sup_{m} \K(\rho_m) > -\infty$.

It remains to show that $\liminf_{m \to +\infty} \Se_m(\rho_m) \geq \Se_\infty(\rho)$. Note that
\[ \sup_{m} \E_m(\rho_m) < +\infty \ \text{ and } \ \sup_{m} \K(\rho_m) > -\infty \implies \sup_{m} \Se_m(\rho_{m}) < +\infty. \]
Therefore, there exists $C>1$ so that $\|\rho_m\|_m \leq C^{1/m}(m-1)^{1/m} \leq C^{1/m}$ for $m \geq m_0$.  By the interpolation of $L^p(\Rd)$ norms, for any $r \in [1, m]$,
	\[
		\|\rho_{m}\|_r \leq \|\rho_{m}\|_m^{1-\theta} \leq C^{(1-\theta)/m} \leq C^{(1-\theta)/r} \leq C^{1/r} ,
	\]
where $0\leq \theta\leq 1$ satisfies $1/r=\theta + (1-\theta)/m$. Since the $L^r(\Rd)$-norm is lower semicontinuous with respect to weak-* convergence, we obtain
	\[
		\|\rho\|_r \leq \liminf_{m\to\infty} \|\rho_{m}\|_r \leq C^{1/r}.
	\]
Sending $r \to +\infty$ then yields
	\[
		\|\rho\|_\infty \leq 1.
	\]
Therefore $\Se_\infty(\rho) =0$, and since $\Se_m$ is positive, we have
	\[
		\liminf_{m \to +\infty} \Se_m(\rho_m) \geq \Se_\infty(\rho).
	\]

	We now turn to the $\limsup$ inequality. 
Without loss of generality, we may assume $\E_\infty(\rho)< +\infty$, so $\|\rho\|_\infty \leq 1$. Applying H\"{o}lder's inequality gives $| \Se_m(\rho) | \leq \frac{1}{m-1}\| \rho \|^{m-1}_\infty \to 0 = \Se_\infty(\rho)$ as $m\to +\infty$, which gives the result.
\end{proof}

We now prove that minimizers of $\E_m$ and $\E_\infty$ exist and are compactly supported.

\bprop[existence of minimizers in $\P_2(\Rd)$]	\label{prop:exist_min}
Suppose $K$ satisfies \ref{GFcont} and \ref{GFgrowth}. Then   $\E_\infty$ and $\E_m$ admit compactly supported minimizers in $\P_2(\Rd)$ for all $m>\max\{m_0,2\}$.
\eprop

\medskip
\begin{remark}
Although it is possible to prove the existence of minimizers in the regime $1<m<2$, this requires additional assumptions in the hypothesis \ref{GFgrowth} (cf. hypothesis (K6) in \cite{CaHiVoYa16}). Since we are interested in the large $m$ regime we choose to prove the above theorem for $m>\max\{m_0,2\}$.
\end{remark}

\begin{proof}
If $K$ satisfies  \ref{GFgrowth}(i), then the existence of minimizers of $\E_\infty$   follows from the fact that the energy is decreasing under symmetric decreasing rearrangements of $\rho$ (see e.g. \cite[Proposition 3.1]{BuChTo2016}). For $\E_m$, existence of compactly supported minimizers is established in \cite[Theorem 3.1 and Lemma 3.7]{CaHiVoYa16}.

If $K$ satisfies  \ref{GFgrowth}(ii), then the interaction potential $K$ is strictly increasing in every coordinate outside of some fixed set, and the existence of a minimizer in $\P(\Rd)$ for the energy $\E_\infty$ over $\P(\Rd)$ follows simply by \cite[Proposition 4.1]{BuChTo2016}. The minimizer in $\P(\Rd)$ is in fact compactly supported by \cite[Lemma 4.4]{BuChTo2016}, and therefore is in $\P_2(\Rd)$.

The existence of a minimizer in $\P(\Rd)$ for $\E_m$ also follows by using the growth of $K$ given by \ref{GFgrowth}(ii), and by adapting the arguments in \cite[Theorem 3.1]{SiSlTo2014}. This strong coercive behavior of $K$ is sufficient to obtain the existence of a minimizer in $\P(\Rd)$ as the diffusion term is bounded from below (for $m>1$); hence, can be controlled by the growth of $K$. 

In order to conclude that the minimizer is indeed in $\P_2(\Rd)$ we need to show that it is compactly supported in $\Rd$. To this end, note that, if $\rho$ minimizes $\E_m$ in $\P(\Rd)$ then a simple calculation shows that it satisfies the first-order variational inequality
	\[
		(K * \rho)(x) \leq (K * \rho)(x) + \frac{m}{m-1}\rho^{m-1}(x) \leq \lambda
	\]
for some $\lambda\in\R$ and for all $x\in\supp\rho$. Note that
	\[
		(K*\rho)(x) \geq \int_{|y|\leq R} K(x-y)\,d\rho(y) \geq C_R \inf\big\{K(z) \colon |z| > |x|-R\big\}
	\]
where $R>0$ is chosen large enough so that $C_R := \rho(\{y \colon |y|<R\})>0$. Thus $\lim_{|x|\to\infty} (K*\rho)(x) \to \infty$; hence, $\big\{x\in\Rd \colon (K*\rho)(x) \leq \lambda \big\}$ is bounded, and so $\supp\rho$ is compact.
\end{proof}

An important step in the proof of Theorem \ref{thm:gamma_conv} is the compactness of a sequence of admissible measures whose $\E_m$-energy is uniformly bounded. The main idea in proving such a compactness theorem is to utilize Lions' concentration compactness lemma \cite{Lions84} in order to show that any sequence of probability measures with uniformly bounded energy is tight up to translations. The proof of the following lemma follows by arguing as in \cite{Bed2010} and \cite{CraigTopaloglu}.

\blemma[Compactness in $\P(\Rd)$] \label{lem:compactness} Suppose $K$ satisfies \ref{GFcont} and \ref{GFgrowth}. Let $\{\rho_m\}_{m>1} \subset \P(\Rd)$ be a sequence so that $\sup_{m  } \E_m(\rho_m) <+\infty$. Then, up to translations, a subsequence of $\{\rho_m\}_{m>1}$ converges to a measure $\rho \in \P(\Rd)$ with respect to the weak-* topology.
\elemma

\bigskip

Now we turn to our convergence result.

\begin{proof}[Proof of Theorem \ref{thm:gamma_conv}] 
The convergence of minimizers is a classical consequence of the $\Gamma$-convergence result in Theorem \ref{liminfgamma}, when the sequence of energies satisfies a sequential compactness property. Let $\{\rho_m\}_{m>1}\subset\P_2(\Rd)$ be a sequence of minimizers of $\E_m$ over $\P_2(\Rd)$. Then there exists $C>0$ such that $\E_m(\rho_m) \leq C$ for $m>1$ sufficiently large. Hence, by Lemma \ref{lem:compactness}, there exists $\rho\in\P(\Rd)$ such that, up to a subsequence, $\rho_m \wsto \rho$ as $m\to+\infty$ in the weak-* topology of $\P(\Rd)$.

Now let $\nu\in P_2(\Rd)$ be arbitrary. Then
	\[
		\E_\infty(\rho) \leq \liminf_{m\to+\infty} \E_m(\rho_m) \leq \liminf_{m\to+\infty} \E_m(\nu) = \E_\infty(\nu).
	\]
However, since minimizers of $\E_\infty$ are compactly supported we have that
	\[
		\inf_{\P(\Rd)} \E_\infty = \inf_{\P_2(\Rd)} \E_\infty.
	\]
Therefore $\rho$ minimizes $\E_\infty$ over $\P(\Rd)$, and since it is compactly supported, we have $\rho\in \P_2(\Rd)$.
\end{proof}

Next we show that, for particular choices of  $K$, a sequence of minimizers $\{\rho_m\}_{m>1}$ has compact support uniform in $m$. In order to establish this we adapt the arguments by Rein \cite{Re2001} to purely attractive interaction potentials, and follow the method by Frank and Lieb \cite{FrLi2016} to handle repulsive-attractive interactions. This allows us to conclude  the sequence of minimizers converges in the 2-Wasserstein metric, despite fact that the compactness result, Lemma \ref{lem:compactness}, holds only in $\P(\Rd)$. This gives Corollary \ref{cor:P2_conv}.

\begin{proof}[Proof of Theorem \ref{thm:unif_comp_supp}]
In order to prove this theorem for interaction potentials  $K(x)= \frac{1}{p}|x|^{p}$, $-d<p<0$, we proceed similarly to \cite{Re2001} where in the regime $d=3$ and $p=-1$, Rein obtains a bound on the support of minimizers of $\E_m$ independent of the diffusion term. Let 
	\[
		I_M := \inf \left\{ \E_m(\rho) \colon \int \rho\,dx  = M \right\}.
	\]
Using the scaling properties of the energy functional $\E_m$ under the transformations of the form $\rho(x) \mapsto l_1\rho(l_2 x)$ it is easy to see (cf. \cite[Lemma 3.5]{Re2001}) that  $I_M<0$ for all $M>0$ by taking $l_1=l_2^d$, and for $0<\bar{M}\leq M$, we have 
	\beqn \label{eqn:inf_scaling}
		I_{\bar{M}} \geq \big(\bar{M}/M\big)^{(2d+p)/d}I_M
	\eeqn
by taking $l_1=1$ and $l_2=(M/\bar{M})^{1/d}$.
	
Now, let $\rho\in L^1(\Rd)\cap L^m(\Rd)$ be any spherically symmetric, nonnegative function with $\|\rho\|_1=1$, and define $M_R:=\int_{|x|\geq R}\rho\,dx$ for any $R>0$. Then the splitting of the energy
	\[
		\E_m(\rho) = \E_m(\rho\chi_{B_R}) + \E_m(\rho\chi_{B_R^c}) + \iint K(x-y)\chi_{B_R}(x)\chi_{B_R^c}(y)\,d\rho(x)d\rho(y),
	\]
combined with the estimate (which follows due to the spherical symmetry of $\rho$) 
	\[
		\left| \iint K(x-y)\chi_{B_R}(x)\chi_{B_R^c}(y)\,d\rho(x)d\rho(y) \right| \leq C (1-M_R)M_R R^{p}
	\]
implies that 
	\beqns
		\begin{aligned}
			\E_m(\rho) &\geq I_{M_R} + I_{1-M_R} - (1-M_R)M_R R^{p} \\
								&\geq \left((1-M_R)^{(2d+p)/d}+M_R^{(2d+p)/d}\right) I_1 -    (1-M_R)M_R R^{p} \\
								&\geq \left(1-\frac{2d+p}{d}(1-M_R)M_R\right) I_1 - (1-M_R)M_R R^{p}.
		\end{aligned}
	\eeqns
where we have used \eqref{eqn:inf_scaling} in the second line and Taylor's expansion in the third line. Defining
	\[
		R_0 := -\frac{d}{(2d+p)I_1}
	\]
the above estimate becomes
	\beqn  \label{eqn:split_estimate}
		\E_m(\rho) \geq I_1 + \left(\frac{1}{R_0}-R^{p}\right)(1-M_R)M_R.
	\eeqn

Take $R>R_0^{-1/p}$, and assume that for any spherically symmetric minimizing sequence $\{\rho_k\}_{k\in\N}$ for $I_1$ we have, up to a subsequence, that $\lim_{k\to\infty} \int_{|x|\geq R} \rho_k\,dx = M_R >0$. Choose $R_k>R$ such that $M_k = \int_{|x|\geq R_k} \rho_k\,dx = 1/2 \int_{|x|\geq R} \rho_k\,dx$. Then, by \eqref{eqn:split_estimate},
	\[
		\E_m(\rho_k) \geq I_1 + \left(\frac{1}{R_0}-R^{p}_k\right)(1-M_k)M_k \geq I_1 + \left(\frac{1}{R_0}-R^{p}\right)(1-M_k)M_k.
	\]
Sending $k\to\infty$, we get
	\[
		I_1 \geq I_1 + \left(\frac{1}{R_0}-R^{p}\right)\left(1-\frac{M_R}{2}\right)\frac{M_R}{2} > I_1;
	\]
a contradiction. Therefore, the minimizer of $\E_m$, which is the weak-* limit of $\rho_k$, is supported in the ball of radius $R_0^{-1/p}$.

For interaction potentials in the power-law form, given by $K(x) = \frac{1}{q}|x|^q-\frac{1}{p}|x|^p$ with $-d<p<0<q$, we adapt the arguments by Frank and Lieb \cite{FrLi2016} to our case. Let $K^a(x)=\frac{1}{q}|x|^q$ and $K^r(x)=-\frac{1}{p}|x|^p$ denote the attractive and repulsive parts of the interaction potential, and $\K^a$ and $\K^r$ the corresponding interaction energies, respectively, so that
	\[
		\E_m(\rho) = \Se_m(\rho) + \K^a(\rho) + \K^r(\rho).
	\]
	
Let $\rho\in L^1(\Rd)\cap L^m(\Rd)$ be any nonnegative function with $\|\rho\|_1=1$. Then
	\beqns
		\begin{aligned}
			(K * \rho)(x) &\geq \int_{\Rd \setminus B_R(x)} K^a(x-y)\,d\rho(y) \geq \frac{R^q}{q} \left(1-\int_{B_R(x)}\rho\,dy\right) \\
								 & \geq \frac{R^q}{q} \left( 1-\sup_{a\in\Rd} \int_{B_R(a)} \rho\,dy \right).
		\end{aligned}
	\eeqns
Together with the positivity of the diffusion term, this implies that
	\[
		\E_m(\rho) \geq \frac{1}{2} \int (K*\rho)\,d\rho \geq \frac{R^q}{2q}\left( 1-\sup_{a\in\Rd} \int_{B_R(a)} \rho\,dy \right).
	\]
Therefore we get
	\beqn \label{eqn:FL_Lem8}
		\sup_{a\in\Rd} \int_{B_R(a)} \rho\,dy \geq 1-\frac{2q\E_m(\rho)}{R^q}.
	\eeqn
	
Now let $\rho$ be a minimizer of the energy $\E_m$, and let $x\in\Rd$ be a given Lebesgue point of $\rho$ with $\rho(x)>0$. Let $r>0$ be arbitrary, and define
	\[
		\tilde{\rho}(y) := \rho\big(y/l_r\big) \chi_{B^c_r}\big(y/l_r\big) \qquad \text{with} \qquad l_r := \left(\int_{\Rd} \chi_{B_r^c}\,d\rho\right)^{-1/d}
	\]
so that $\|\tilde{\rho}\|_1=1$. We now suppress the dependence on $x$ and denote by $B_r$ the ball of radius $r$ centered at  $x\in\Rd$.   Since $\chi_{B_r}+\chi_{B_r^c}\equiv 1$, we have
	\beqns	
		\begin{aligned}
			\E_m(\tilde{\rho}) &= l_r^{2d+q} \K^a(\rho\chi_{B_r^c}) + l_r^{2d+p} \K^r(\rho\chi_{B_r^c})+\frac{l_r^d}{m-1}\int \rho^m \chi_{B_r^c}\,dy \\
										   &\leq l_r^{2d+q} \left( \K^a(\rho) - \int_{B_r} K^a * \rho \,d\rho(y) + \K^a (\rho\chi_{B_r})  \right) \\
										   &\qquad + l_r^{2d+p} \left( \K^r(\rho) - \int_{B_r} K^r * \rho \,d\rho(y) + \K^r (\rho\chi_{B_r})  \right) +  l_r^d \Big( \Se_m(\rho) - \Se_m(\rho\chi_{B_r}) \Big).
		\end{aligned}
	\eeqns
Since $\E_m(\rho) \leq \E_m(\tilde{\rho})$ and $l_r \geq 1$, we have
	\beqn \label{eqn:energy_in_ball}
		\begin{aligned}
			\int_{B_r} K*\rho \,d\rho(y) &\leq \int_{B_r} K*\rho \,d\rho(y) + \frac{1}{m-1} \int_{B_r} \rho^m \,dy \\
														  & \leq \big(l_r^{2d+q}-1\big)\K^a(\rho) + \big(l_r^{2d+p}-1\big) \K^r(\rho) + \big(l_r^d-1\big) \Se_m(\rho) \\
														  & \qquad\qquad\qquad + l_r^{2d+q} \,\K^a(\rho\chi_{B_r}) + l_r^{2d+p}\, \K^r(\rho\chi_{B_r}).
		\end{aligned}
	\eeqn

First, we will estimate the last two terms and show that they are of order $o(r^d)$. Since $-d<p<0$, we have, by the Hardy-Littlewood-Sobolev inequality,
	\beqn \label{eqn:HLS_repulsive_part}
		\K^r(\rho\chi_{B_r}) \leq C \|\rho\|_{L^1(B_r)} \|\rho\|_{L^{d/(d+p)}(B_r)}.
	\eeqn
Using interpolation of $L^p(\Rd)$ norms, for $m \geq d/(d+p)$,
	\[
		\|\rho\|_{L^{d/(d+p)}(B_r)} \leq \|\rho\|^\theta_{L^1(B_r)} \|\rho\|^{1-\theta}_{L^m(B_r)} \leq C \|\rho\|^\theta_{L^1(B_r)}
	\]
where
	\[
		\theta = \frac{\frac{d+p}{d}-\frac{1}{m}}{1-\frac{1}{m}}.
	\]
On the other hand, by H\"older's inequality, $\|\rho\|_{L^1(B_r)}\leq C \|\rho\|_{L^m(\Rd)}r^{d-d/m}$. Combining these, \eqref{eqn:HLS_repulsive_part} implies
	\beqn \label{eqn:HLS_repulsive_part2}
		\K^r(\rho\chi_{B_r}) \leq C \|\rho\|^{\theta+1}_{L^1(B_r)} \leq C r^{2d(m-1)/m + p}.
	\eeqn 
	
Since the attractive part of the potential is strictly increasing a direct calculation shows that
	\beqn \label{eqn:est_attractive_part}
		\begin{aligned}
			\K^a(\rho\chi_{B_r}) &\leq  \frac{1}{2q} \iint (K^a \chi_{B_{2r}(0)})(y-y^\pr) \,d\rho(y)d\rho(y^\pr) \leq  C r^{d+q},
		\end{aligned}
	\eeqn
where in the last step we use Young's inequality, and the embedding of $L^m(B_{2r}(0))$ into $L^2(B_{2r}(0))$ for $m \geq 2$.

Note that, since $-d<p$, we have $(d-p)/(2d) <1$. Therefore, for $m$ sufficiently large, $(m-1)/m > (d-p)/(2d)$, and $2d(m-1)/m + p - d >0$. By the estimate \eqref{eqn:HLS_repulsive_part2},
	\[
		\frac{1}{|B_r|} \K^r(\rho\chi_{B_r}) \leq C r^{2d(m-1)/m + p - d} \xrightarrow{r\to 0} 0.
	\]
Moreover, the estimate \eqref{eqn:est_attractive_part} implies that
	\[
		\frac{1}{|B_r|} \K^a(\rho\chi_{B_r}) \leq C r^q \xrightarrow{r\to 0} 0.
	\]
Consequently, $\K^a(\rho\chi_{B_r}) + \K^r(\rho\chi_{B_r})=o(r^d)$ as $r\to 0$; and since $l_r\to 1$ as $r\to 0$, the last two terms in \eqref{eqn:energy_in_ball} are of order $o(r^d)$, as well. 

Since $x\in\Rd$ is a Lebesgue point of $\rho$, recalling the definition of $l_r$, we have
	\[
		\frac{l_r^d-1}{|B_r|} = \frac{1}{|B_r|(1-\int_{B_r}\rho\,dy)}  \int_{B_r}\rho\,dy \xrightarrow{r\to 0} \rho(x).
	\]
Also,
	\[
		\frac{l_r^{2d+q}-1}{l_r^d-1}\to \frac{2d+q}{d} \quad \text{ and } \quad \frac{l_r^d-1}{l_r^d-1}\to 1,
	\]
as $r\to 0$. Therefore, dividing both sides of \eqref{eqn:energy_in_ball} by $|B_r|$ and sending $r$ to zero, we get
	\beqn\label{eqn:FL_Lem9}
		K*\rho(x) \leq \frac{2d+q}{d} \E_m(\rho).
	\eeqn
	
Now, let $R=(4q\E_m(\rho))^{1/q}$ so that $1/2=1-2qR^{-q}\E_m(\rho)$. By \eqref{eqn:FL_Lem8}, there exists $a\in\Rd$ such that $\int_{B_R(a)}\rho  \geq 1/2$. 
This implies that for every $y\in \Rd$ such that $|y-a| > (\sigma+1)R$ (with $\sigma$ to be chosen shortly), we have  
	\beqns
		\begin{aligned}
			K*\rho(y) &\geq \int_{B_R(a)} K^a(y-y^\pr)\,d\rho(y^\pr) \geq \frac{(|y-a|-R)^q}{q} \int_{B_R(a)} \rho\,dy^\pr \geq \frac{(|y-a|-R)^q}{2q} \\
							&> \frac{\sigma^q R^q}{2q} = 2\sigma^q \E_m(\rho).
		\end{aligned}
	\eeqns
Let $\sigma:=[(2d+q)/(2d)]^{1/q}$. Now, for $x\in\Rd$ is Lebesgue point of $\rho$ such that $\rho(x)>0$, combining \eqref{eqn:FL_Lem9} with the above estimate yields a contradiction if $|x-a|>(\sigma+1)R$. Therefore $\rho(x)=0$ for all $x\in\Rd$ with $|x-a|>(\sigma+1)R$, 
and consequently,
	\[
		\diam \supp \rho \leq 2(\sigma+1)R = C (\E_m(\rho))^{1/q} \leq C^\pr
	\]
for some constant $C^\pr>0$ independent of $m$ when $m$ is sufficiently large.
\end{proof}

\section{Convergence of Gradient Flows}\label{sec:proof2}

In this section, we prove our main result on the convergence of gradient flows, Theorem \ref{thm:grad_flow_diff}.  Throughout this section, we   impose the following assumptions on our interaction potential $K$, diffusion exponent $m$, and the initial data of the gradient flows $\rho_m^{(0)}$.
\begin{as}[Interaction potential, diffusion exponent, and initial data] \label{GFass}
Suppose $K$ satisfies hypotheses \ref{GFcont}, \ref{GFbounds}--\ref{GFdualsobolev}, $m \geq m_0$, and   $\sup_{m \geq m_0} \big( \E_m(\rho_m^{(0)}) + M_2(\rho_m^{(0)}) \big) < +\infty$.
\end{as}

In the next proposition, we prove that several key quantities remain bounded along the gradient flow, uniformly in $m \geq m_0$. This plays a key role in our proof of $\Gamma$-convergence, since it provides weak compactness of the sequences $\nabla K * \rho_m$ and $ \quad  \nabla \rho_m^m / \rho_m$
with respect to $\rho_m$, as well as weak compactness of $\rho_m$ in arbitrarily large $L^p(\Rd)$ spaces.

\medskip

\begin{proposition}[Uniform bounds along gradient flow] \label{uniformRenyislopes}
Fix $m <+\infty$. Suppose Assumption \ref{GFass} holds and $\rho_m(t)$ is a gradient flow of $\E_m$ with initial data $\rho_m^{(0)}$.
Then  $\rho_m^m(t) \in W^{1,1}(\Rd)$ for a.e. $t  >0$, and for all $T>0$,
\begin{align} \label{GFslopebound}
\sup_{m \geq m_0} \int_0^T |\rho^{\pr}_m|^2(t) +  \| \grad K * \rho_m \|^2_{L^2(\rho_m(t))} + \left\|\grad \rho_m^{m}(t)/\rho_m(t) \right\|^2_{L^2(\rho_m(t))}  \,   dt < +\infty ,
\end{align}
\begin{align}
\sup_{m \geq m_0,\, t \in [0,T]} \|\rho_m(t) \|_m <+\infty , \quad \text{ and } \  \sup_{m \geq \max\{m_0, d/2 \}} \|\rho_m^m(t)\|_{L^2([0,T]\times \Rd)} < +\infty. \label{prelimLm}
\end{align}
\end{proposition}

\begin{proof}
We begin by showing the first inequality in (\ref{prelimLm}). By Propositions \ref{firstLmbound} and  \ref{GFpropprop}, there exists $C_r >0$ so that
\[ \|\rho_m(t)\|_m^{1+\frac{1}{r}} \leq \E_m(\rho_m(t)) +C_r \leq \E_m(\rho_m^{(0)}) +C_r\text{ for all }t >0 . \] 
By Assumption \ref{GFass}, $\sup_{m \geq m_0} \E_m(\rho_m^{(0)}) < +\infty$, which gives the result.

We now consider inequality \eqref{GFslopebound}. By  Proposition \ref{firstLmbound},  $\E_m$ is uniformly bounded below. By Assumption \ref{GFass}, $  \E_m(\rho_m^{(0)})  $ is uniformly bounded above.  Combining this with Proposition \ref{GFpropprop},
\begin{align} \label{metricslopineq} \sup_{m \geq m_0} \int_0^T |\rho_m^{\pr}|^2(t)\, dt = \sup_{m \geq m_0} \int_0^T |\partial \E_m|^2(\rho_m(t)) \, dt < +\infty, \end{align}
which is the first term in \eqref{GFslopebound}. 

Next, we apply inequality (\ref{metricslopineq}) to obtain a uniform bound on the second moment of $\rho_m(t)$. By H\"older's inequality and the definition of the metric slope, 
\[ \sup_{m \geq m_0,\, t \in [0,T]} d_W(\rho_m^{(0)}, \rho_m(t)) \leq \sup_{m \geq m_0}\left(\int_0^T |\rho_m^\pr|^2(t)\, dt \right)^{1/2} \sqrt{T} < +\infty.  \]
Since Assumption \ref{GFass} ensures $M_2(\rho_m^{(0)})$ is uniformly bounded, this implies
\begin{align} \label{secondmomentbd}
\sup_{m \geq m_0, \, t \in [0,T]} M_2(\rho_m(t)) < +\infty . 
\end{align}

We now turn to the second two terms in \eqref{GFslopebound}. By  Proposition \ref{subdifftheorem} and Assumption \ref{GFass}, for almost every $t >0$, we have $\rho_m^m(t) \in W^{1,1}(\Rd)$ and 
\begin{align} \label{subdiffchar4p2}
 |\partial \E_m|(\rho_m(t)) = \| \grad K * \rho_m(t) + \grad \rho_m^m(t)/\rho_m(t) \|_{L^2(\rho_m(t))}. 
\end{align}
By hypothesis \ref{GFbounds}, the uniform bound on $\|\rho_m(t)\|_m$ from inequality (\ref{prelimLm}), and the uniform bound on $M_2(\rho_m(t))$ from inequality (\ref{secondmomentbd}), we have 
\begin{align}  \sup_{m \geq m_0} \int_0^T \| \grad K*\rho_m(t)\|_{L^2(\rho_m(t))}^2 \,  dt  <+\infty  ,\label{firstterm4p2}
\end{align}
which is the second term in  \eqref{GFslopebound}. 
 Then combining equations (\ref{metricslopineq}), (\ref{subdiffchar4p2}),  and (\ref{firstterm4p2}) with the triangle inequality gives 
\begin{align} \label{secondterm4p21}
 \sup_{m \geq m_0} \int_0^T \| \grad \rho_m^{m}(t)/\rho_m(t) \|_{L^2(\rho_m(t))}^2 \,  dt  <+\infty ,
\end{align}
which is the third term in  \eqref{GFslopebound}.  This completes the proof of inequality  \eqref{GFslopebound}.

We finally consider the second term in inequality (\ref{prelimLm}). We proceed by using inequality \eqref{secondterm4p21} and our uniform bound on $\|\rho_m(t)\|_m$  to obtain improved estimates on $\rho_m$. Since $ \rho_m^m(t) \in W^{1,1}(\Rd)$ for almost every $t>0$, for any $\xi \in C^\infty_c([0,T] \times\Rd)$, there exists $C^{\pr}>0$ so that
\begin{align*}
\int_0^T \int \big\la \grad \rho_m^m(t), \xi(t) \big\ra\,dx dt &= \int_0^T \int \big\la \grad \rho_m^m(t)/\rho_m(t), \xi(t) \big\ra\, d \rho_m(t) dt \\
&\leq \int_0^T \|\grad \rho_m^{m}(t)/\rho_m(t)\|_{L^2(\rho_m(t))} \|\xi(t) \|_{L^2(\rho_m(t))} \, dt \\
&\leq \int_0^T \|\grad \rho_m^{m}(t)/\rho_m(t)\|_{L^2(\rho_m(t))} \|\rho_m(t)\|_{m}^{1/2} \|\xi(t)\|_{2m/(m-1)} \, dt \\
&\leq  C^{\pr} \left( \int_0^T \|\xi(t)\|^2_{2m/(m-1)} \, dt \right)^{1/2}.
\end{align*}
Thus, 
\begin{align} \label{lastterm4p21}
 \sup_{m \geq m_0} \int_0^T  \|\grad \rho^m_m(t) \|_{2m/(m+1)}^2 \, dt < +\infty.
\end{align}

Since $ \rho_m^m(t) \in W^{1,1}(\Rd)$, $\rho_m^m(t)$ vanishes at $+\infty$. Thus, because we have  $2m/(m+1) < 2  $, we may apply the Sobolev embedding   \cite[Theorem 11.2]{Leoni}, with 
\[ q= \begin{cases} +\infty &\text{ for } d = 1, \\
2m^2d/(m(d-2)+d) &\text{ for } d \geq 2, \end{cases} \]
which gives
\begin{align} \label{Sobolevemb}
 \|\rho_m(t)\|_{q}^{m} = \|\rho_m^m(t)\|_{q/m} \leq C^{\pr\pr} \|\grad \rho^m_m(t) \|_{2m/(m+1)} . 
 \end{align}
 
 For   $m \geq d/2$, we have $q \geq 2m$.  Interpolating $L^p(\Rd)$-norms gives  
\[ \|\rho_m^m(t)\|_2^2 = \| \rho_m(t)\|_{2m}^{2m} \leq  \|\rho_m(t)\|_{q}^{q^\pr(2m-1)}  , \quad \text{ for } q^\pr = q/(q-1). \]
Integrating in time and applying inequality (\ref{Sobolevemb}),
\begin{align} \label{lastterm4p22}
 \int_0^T \|\rho_m^m(t)\|_{2}^{2} \, dt \leq \int_0^T \|\rho_m(t)\|_{q}^{q^\pr(2m-1)} \, dt    \leq C^{\pr\pr}  \int_0^T \|\grad \rho^m_m(t)\|_{2m/(m+1)}^{2q^\pr/(2m)^\pr} \, dt,
\end{align}
where $(2m)^\pr := 2m/(2m-1)$.
Since $q \geq 2m$, we have $q^\pr \leq (2m)^\pr$, and $s \mapsto s^{q^\pr/(2m)^\pr}$ is concave.  Thus, applying  Jensen's inequality to inequality (\ref{lastterm4p22}) gives
\begin{align*}
  \int_0^T \|\rho_m^m(t)\|_{2}^{2}\, dt & \leq C^{\pr\pr}T \left( \frac{1}{T} \int_0^T \|\grad \rho^m_m(t)\|_{2m/(m+1)}^{2q^\pr/(2m)^\pr} \, dt \right) \\
  &\leq C^{\pr\pr}T \left(\frac{1}{T} \int_0^T  \|\grad \rho^m_m(t)\|_{2m/(m+1)}^2 \, dt \right)^{q^\pr/(2m)^\pr} .
\end{align*}
By inequality (\ref{lastterm4p21}), the right hand side is bounded uniformly in $m$, which gives the result. 
\end{proof}

\bigskip

We conclude with the proof of Theorem \ref{thm:grad_flow_diff}.

\begin{proof}[Proof of Theorem \ref{thm:grad_flow_diff}]
First we show that there exists $\rho(t) \in \P_{2}(\Rd)$ so that, up to a subsequence,  $\rho_{m}(t) \wsto \rho(t) $ for all $t \in [0,T]$. 
By  definition of the metric derivative and Proposition \ref{uniformRenyislopes}, there exists $C^\pr>0$ so that for all $0 \leq s \leq t \leq T$ and $m \geq m_0$,
	\beqn \label{equicontinuity}
			d_W(\rho_{m}(s),\rho_{m}(t)) \leq  \int_s^t |\rho_{m}^\pr|(r)\,dr  \leq \sqrt{t-s} \left( \int_s^t |\rho_{m}^\pr|^2(t)\,dt  \right)^{1/2}  \leq C^\pr \sqrt{t-s} .
	\eeqn
	In particular, taking $s =0$ and recalling that $\sup_{m \geq m_0} M_2(\rho_m(0))<+\infty$, we see that $\{\rho_{m}(t)\}_{m \geq m_0, t \in [0,T]}$ is uniformly bounded in $\P_2(\Rd)$, hence $\rho_m(t)$ is sequentially compact in the $b$-Wasserstein metric for all $1 \leq b < 2$ and $t \in [0,T]$ \cite[Proposition 7.1.5]{AGS}. Take $b = 2-\epsilon$ for $\epsilon \in (0,1)$ as in hypothesis \ref{GFdualsobolev}. Then using the equicontinuity from inequality (\ref{equicontinuity}), the generalized Arzela-Ascoli/Aubin-Lions theorem  \cite[Proposition 3.3.1]{AGS} implies that there exists $\rho(t) \in \P_{2}(\Rd)$ so that, up to a subsequence,  
	\begin{align} \label{Wbconvergence} \lim_{m \to +\infty} d_{W_{2-\eps}}(\rho_{m}(t), \rho(t)) = 0 \text{ for all } t \in [0,T] . 
	\end{align}
	 In particular, up to a subsequence, we have $\rho_{m}(t) \wsto \rho(t)$ for all $t \in [0,T]$.

It remains to verify criteria (i)--(iii) of Theorem \ref{thm:Serfatygammagradflow}  to conclude that $\rho(t)$ is a gradient flow of $\E_\infty$ and that the corresponding energies, local slopes, and metric derivatives converge as $m\to +\infty$. 

Criterion (i) follows immediately from   Theorem \ref{liminfgamma}, and as a consequence, we conclude that $\|\rho(t) \|_\infty \leq 1$ for all $t >0$. Criterion (ii) is proved in Lemma \ref{lem:lower_bd_met_der}.
Thus, it remains to show criterion \ref{Serfatyslopes}. By Fatou's lemma, it suffices to show that 
\begin{align} \label{Fatoureduceiii}
\liminf_{m \to +\infty}  |\partial \E_m|(\rho_m(s)) \geq  |\partial \E_\infty|(\rho(s)) , \quad \text{ for a.e. }  s \in [0,t] . 
\end{align}

By Proposition \ref{uniformRenyislopes} and Fatou's lemma,
\begin{gather} \label{keybound}
\int_0^t \liminf_{m \to +\infty} I_m(s) \, ds  \leq \liminf_{m \to +\infty} \int_0^t I_m(s) \, ds < +\infty  \\
\text{for } I_m(s) :=  |\partial \E_m|^2(\rho_m(s)) + \|v^1_m(s)\|_{L^2(\rho_m(s))}^2 + \|v^2_m(s) \|_{L^2(\rho_m(s))}^2 + \|\rho_m^m(s)\|_2^2 \nonumber \\
v_m^1(s) := (\grad K*\rho_m)(s) \text{ and } v_m^2(s) :=  \grad \rho_m^m(s)/ \rho_m(s). \nonumber
\end{gather}
In particular, $\liminf_{m \to +\infty} I_m(s)<+\infty$ for almost every $s \in [0,T]$. Fix such an $s \in[0,T]$. We will now show that (\ref{Fatoureduceiii}) holds, which completes the proof. As $s \in[0,T]$ is fixed, in what follows, we will suppress the dependence on time to ease notation.

Up to a subsequence, we have
\begin{align*}
\liminf_{m \to +\infty} I_m = \lim_{m \to +\infty} I_m  <+\infty ,
\end{align*}
so that $I_m$ is bounded uniformly in $m$.
 By Proposition \ref{subdifftheorem},  this ensures $\rho_m \in W^{1,1}(\Rd)$ and
\beqn
\begin{gathered} \label{vchareqn}
|\partial \E_m|(\rho_m)  =  \|v_m^1+v_m^2\|_{L^2(\rho_m)}.
\end{gathered}
\eeqn
Since   $\|v^1_m\|_{L^2(\rho_m)}$ and $\|v^2_m\|_{L^2(\rho_m)}$  are   bounded uniformly in $m$, up to another subsequence, there exist $v^1,v^2 \in L^2(\rho)$ so that for all $f \in C^\infty_c(\Rd)$ \cite[Theorem 5.4.4]{AGS},
\beqn
\begin{aligned} \label{weakvsconvergence}
\lim_{m \to +\infty} \int f v^1_m   \, d \rho_m    &=   \int f v^1 \, d \rho    , \quad 
\lim_{m \to +\infty}  \int f v^2_m \, d \rho_m   &=   \int f v^2  \, d \rho   ,
\end{aligned}
\eeqn
 and
\begin{align} \label{futuresubdiffchar}
	\liminf_{m \to +\infty}  |\partial \E_m|(\rho_m ) = \lim_{m \to +\infty}   \|v_m^1 +v_m^2 \|_{L^2(\rho_m )}  \,   \geq  \|v^1 +v^2 \|_{L^2(\rho )}  .
\end{align}
By Remark \ref{subdiffmetricslope}, to complete the proof it suffices to show that $v^1+v^2 \in \partial \E_\infty(\rho)$.

First, we will show that $v^1 = \grad K*\rho$ $\rho$-almost everywhere. For all $f \in C^\infty_c( \Rd)$,
\begin{align*}
 &\left| \int  f  \big[v^1  - (\nabla K*\rho)\big]\, d \rho  \right|  = \lim_{m \to +\infty}  \left| \int  f  \big[(\nabla K * \rho_m)\, d\rho_m   - (\nabla K*\rho) \, d \rho\big] \right| \\
  &\qquad \leq \lim_{m \to +\infty} \left| \int  f  \big[(\nabla K * \rho )\, d \rho_m  - (\nabla K*\rho ) d \rho \big] \right|  +  \left| \int  f  \big[(\nabla K * \rho_m)   - (\nabla K*\rho) \big]\, d \rho_m \right|  \\
  &\qquad = \lim_{m \to +\infty} A_m + B_m .
\end{align*}

By hypothesis \ref{GFcty}, $(\grad K*\rho)(x)$ is continuous in $x$. Thus $\rho_m \wsto \rho$ implies   $\lim_{m \to +\infty} A_m = 0$. 
We now consider $B_m$. By H\"older's inequality and hypothesis \ref{GFdualsobolev}, there exists $\eps >0$ so 
\begin{align*}
B_m \leq  \|f\|_\infty  \|\grad K*(\rho_m - \rho)\|_{L^2( \rho_m)}  \leq C \|f\|_\infty  \psi( d_{W_{2-\eps}}(\rho_m,\rho))  .
\end{align*}
Since $\lim_{m \to +\infty} d_{W_{2-\eps}}(\rho_m,\rho) = 0 $, $\lim_{m \to +\infty} B_m = 0$ and  $v^1 = \grad K*\rho$, $\rho $-almost everywhere.

 By Proposition \ref{Einftysubdiff},   $v^1= \grad K*\rho \in \partial \E_\infty(\rho)$. Thus, by the definition of the subdifferential  \eqref{subdiffdef}, to complete our proof that $v^1 + v^2 \in \partial \E_\infty$, it suffices to show that 
\begin{align} \label{v2subdiff}
\int \big\la v^2, \bt_{\rho}^\nu - \id \big\ra \, d \rho \leq 0 \ , \quad \forall \nu \in D(\E_\infty)  .
\end{align}

We claim that 
\begin{align} \label{mainclaim}
v^2 \rho = \grad \sigma \text{ for } \sigma \in H^1(\Rd) \text{ satisfying }\sigma \geq 0 \text{ and }\sigma =0 \text{ almost everywhere on } \{ \rho < 1\} \ . 
\end{align} 
We now show that this implies (\ref{v2subdiff}).
 Fix $\nu \in D(\E_\infty)$ and let $\bt_{\rho}^\nu$ be the 2-Wasserstein optimal transport map from $\rho$ to $\nu$. Define $\bt_\alpha = (1-\alpha) \id + \alpha \bt_{\rho}^\nu$, $\alpha \in [0,1]$, so $\rho_\alpha = {\bt_\alpha}_{\#} \rho$ is the geodesic from $\rho$ to $\nu$. Since $\rho,\, \nu \in D(\E_\infty)$,  $\|\rho_\alpha\|_\infty \leq \max \{ \|\rho\|_\infty, \|\nu\|_\infty\} \leq 1$ for all $\alpha \in [0,1]$. The geodesic $\rho_\alpha$ satisfies the following weak form of the continuity equation for all $\varphi \in C^\infty_c(\Rd)$ (cf. \cite[equation (8.1.4)]{AGS}, \cite{CarrilloCraigWangWei})
\begin{align} \label{ctygeodesic}
\int \varphi \, d \rho_\alpha - \int \varphi \, d \rho - \int_0^\alpha \! \int \big\la \grad \varphi , \bt_\rho^\nu \circ \bt_\beta^{-1}  - \bt_\beta^{-1} \big\ra \,  d \rho_\beta d \beta = 0.
\end{align}

Note that, for all $\beta \in [0,1]$, we have $\|\rho_\beta\|_2 \leq \|\rho_\beta\|_1^{1/2} \|\rho_\beta\|_\infty^{1/2} \leq 1$ and
\[   \|( \bt_\rho^\nu \circ \bt_\beta^{-1}  - \bt_\beta^{-1} ) \rho_\beta\|_{2} \leq \| \bt_\rho^\nu \circ \bt_\beta^{-1}  - \bt_\beta^{-1} \|_{L^2(\rho_\beta)} \|\rho_\beta\|_\infty^{1/2} \leq W_2(\rho,\nu) .\]
Furthermore, we also have $\lim_{\beta \to 0} (\bt_\rho^\nu \circ \bt_\beta^{-1}  - \bt_\beta^{-1}) \rho_\beta = (\bt_\rho^\nu - \id) \rho$ in distribution. Finally,  since $(\bt_\rho^\nu \circ \bt_\beta^{-1}  - \bt_\beta^{-1}) \rho_\beta$ is uniformly bounded in $L^2(\Rd)$ for all $\beta \in [0,1]$, compactness with respect to the weak $L^2(\Rd)$ topology and uniqueness of limits implies 
\begin{align} \label{momentumconvergence} \text{ as } \beta \to 0 , \quad (\bt_\rho^\nu \circ \bt_\beta^{-1}  - \bt_\beta^{-1}) \rho_\beta \wto (\bt_\rho^\nu - \id) \rho \quad \text{ in } L^2(\Rd) .  
\end{align}

Since $\sigma \in H^1(\Rd)$, approximating it by a sequence in $C^\infty_c(\Rd)$  and applying equation (\ref{ctygeodesic}),
\begin{align*}
\int \sigma \, d \rho_\alpha - \int \sigma \, d \rho - \int_0^\alpha \! \int  \la \grad \sigma , \bt_\rho^\nu \circ \bt_\beta^{-1}  - \bt_\beta^{-1} \ra \, d \rho_\beta d \beta = 0.
\end{align*}
As $\rho = 1$ wherever $\sigma \neq 0$, this is equivalent to
\begin{align*}
\frac{1}{\alpha} \! \int \sigma (\rho_\alpha -1)\,dx =  \frac{1}{\alpha} \int_0^\alpha\! \int \big\la \grad \sigma , \bt_\rho^\nu \circ \bt_\beta^{-1}  - \bt_\beta^{-1} \big\ra\, d \rho_\beta d \beta.
\end{align*}
Since $\sigma \geq 0$ and $\rho_\alpha \leq 1$, this implies that the left hand side is nonpositive for all $\alpha \in (0,1)$. Thus, sending $\alpha \to 0$ and using (\ref{momentumconvergence}) gives
\begin{align*}
0 \geq \limsup_{\alpha \to 0} \frac{1}{\alpha} \int_0^\alpha\! \int \big\la \grad \sigma , \bt_{\rho_\beta}^{\nu} - \bt_{\rho_\beta}^{\rho} \big\ra\, d \rho_\beta d \beta   = \int \big\la \grad \sigma ,\bt_\rho^\nu - \id \big\ra \,d \rho  .
\end{align*}
Since $\rho \geq 0$, the integrand is nonpositive $\rho$-almost everywhere. Since $\rho \leq 1$, we obtain
\begin{align*}
\ird \la  \grad  \sigma , \bt_{\rho}^{\nu} - \id \ra \,dx \leq  \int \big\la \grad \sigma ,\bt_\rho^\nu - \id \big\ra \,d \rho  \leq 0 ,
\end{align*}
which shows (\ref{v2subdiff}).

It remains to show that the claim in equation (\ref{mainclaim}).  Since $I_m$ is bounded uniformly in $m$, so is  $ \|\rho_m^m\|_{2}$, and  there exists $\sigma \in L^2( \Rd)$ with $\sigma \geq 0$   so that, up to a subsequence, $\rho_m^m \wto \sigma$ in $L^2( \Rd)$. Combining this fact with   \eqref{weakvsconvergence}, we have that for all $f \in C^\infty_c(\Rd)$,
\begin{align*}
 -  \int  \grad f \sigma \,dx  &= - \lim_{m \to +\infty} \! \int  \grad f \rho_m^m \,dx   = \lim_{m \to +\infty}  \! \int  f \grad \rho_m^m \,dx   = \lim_{m \to +\infty}  \! \int   f v^2_m \, d \rho_m  =  \! \int  f v^2 \, d\rho .
\end{align*}
Since $\|v^2 \rho\|_{2} \leq \|v^2\|_{L^2( \rho)} \|\rho\|_\infty^{1/2} < +\infty$, we have that $v^2 \rho = \sigma \in H^1(\Rd)$. 

We conclude by showing that $\sigma = 0$ almost everywhere on $\{\rho < 1\}$, which is equivalent to
\begin{align} \label{finalintegral}  \int   \sigma(\rho-1)\,dx  = 0.
\end{align}
By  H\"older's inequality and the uniform bound on $I_m$,
\begin{align*}
 \sup_{m  }   \| \grad \rho_m^{m}\|_{1}   \leq  \sup_{m }  \| \grad \rho_m^{m}(t)/\rho_m \|_{L^2(\rho_m)} \|\rho_m\|_1^{1/2}    =  \sup_{m }  \| v^2_m \|_{L^2(\rho_m)}   <+\infty .
\end{align*}
Fix $R>0$ and let $\eta_R$ be a smooth, radially decreasing cutoff function,
\begin{align*}
\eta_R(x) = \eta(x/R) \text{ for } \eta \in C^\infty(\Rd) \text{ such that } \eta(x) \equiv 1 \text{ for } |x| \leq \frac12 \text{ and } \eta(x) \equiv 0 \text{ for } |x| > 1. 
\end{align*}
Then we have
\begin{align*}
\|\grad (\rho_m^m \eta_R) \|_1 \leq \|\grad \rho_m^m \|_1 + \frac{1}{R} \| \grad \eta \|_2 \| \rho_m^m\|_2 , \quad \| \rho^m_m \eta_R \|_1 \leq \|\rho_m^m\|_2 \|\eta_R\|_2 , \quad \text{ and } \|\rho^m_m \eta_R\|_2 \leq \|\rho_m^m\|_2
\end{align*}
each of which is bounded uniformly in $m$. By Rellich-Kondrachov (cf. \cite[Theorem 13.32]{Leoni}), there exists a subsequence so that $\rho^m_m \eta_R \to \sigma \eta_R$ strongly in $L^1(\Rd)$. Since $\rho^m_m \eta_R$ is uniformly bounded in $L^2(\Rd)$, we also obtain $\rho^m_m \eta_R \to \sigma \eta_R$ strongly in $L^p(\Rd)$ for all $1 \leq p < 2$.

Similarly, by interpolating $L^p(\Rd)$ norms, for any $p^\pr >2$ and $m \geq p^\pr/2$,
\[ \|\rho_m\|_{p^\pr} \leq  1+\|\rho_m\|_{2m} , \]
which is bounded uniformly in $m$. Thus, up to a subsequence, $\rho_m \wto \rho$ weakly in $L^{p^\pr}(\Rd)$. Combining these two facts, we obtain,
\begin{align} \label{almostdone1}
\lim_{m \to +\infty} \int \eta_{R} \rho_m^m(\rho_m - 1) \,dx  \to   \int \eta_{R} \sigma(\rho-1) \,dx 
\end{align}
Rewriting the left hand side of \eqref{almostdone1},
\begin{align*} 
 \int \eta_R \rho_m^m(\rho_m - 1) \,dx &=  \int \eta_R (\rho_m^{m} - \rho_m^{m-1}) \,d\rho_m =  \int \eta_R (\rho_m^{m} - \rho_m^{m(1-1/m)}) \,d\rho_m \\
 													   &= \int \eta_R (\psi(\rho_m^{m},0) - \psi(\rho_m^{m},b))\,d\rho_m,
\end{align*}
where $\psi(s,a):= s^{(1-a)}$ and $b=1/m$. By Lemma \ref{elemineq}, we may control the right hand side by
\begin{align*}
 \frac{1}{m}  \int \eta_R \rho_m \left|1+\rho_m^{2m-1} \right|\,dx   \leq \frac{1}{m}  \int \eta_R \,d\rho_m + \frac{1}{m} \int \eta_R \,\rho_m^{2m} \, dx
\end{align*}
which goes to zero uniformly 
as $m\to\infty$. Thus, for all $R>0$,
\[  \int \eta_R \sigma(\rho-1)\,dx  = 0. \]
Sending $R \to +\infty$ via the monotone convergence theorem, we obtain (\ref{finalintegral}).

This concludes the proof of the criteria from  Theorem \ref{thm:Serfatygammagradflow}. In particular, we have
\[ |\pt \E_m|(\rho_{m}) \to |\pt \E_\infty|(\rho) \text{ in } L^2([0,T]).\]
Integrating  inequality (\ref{futuresubdiffchar}), we obtain that for all $f \in C^\infty([0,T])$ with $f \geq 0$,
\begin{align*}
	 \int_0^T  | \partial \E_\infty(\rho(t))| f(t)\, dt &= \liminf_{m \to +\infty} \int_0^T |\partial \E_m|(\rho_m(t) ) f(t)\, dt \\
	 &\geq  \int_0^T \lim_{m \to +\infty}   \|v_m^1(t) +v_m^2(t) \|_{L^2(\rho_m(t) )} f(t) \,dt \\
	 &  \geq \int_0^T  \|v^1(t) +v^2(t) \|_{L^2(\rho(t) )} f(t) \, dt .
\end{align*}
This gives $ | \partial \E_\infty(\rho(t))|  \geq   \|v^1(t) +v^2(t) \|_{L^2(\rho(t) )}$ for almost every $t \in [0,T]$. On the other hand, by Remark \ref{subdiffmetricslope} and the fact that  $v^1(t) + v^2(t) \in \partial \E_\infty(\rho(t))$, we have the opposite inequality. Therefore, equality holds and
\[ v^1(t) +v^2(t) =  \grad K*\rho(t)+ \grad \sigma(t)/\rho(t) \]
 is the element of $\partial E_\infty(\rho(t))$ with minimal $L^2(\rho(t))$ norm.
\end{proof}

%

\section{  Numerical Results}\label{sec:numerics}
In this section, we apply our theoretical results on the slow diffusion limit to develop a numerical method  to simulate gradient flows and minimizers of the constrained interaction energy $\E_\infty$. For $m \geq m_0$ and $\rho_m^{(0)} \in D(\E_m)$, Theorem \ref{GFexistthm} ensures that $\rho_m(t)$ solves an aggregation-diffusion equation
\begin{align} \partial_t \rho_m  - \grad \cdot ((\grad K*\rho_m)\rho_m) = \Delta \rho_m^m  , \quad \rho_m(0) = \rho_m^{(0)} \label{aggdiffeqn2} 
\end{align}
 if and only if it is a gradient flow of $\E_m$. By Theorem \ref{thm:grad_flow_diff}, gradient flows of $\E_m$ converge, up to a subsequence, to a gradient flow of $\E_\infty$. Thus,  one may approximate the dynamics of gradient flows of $\E_\infty$ by numerically simulating solutions of equation (\ref{aggdiffeqn2}) for $m$ large. By Theorem \ref{thm:gamma_conv}, minimizers of $\E_m$ converge, up to a subsequence and translation, to minimizers of $\E_\infty$. Thus, to numerically approximate minimizers of $\E_\infty$, one may simulate gradient flows of $\E_m$ for $m$ large, in the long time limit.
 
We now give several examples illustrating this approach, computed using Carrillo, Patacchini, and the first author's blob method for diffusion \cite{CarrilloCraigPatacchini}. For present purposes, we merely consider  simulations in one dimension, though our method extends naturally to all dimensions $d \geq 1$. As the primary goal of present work is the rigorous analysis of the slow diffusion limit, we defer a more comprehensive numerical study to future work. Throughout, we take the regularization parameter $\epsilon$ in the blob method for diffusion to depend on the spatial grid spacing $h$ according to $\epsilon = h^{.999}$. The initial data for our simulations is either patch initial data, 
\begin{align} \label{patch initial}
\rho_m^{(0)} = \chi_\Omega \text{ for some } \Omega \subseteq \R^d ,
\end{align}
 or Barenblatt  profiles, for $m_*>1$ and $\tau >0$
\begin{align} \label{Barenblatt}
	\rho_{m}^{(0)}(x) =  
 \tau^{-d\beta}(K - \kappa \tau^{-2\beta} |x|^2)_+^{1/(m_*-1)} ,  \quad \beta = \frac{1}{2+d(m_*-1)} , \quad \kappa = \frac{\beta}{2} \left( \frac{m_*-1}{m_*} \right),
\end{align}
with $K = K(m_*,d) >0$ chosen so that $\int \rho_{m}^{(0)} \,dx= 1$.

\begin{figure}[h]
\begin{center}
 {\bf Equilibria of Aggregation-Diffusion Equation for Varying $m >1$ and mass}  

\vspace{-.3cm}
\hspace{-.7cm}

\hspace{-1.55cm} \rotatebox{90}{\hspace{.5cm} {\bf \footnotesize Mass = 1}} \hspace{.43cm}
\includegraphics[height=3cm,trim={.6cm 10.5cm .6cm 1.2cm},clip]{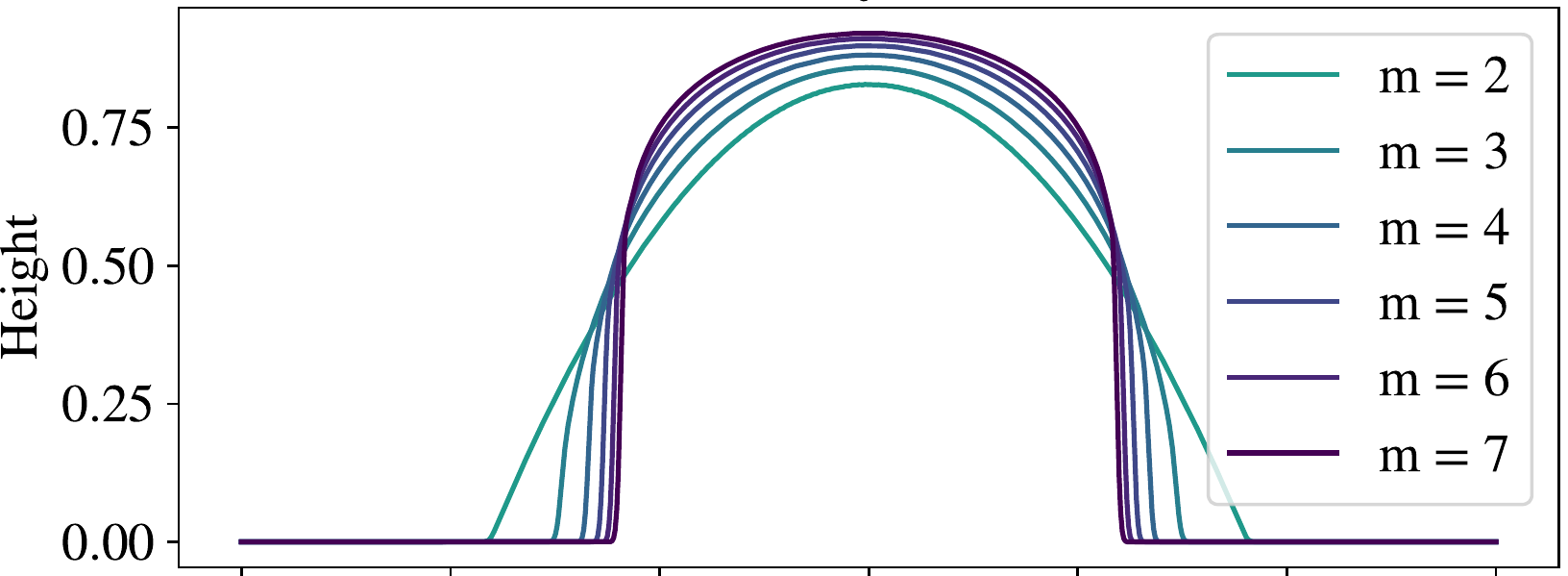}

\hspace{-1.45cm} \rotatebox{90}{\hspace{.5cm} {\bf \footnotesize Mass = 2}} \hspace{.75cm}
\includegraphics[height=2.93cm,trim={1.2cm 10.5cm .6cm 1.25cm},clip]{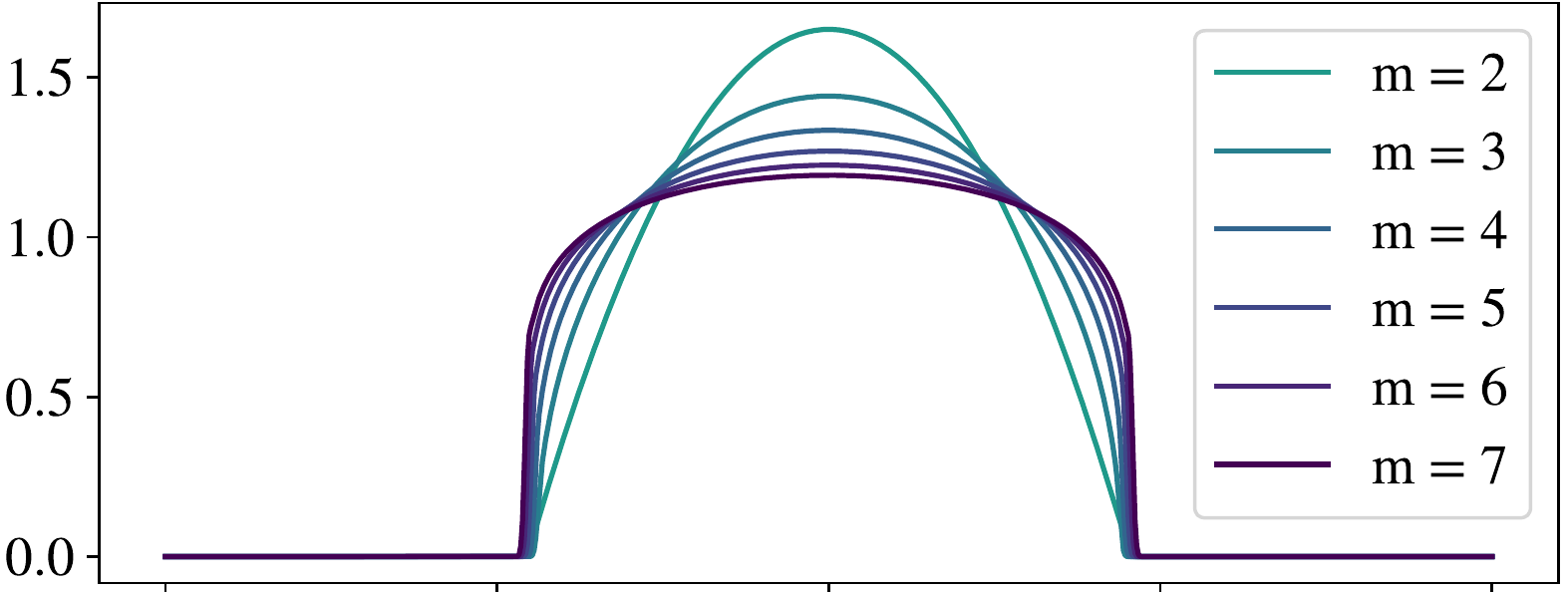}

\hspace{-.1cm} \vspace{-.45cm}

\hspace{-1.4cm} \rotatebox{90}{\hspace{1.1cm} {\bf \footnotesize Mass = 3}} \hspace{.9cm}
\includegraphics[height=3.67cm,trim={1.6cm 9cm .6cm 1.25cm},clip]{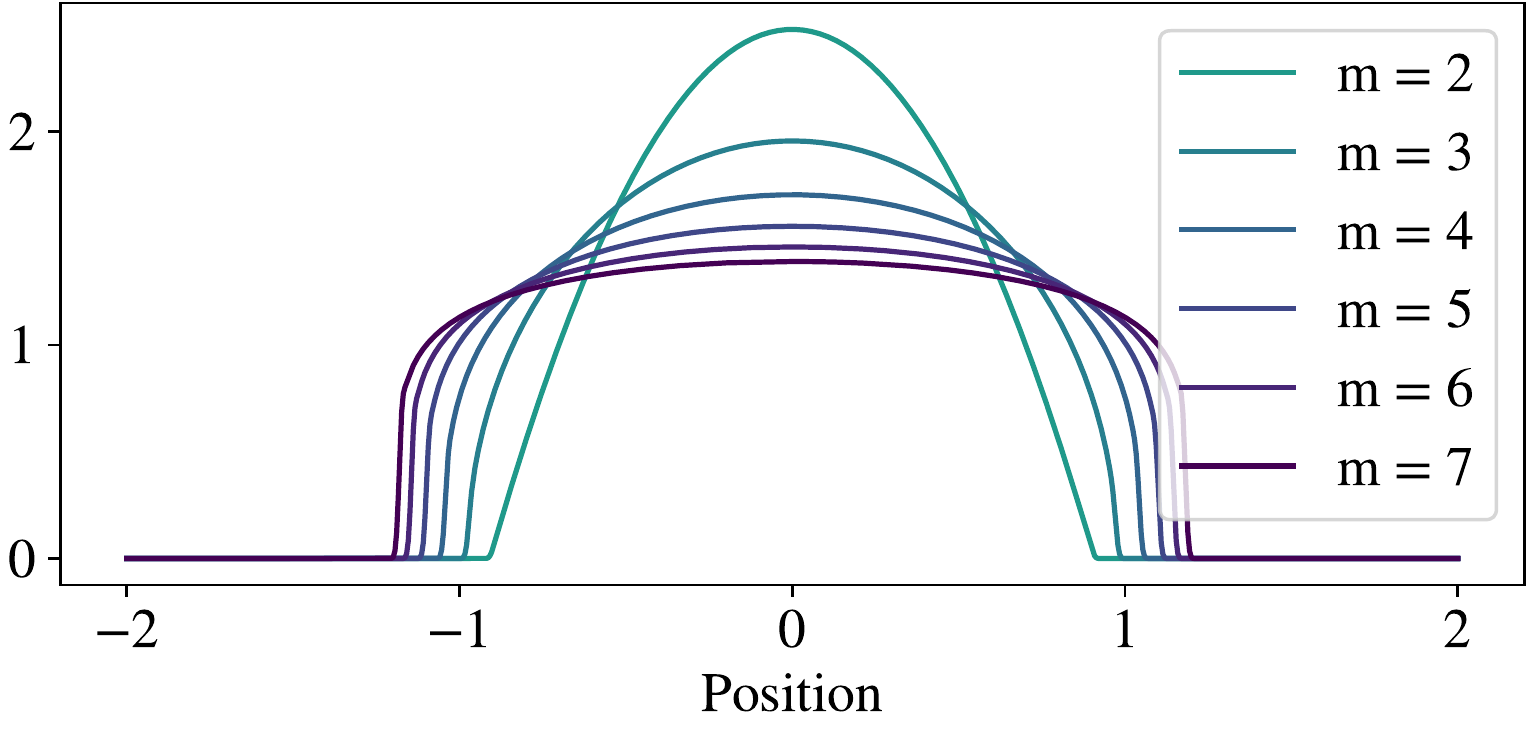}

\vspace{-.2cm}
\caption{Minimizers of $\E_m$, for $K(x) = 2|x|^2$ and varying  $m$ and  $\int \rho_m\,dx$. Solutions of aggregation-diffusion equations 
are simulated at time $T = 10$ with  spatial  and  temporal discretizations  $h = 0.007$, $k = 10^{-4}$. The initial data is a constant multiple of a  Barenblatt profile (equation (\ref{Barenblatt}), $m_*=2$, $\tau = 0.15$).}
\label{supportofminimizers}
\end{center}

\vspace{-.3cm}
\end{figure}

In Figure \ref{supportofminimizers}, we simulate minimizers of   $\E_m$ to study how the  support of minimizers depends on the diffusion exponent $m$ and the mass of the initial data  $\int \rho^{(0)}_m\,dx$. This complements our theoretical result from  Theorem \ref{thm:unif_comp_supp} that the  support of minimizers of $\E_m$ is uniformly bounded for $m$ sufficiently large. For the purely attractive quadratic interaction potential $K(x) = 2|x|^2$, 
the size of the support of minimizers is decreasing in $m$ for $\int \rho^{(0)}_m\,dx =1$, constant for $\int \rho^{(0)}_m\,dx=2$, and increasing in $m$ for $\int \rho^{(0)}_m\,dx = 3$.  This simulation demonstrates that monotonicity properties of the size of the support of minimizers  strongly depend  on the choice of interaction potential $K$ and the mass of the initial data.

\begin{figure}[h]
\begin{center}
 {\bf Critical Mass of Set-Valued Minimizers of $\E_\infty$ for Varying $q$}  
\vspace{.1cm}

 \rotatebox{90}{ \hspace{1.3cm} {\footnotesize  Critical Mass } } \hspace{.1cm}
\includegraphics[height=4.5cm,trim={1.3cm 1.5cm .6cm 1.4cm},clip]{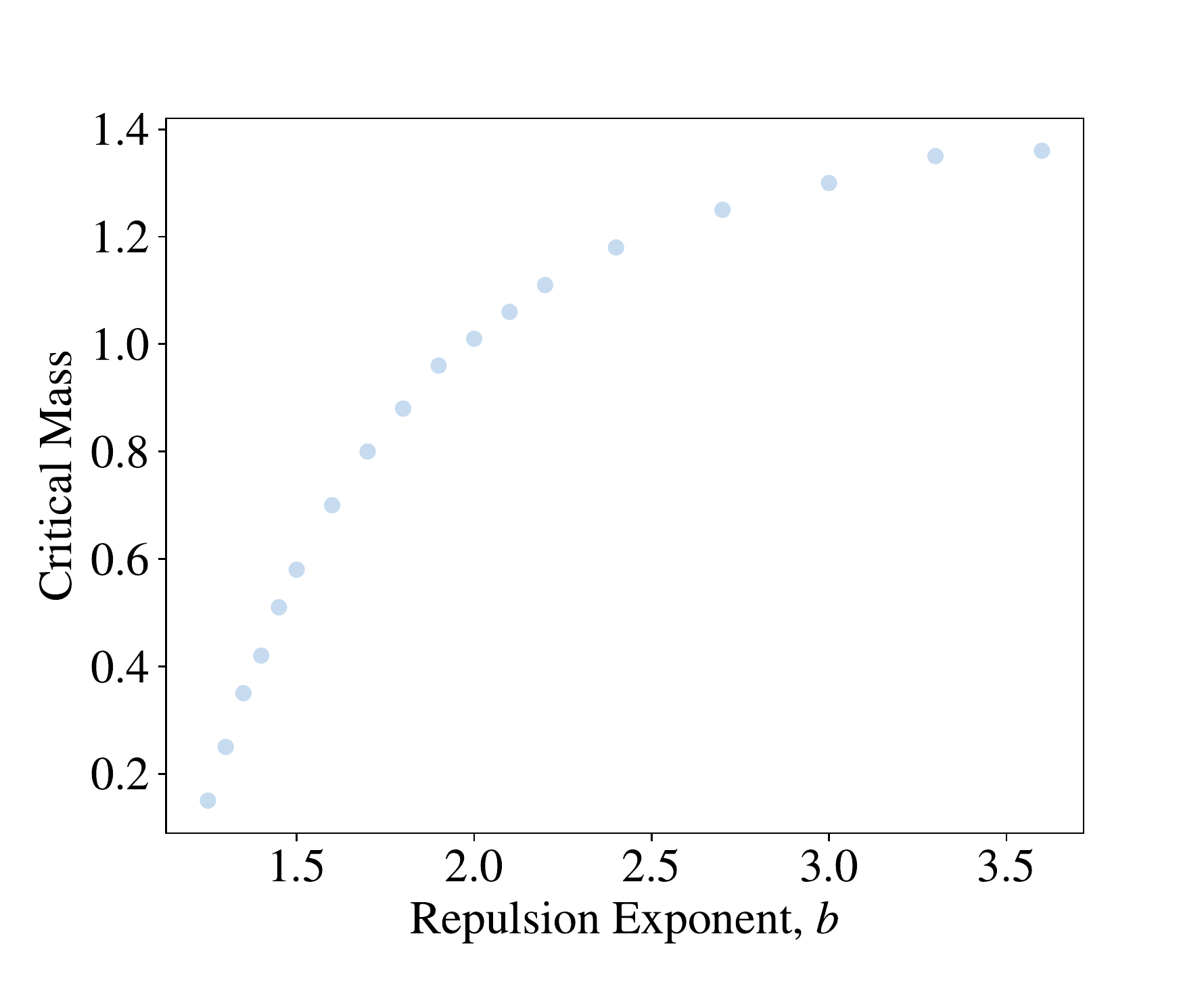}

\hspace{.5cm} {\footnotesize  Attraction Exponent, $q$}

\caption{We approximate the critical mass that determines existence of set-valued minimizers of $\E_\infty$ by simulating solutions of aggregation-diffusion equations for $m = 800$, $K(x) = |x|^q/q - |x|$,  $h = 0.004$, $k = 10^{-4}$, and various choices of initial data. For both initial data given by constant multiples of Barenblatt profiles (equation (\ref{Barenblatt}), $m_* =2$, $\tau = 0.1$) and patch functions (equation (\ref{patch initial}),  $\Omega = [-1,1]$), we observe the existence of set valued minimizers occurs for the same value of mass. Above, we plot how this critical mass value  depends on  $q$.}
\label{critical mass figure}

\vspace{-.7cm}

\end{center}
\end{figure}

Next, we consider gradient flows and minimizers of the constrained interaction energy $\E_\infty$ with repulsive-attractive power-law interaction potentials of the form 
\[ K(x) = |x|^q/q - |x|^p/p , \quad d-2 \leq p \leq q  . \] 
We take the repulsion exponent $p =1 $ (the Newtonian singularity in one dimension) and allow the attraction exponent $q$ to vary.  We apply our numerical method for constrained interaction energies to explore  open questions related to  minimizers of $\E_\infty$, as described in the introduction. 

In Figure \ref{critical mass figure}, we investigate the value of the  critical mass that determines existence versus nonexistence of set valued  minimizers of $\E_\infty$. In particular, for initial data that is either a constant multiple of a Barenblatt profile or patch function, we observe that there exists a single value of the critical mass that separates existence and non-existence of set valued minimizers: in the notation from the introduction, we find $M_1 = M_2$ for all $q \geq 1$.  We plot how the value of the critical mass depends on the attraction exponent $q \geq 1$.

\begin{figure}[h]
\begin{center}
 {\bf Constrained Aggregation: Approaching Equilibrium}  \\ 
\vspace{.1cm}

Barenblatt Initial Data \hspace{2cm} Patch Initial Data

\hspace{-1cm} \rotatebox{90}{\hspace{1.3cm} {\bf \footnotesize $q=1.6$}} \hspace{.25cm} \includegraphics[height=4cm,trim={.6cm 2cm .6cm .6cm},clip]{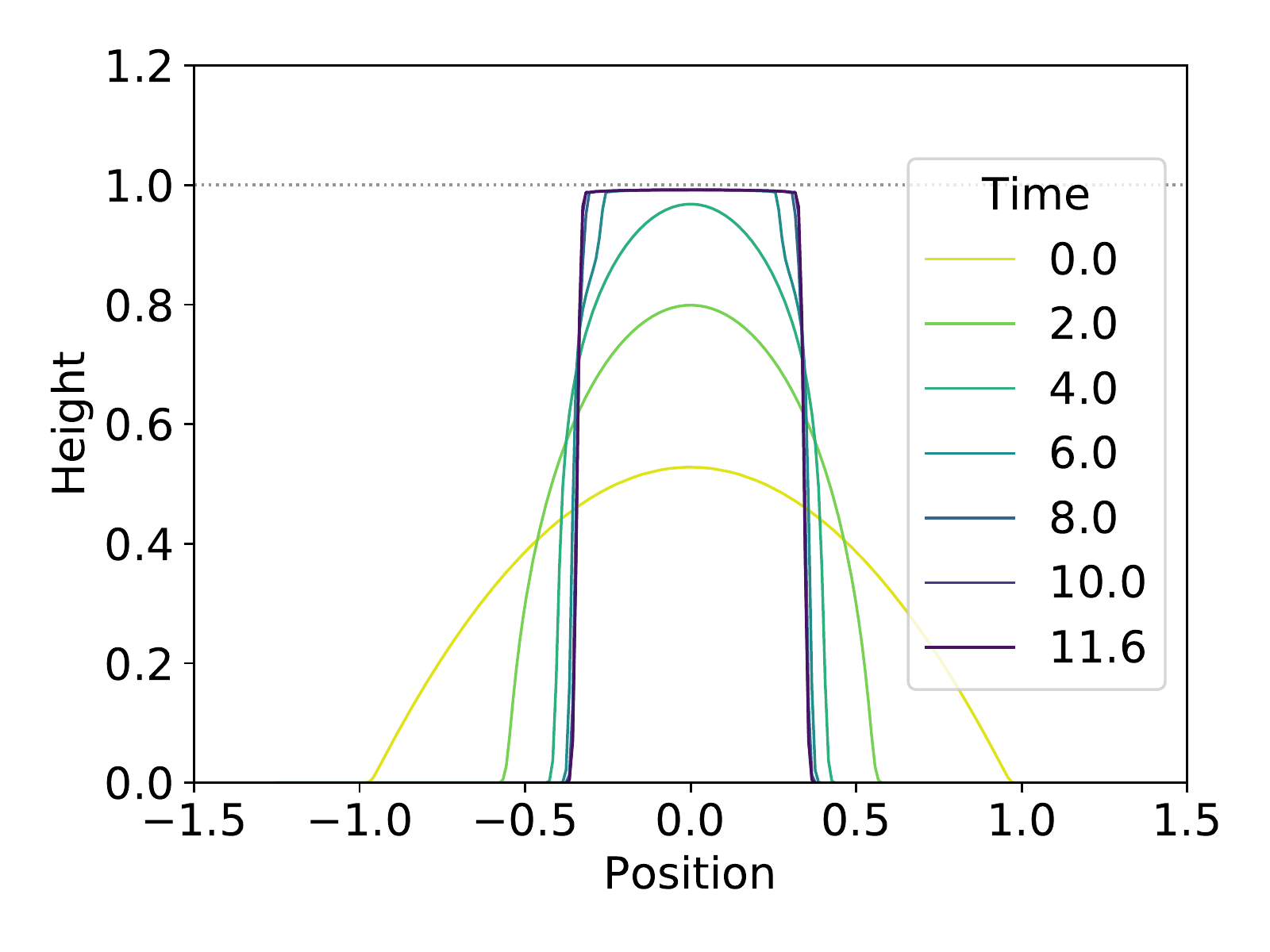}
\includegraphics[height=4cm,trim={2.2cm 2cm .6cm .6cm},clip]{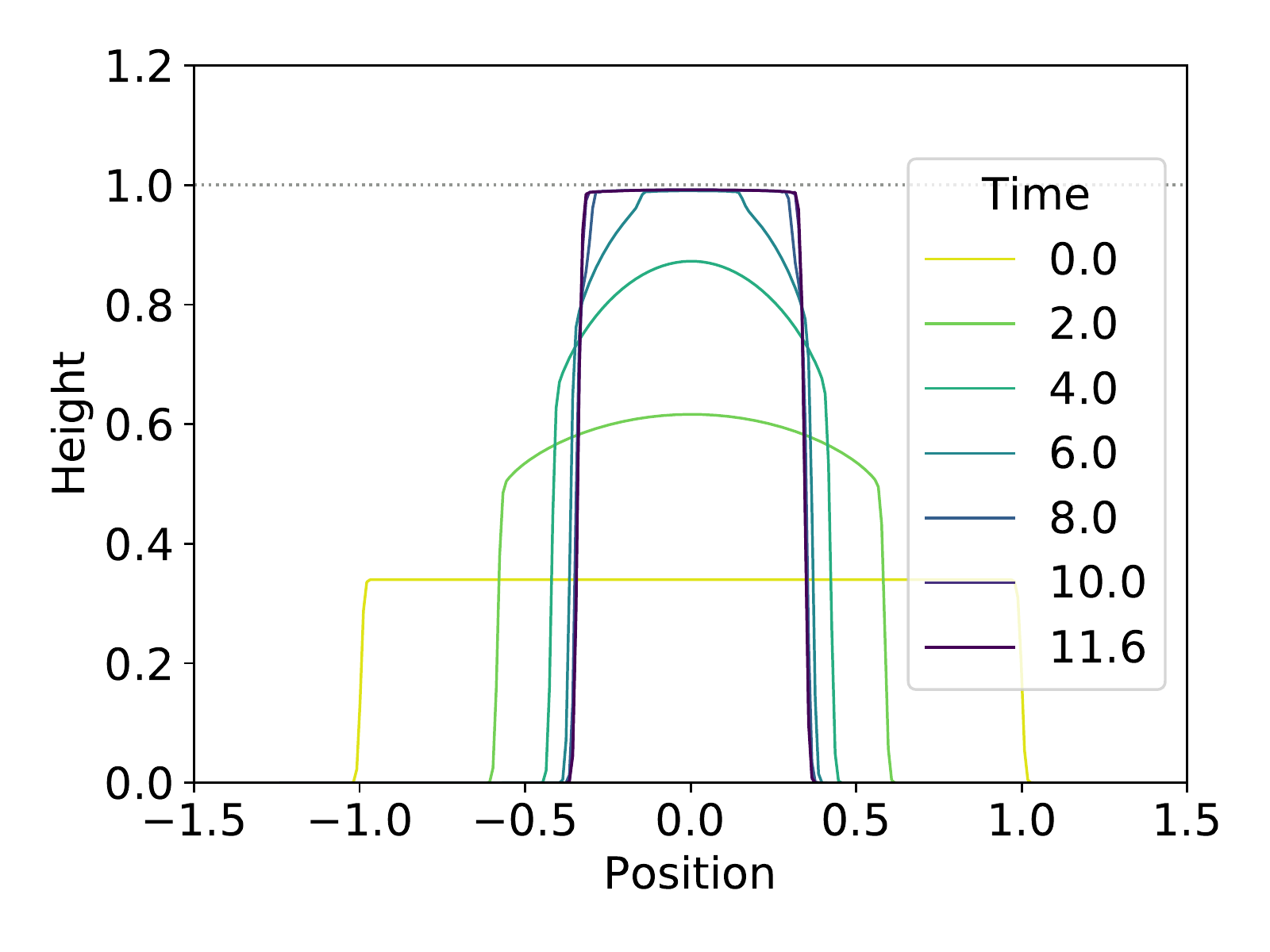}

\hspace{-1cm} \rotatebox{90}{\hspace{1.3cm} {\bf \footnotesize $q=2$}} \hspace{.25cm}
\hspace{.1cm} \includegraphics[height=4cm,trim={1.19cm 2cm .6cm .6cm},clip]{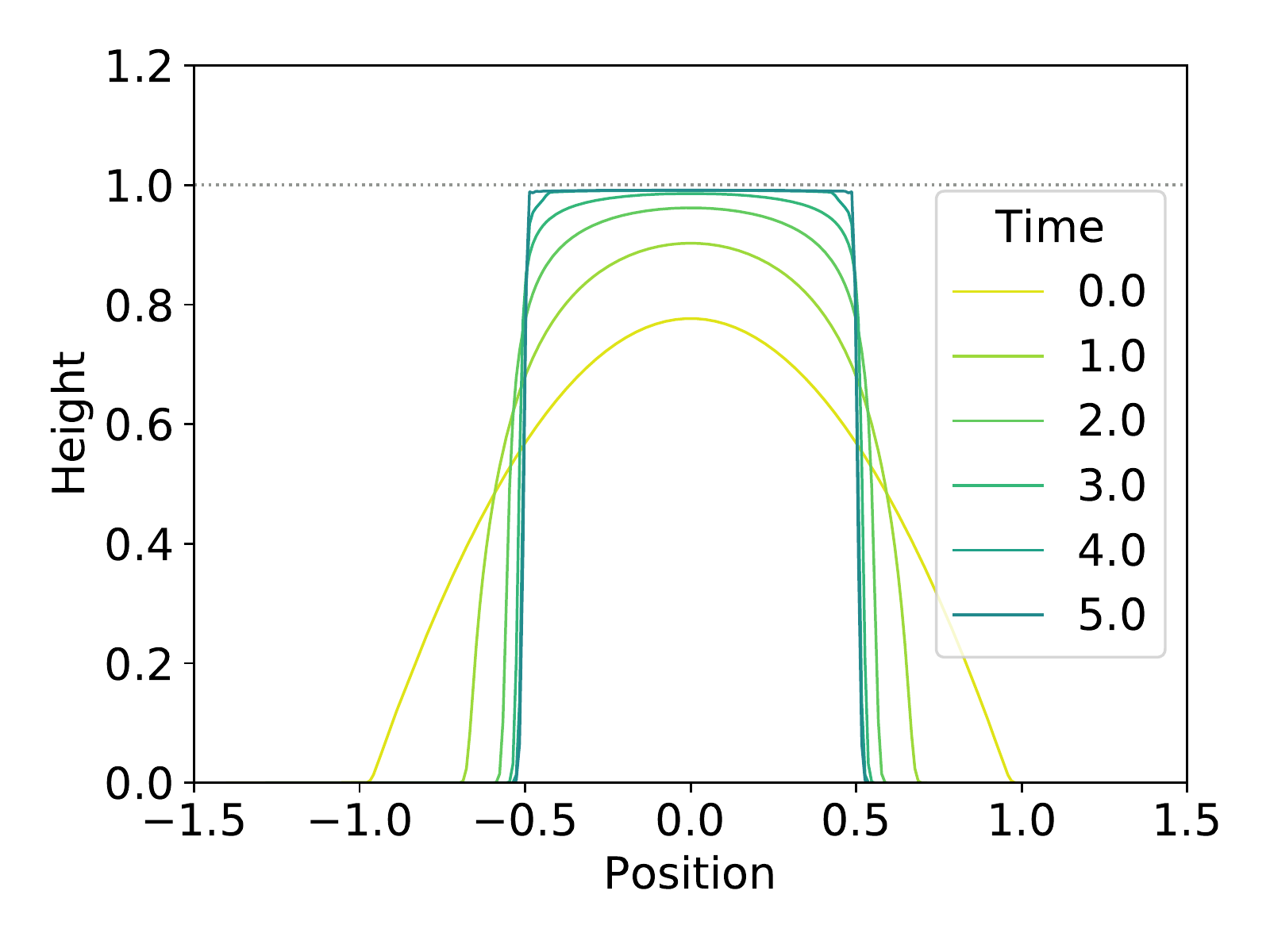}
\includegraphics[height=4cm,trim={2.2cm 2cm .6cm .6cm},clip]{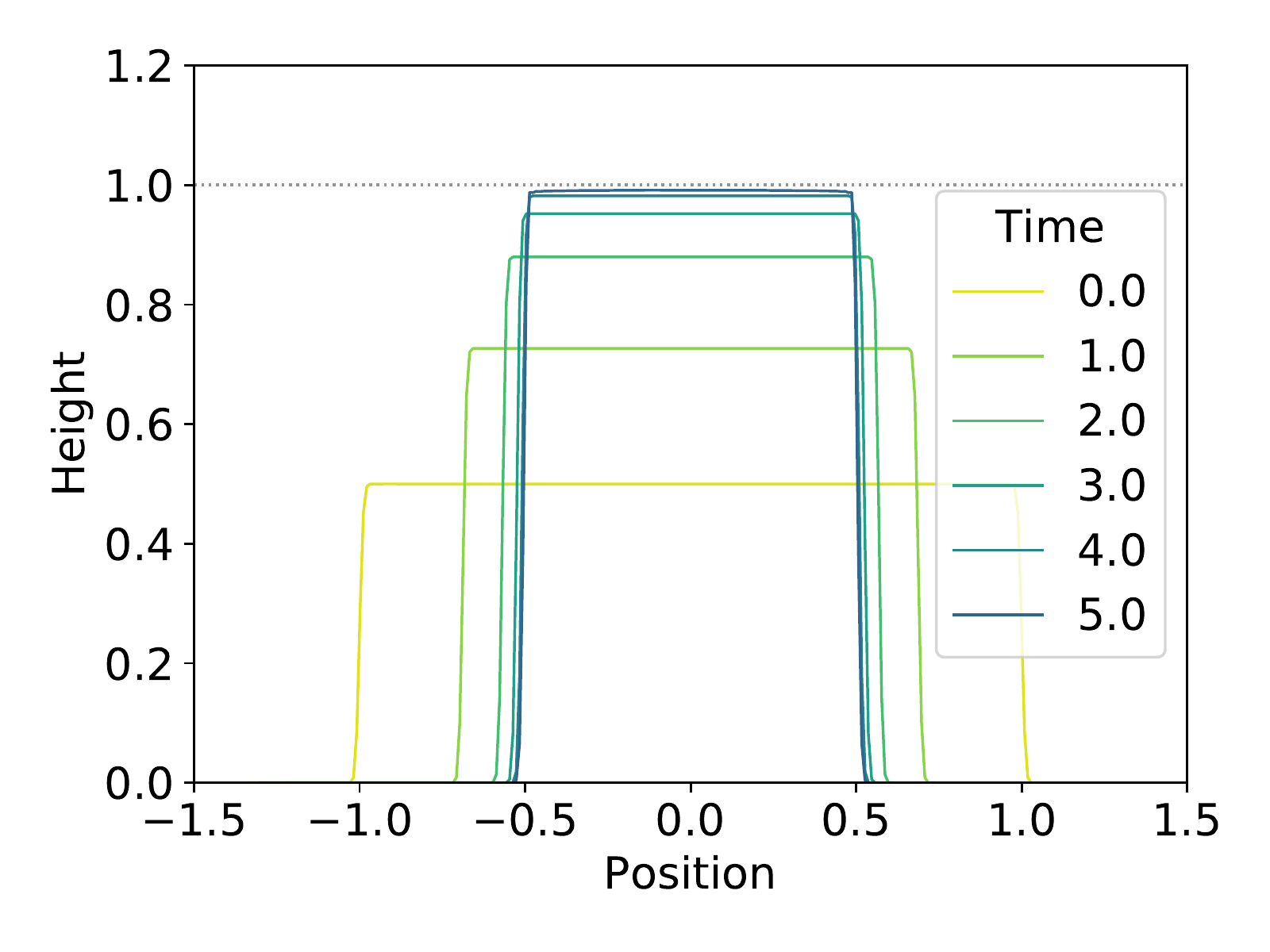}

\adjustbox{valign=t}{\hspace{-1cm} \rotatebox{90}{ {\bf \footnotesize $q=2.4$ \qquad \ \ }}} \hspace{.4cm}
 \includegraphics[height=4.6cm,trim={1.19cm .6cm .6cm .6cm},clip, valign = t]{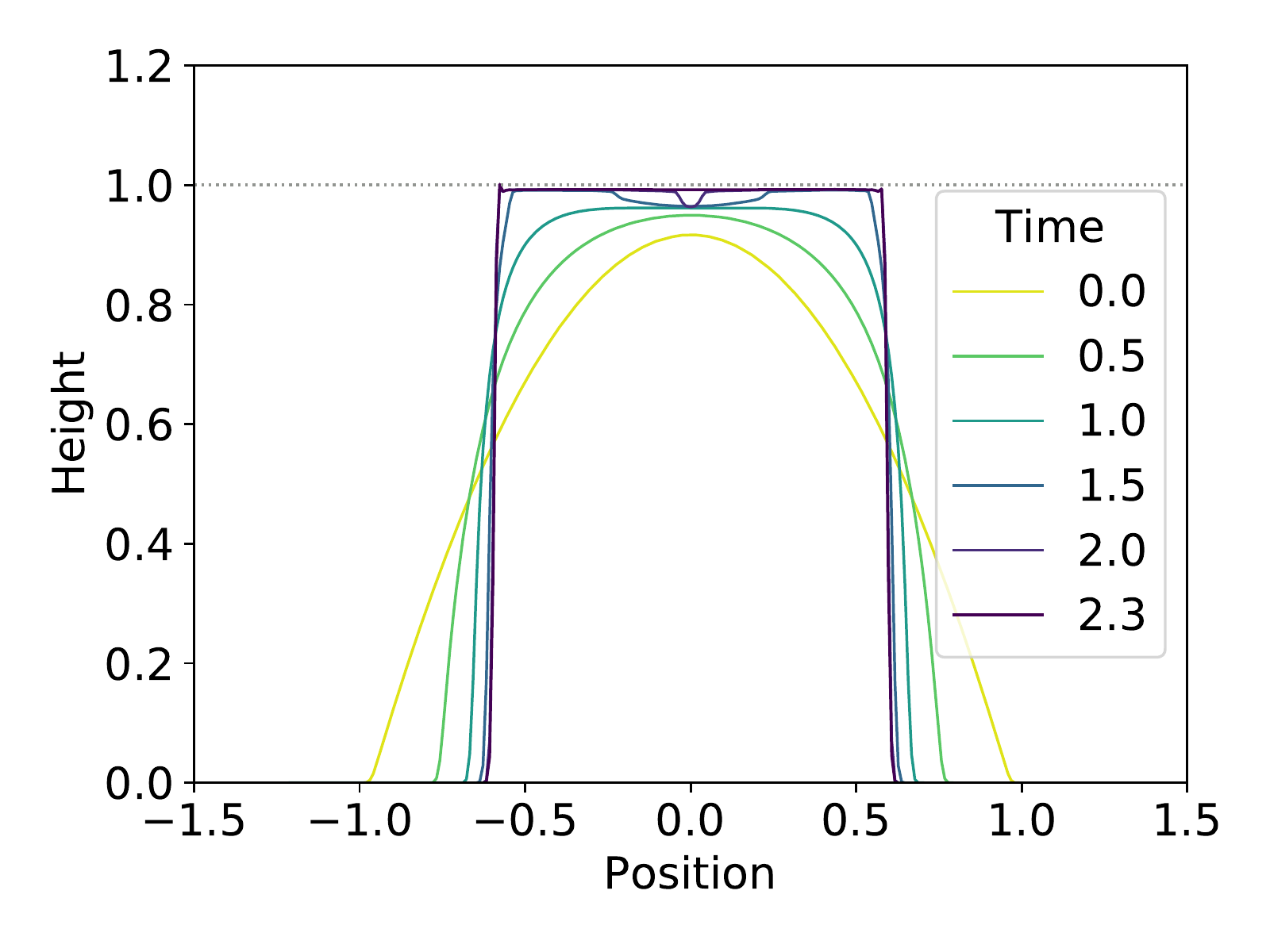}
\includegraphics[height=4.35cm,trim={2.2cm 1.25cm .6cm .6cm},clip, valign=t]{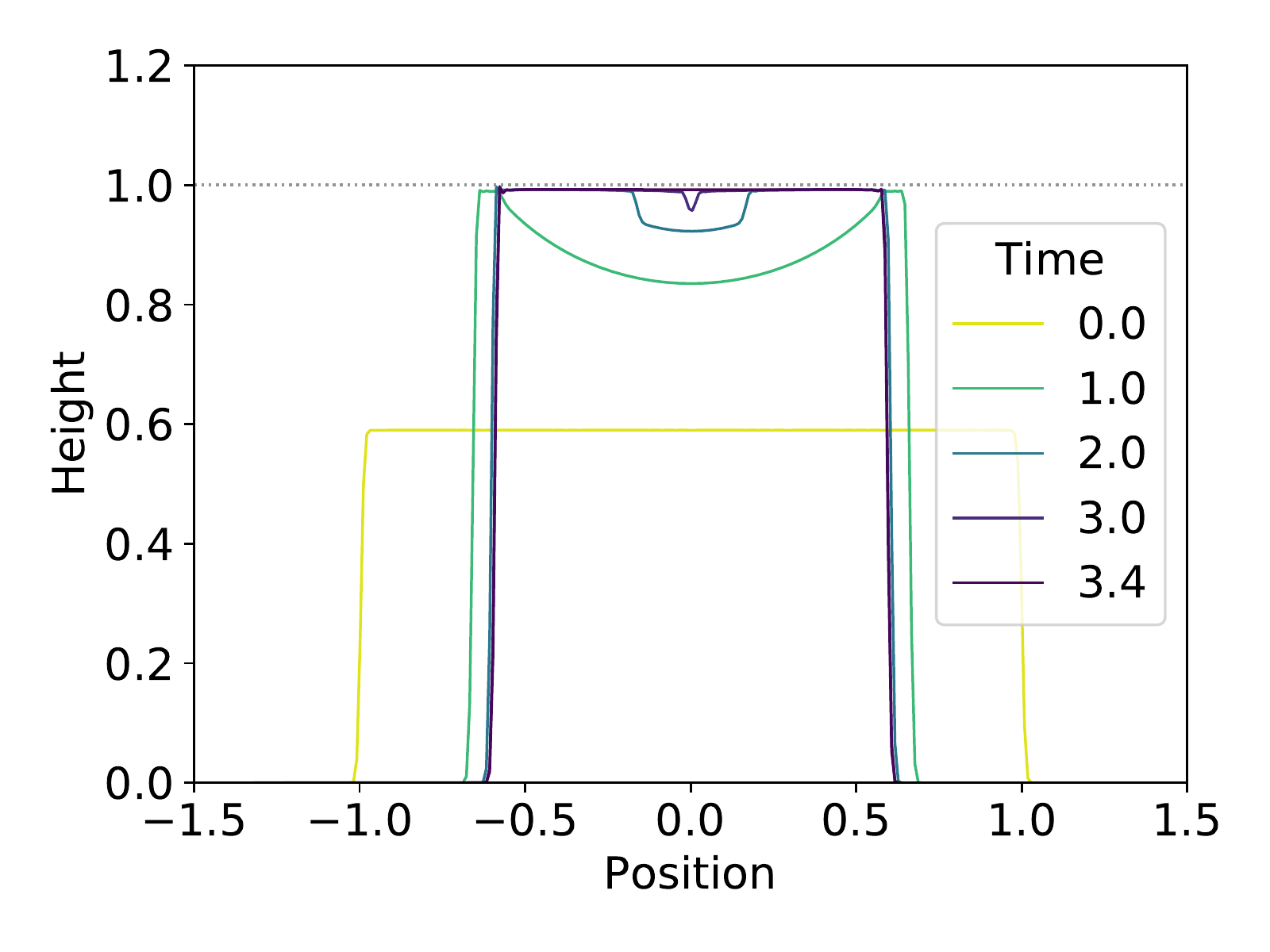}

\caption{We approximate gradient flows of $\E_\infty$ by simulating solutions of aggregation-diffusion equations for $m = 800$, $K(x) = |x|^q/q - |x|$,  $h = 0.004$, $k = 10^{-4}$, and  initial data given by constant multiples of Barenblatt profiles (equation (\ref{Barenblatt}), $m_* =2$, $\tau = 0.1$) and patch functions (equation (\ref{patch initial}),  $\Omega = [-1,1]$). We multiply the initial data by the desired   mass for each $q$, i.e.  mass 0.68, 1.00, and 1.18,  for   $q = 1.6$, 2.0, and 2.4.}
\label{varyingtimes}
\end{center}
\end{figure}

In Figure \ref{varyingtimes}, we approximate gradient flows of the constrained interaction energy $\E_\infty$ for various choices of  attraction parameter $q \geq 1$, with initial data at the critical mass from Figure \ref{critical mass figure}, i.e., the smallest value of mass for which solutions approach a set valued equilibrium. We contrast the behavior of initial data that is a constant multiple of a Barenblatt profile (equation (\ref{Barenblatt}), $m_*=2$, $\tau =1$)) with initial data that is a constant multiple of a patch function (equation (\ref{patch initial}), $\Omega = [-1,1]$). In both cases, gradient flows converge to a characteristic function of height one on an interval centered at the origin. When  $q < 2$, solutions initially reach height one at the center of mass of the density and then spread to become a characteristic function. When  $q>2$, solutions initially reach height one at the boundary of the support of the density and then ``fill in'' the interior to become a characteristic function.

\begin{figure}[h]
\begin{center}
 {\bf Constrained Aggregation: Equilibria for Varying Masses}  \\ 
\vspace{.1cm}

\hspace{.75cm} {\bf $q=1.4$} \hspace{3.75cm} {\bf $q=2$} \hspace{3.75cm} {\bf $q=2.6$}

\includegraphics[height=4cm,trim={.6cm .7cm .6cm .6cm},clip, valign=t]{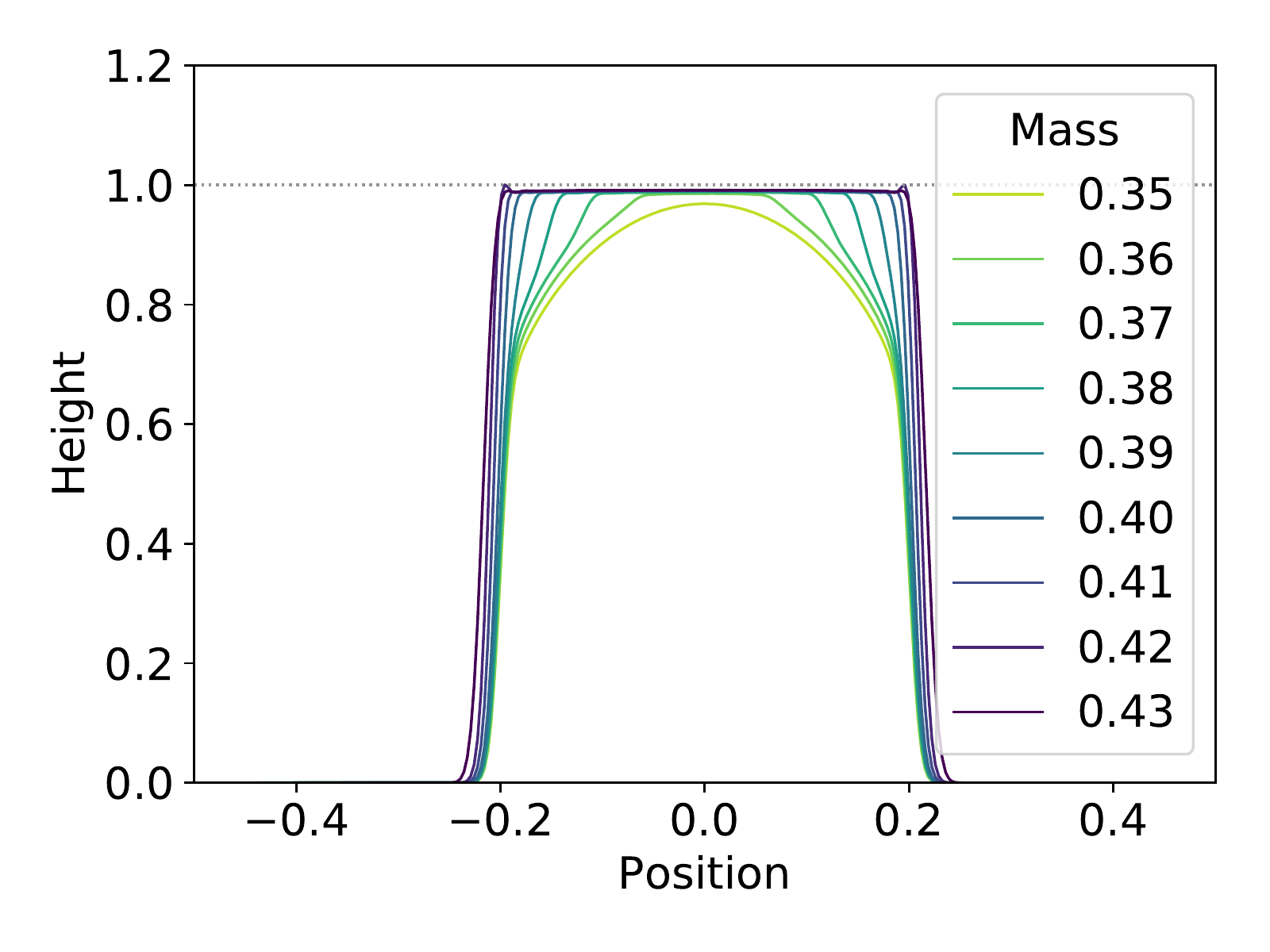} 
\includegraphics[height=3.8cm,trim={2.25cm 1.25cm .6cm .6cm},clip, valign=t]{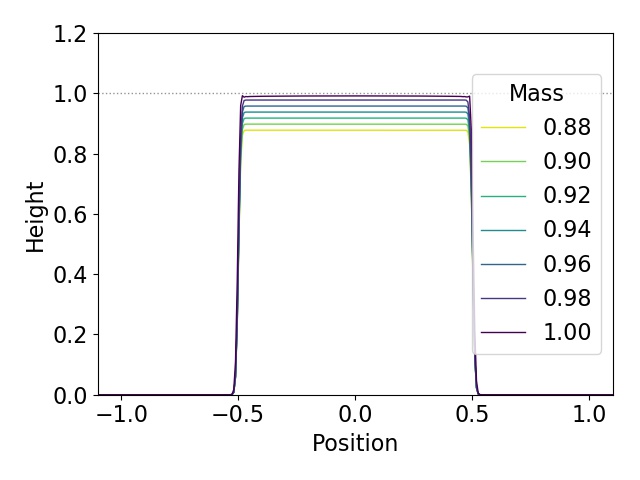}
\includegraphics[height=3.8cm,trim={2.25cm 1.25cm .6cm .6cm},clip, valign=t]{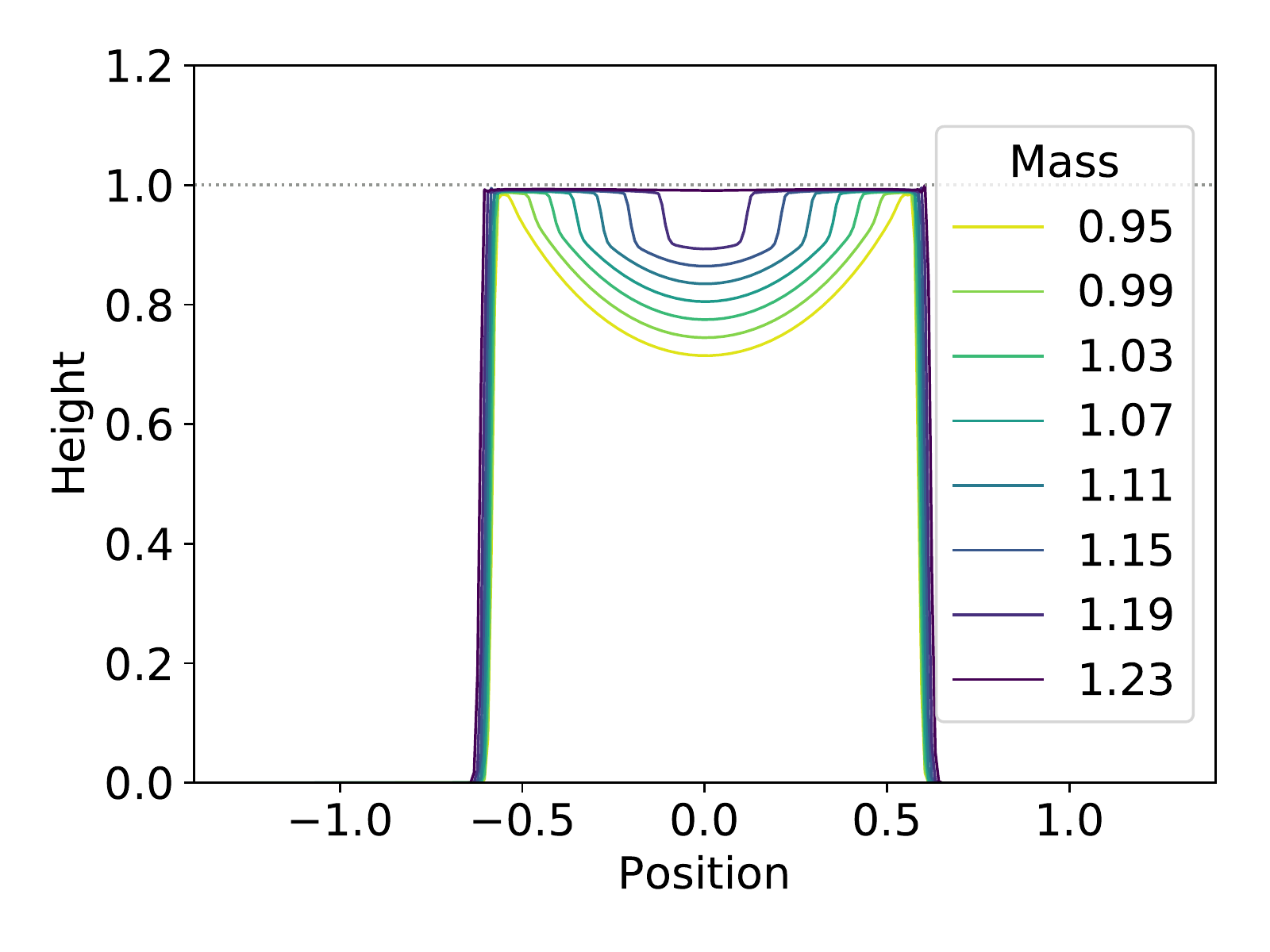}

\caption{We approximate minimizers of $\E_\infty$ by simulating solutions of aggregation-diffusion equations for $m = 800$, $K(x) = |x|^q/q - |x|$,  $h = 0.004$, $k = 10^{-4}$, with  initial data given by constant multiples of Barenblatt profiles (equation (\ref{Barenblatt}), $m_* =2$, $\tau = 0.1$).}
\label{varyingmasses}
\end{center}

\end{figure}

Finally, in Figure \ref{varyingmasses}, we compute minimizers of the interaction energy $\E_\infty$ for varying choices of attraction exponent $q \geq 1$ and masses up to and including the critical mass from Figure \ref{critical mass figure}. These simulations appear to confirm the existence of an \emph{intermediate phase} between the liquid and solid phase as the generic behavior for $q \neq 2$. In the notation of the introduction, these simulations suggest that $M_1^* < M_2^*$ for $q \neq 2$. In particular, we observe that minimizers of mass $M$ satisfy $|\{ \rho = 1\}| \in (0,M)$ for all $M \in [0.36, 0.42]$ when $q= 1.4$ and for all $M \in [0.99,1.19]$ for $q = 2.6$.

\section{Appendix}
\subsection{Regular Functionals}
Consider a functional $\F: \P_2(\Rd) \to (-\infty, +\infty]$ that is proper, lower semicontinuous, and satisfies $D(|\partial \F|) \subseteq \P_{2, ac}(\Rd)$. Then $\xi \in \partial \F(\mu)$ is a \emph{strong subdifferential} if for all Borel measurable functions $\bt : \Rd \to \Rd$ we have
\[ \F(\bt \# \mu) - \F(\mu) \geq \int \la \xi, \bt - \id \ra\, d \mu + o ( \|\bt - \id\|_{L^2(\mu)}) . \]
Following Ambrosio, Gigli, and Savar\'e \cite[Definition 10.1.4]{AGS}, we define the notion of \emph{regular functional} as follows.
\begin{definition} \label{regulardefinition} Given $\F: \P_2(\Rd) \to (-\infty, +\infty]$   proper, lower semicontinuous, and satisfying $D(|\partial \F|) \subseteq \P_{2, ac}(\Rd)$, $\F$ is a \emph{regular} functional if for any 2-Wasserstein convergent sequence $\mu_n  \to \mu $ with strong subdifferentials $\xi_n \in \partial \F(\mu_n)$ satisfying
\begin{enumerate}
\item $\sup_n | \F(\mu_n) | < +\infty $
\item $ \sup_n \|\xi_n\|_{L^2(\mu_n)} < +\infty$
\item there exists $ \xi \in L^2(\mu) \text{ so that }\lim_{n \to +\infty} \int f \xi_n \,d\mu_n = \int f \xi \,d\mu  \text{ for all } f \in C^\infty_c(\Rd),$
\end{enumerate}
we have $\lim_{n \to +\infty} \F(\mu_n) = \F(\mu)$ and $\xi \in \partial \F(\mu)$.
\end{definition}

We now provide a sufficient condition on the subdifferential that ensures the energy is regular and the metric slope is a strong upper gradient. This generalized Ambrosio, Gigli, and Savar\'e's result that $\lambda$-convex energies are regular \cite[Lemma 10.1.3]{AGS}.
\begin{proposition} \label{first regular proposition}
Suppose $\F: \P_2(\Rd) \to (-\infty, +\infty]$ is   proper, lower semicontinuous, and satisfies $D(|\partial \F|) \subseteq \P_{2, ac}(\Rd)$. Furthermore, suppose that there exists  a continuous function $\psi: [0, +\infty) \to [0, +\infty)$ with $\psi(0) = 0$ so that, for any $\xi \in \partial \F(\mu)$,
\begin{align} \label{gensubdiffineq} \F(\nu) - \F(\mu) \geq \int \la \xi , t_{\mu}^\nu - \id \ra\, d \mu - f(\mu,\nu) \psi(d_W(\mu,\nu)) d_W(\mu,\nu) \  \text{ for all } \nu \in D(\F) ,
\end{align}
where $f(\mu,\nu) = C(1+ \F(\mu) + \F(\nu))$ for some $C = C(d_W(\mu,\nu))>0$ which is an increasing function of the distance from $\mu$ to $\nu$.
Then $\F$ is regular and the metric slope $|\partial \F|$ is a strong upper gradient for $\F$.
\end{proposition}
\begin{proof}
We begin by showing that $\F$ is regular.
Consider a 2-Wasserstein convergent sequence $\mu_n \to \mu $ with strong subdifferentials $\xi_n \in \partial \F(\mu_n)$ satisfying criteria (i)-(iii) from Definition \ref{regulardefinition}. First, we show that $\xi \in \partial \F(\mu)$. By assumption,
\begin{multline} \label{reg inequality} \F(\nu) - \F(\mu_n)  \\ \geq \int \la \xi_n , t_{\mu_n}^\nu - \id \ra \, d \mu_n - f(\mu_n,\nu)\psi(d_W(\mu_n,\nu)) d_W(\mu_n,\nu) \  \text{ for all } \nu \in \P_2(\Rd). 
\end{multline}
By the lower semicontinuity of $\F$, we have $\liminf_{n \to +\infty} \F(\mu_n) \geq \F(\mu)$. Since $\sup_n \F(\mu_n) < +\infty$ and  $f$ is  locally
bounded on sublevels of $\F$, we have $\sup_n f(\mu_n,\nu) < +\infty$. Since $\psi$ is continuous,
\[ \lim_{n \to +\infty}  \psi(d_W(\mu_n,\nu)) d_W(\mu_n,\nu) = \psi(d_W(\mu,\nu))d_W(\mu,\nu) . \]
 Furthermore, arguing as in \cite[Lemma 10.1.3]{AGS}, we have
 \[ \lim_{n \to +\infty} \int \la \xi_n , t_{\mu_n}^\nu - \id \ra \, d \mu_n = \int \la \xi , t_{\mu}^\nu - \id \ra \, d \mu . \]
 Therefore, sending $n \to +\infty$ in (\ref{reg inequality}), we conclude that $\xi \in \partial \F(\mu)$. 

Now, we show $\lim_{n \to +\infty} \F(\mu_n) = \F(\mu)$. Taking $\nu = \mu$ in (\ref{reg inequality}) and sending $n \to +\infty$, the previous argument shows that the right hand side converges to 0. Thus,
\[ \liminf_{n \to +\infty} \F(\mu) - \F(\mu_n) \geq 0 \iff \liminf_{n \to +\infty} -\F(\mu_n) \geq -\F(\mu) \iff \limsup_{n \to +\infty} \F(\mu_n) \leq \F(\mu) . \]
Combining this with the lower semicontinuity of $\F$, we obtain that $\lim_{n \to +\infty} \F(\mu_n) = \F(\mu)$. Therefore, $\F$ is regular.

We now show that the metric slope $|\partial \F|$ is a strong upper gradient for $\F$. We argue as in \cite[Corollary 2.4.10]{AGS}. Consider $\mu:(0,T) \to \P_2(\Rd)$ that is \emph{absolutely continuous} in time, i.e., there exists $m \in L^1(0,T)$ so that inequality (\ref{absctsdef}) holds. 
Then its metric derivate $|\mu'|(t)$ is well defined.  It suffices to show that if $\mu(t)$ satisfies $|\partial \F|(\mu) |\mu'| \in L^1(0,T)$, then $\F(\mu(t))$ is absolutely continuous in time. By \cite[Lemma 10.1.5]{AGS}, we have $|\partial \F|(\mu)<+\infty$ if and only if there exists some $\xi \in \partial \F(\mu)$ so that $|\partial \F|(\mu) = \|\xi\|_{L^2(\rho)}$.

We begin by showing that $\F(\mu(t)) |\mu'|(t) \in L^1(0,T)$. Without loss of generality, $\F(\mu(t_*)) < +\infty$ for some $t_* \in (0,T)$, or else $\F(\mu(t)) \equiv +\infty$ is constant, hence absolutely continuous. Furthermore, up to reparametrizing time, we may also assume
\begin{align*}
\sup_{t \in (0,T)}  C \psi(d_W(\mu(t),\mu(t_*)))  d_W(\mu(t), \mu(t_*)) \leq \frac{1}{2} ,
\end{align*}
for $C= C(d_W(\mu,\nu)) >0$
as in the definition of $f$.
Applying inequality (\ref{gensubdiffineq}) with $\mu = \mu(t)$ and $\nu = \mu(t_*)$, we conclude that there exists  $C'>0$ so that
\begin{align*}
\F(\mu(t)) &\leq  \F(\mu(t_*)) + |\partial \F|(\mu(t)) d_W(\mu(t),\mu(t_*)) \\
				&\qquad\qquad\qquad+  f(\mu(t),\mu(t_*)) \psi(d_W(\mu(t),\mu(t_*))) d_W(\mu(t),\mu(t_*)) \\
&\leq  \F(\mu(t_*)) + C'|\partial \F|(\mu(t))    +  \frac{1}{2} \Big( 1+ \F(\mu(t))+\F(\mu(t_*)) \Big)  
\end{align*}
Rearranging and multiplying by $|\mu'|(t)$, we conclude that $\F(\mu(t)) |\mu'|(t) \in L^1(0,T)$.

We now show that $\F(\mu(t))$ is absolutely continuous in time. Consider the compact subset $\S := \{ \mu(t) : t \in [0,T] \} \subset \P_2(\Rd)$ and recall that the global slope on this subset
\[ I_\F(\mu):= \sup_{\nu \in \S, \nu \neq \mu} \frac{(\F(\mu) - \F(\nu))^+}{d_W(\mu,\nu)} \]
is a strong upper gradient \cite[Theorem 1.2.5]{AGS}. 
By inequality (\ref{gensubdiffineq}),
\begin{align} \label{lcsslope1}    \frac{( \F(\mu) - \F(\nu))^+}{d_W(\mu,\nu)} \leq |\partial \F|(\mu) + C \Big(1+\F(\mu) + \F(\nu) \Big) \psi(d_W(\mu,\nu)) . 
\end{align}
Furthermore, we may assume without loss of generality that $\F(\nu) \leq \F(\mu)$. Therefore,
\[ I_\F(\mu(t)) |\mu'|(t) \leq |\partial \F|(\mu(t))|\mu'|(t) + C\Big( 1+ 2\F(\mu(t))  \Big) \psi(\text{diam} \S) |\mu'|(t)  . \]
Therefore, we conclude that $ I_\F(\mu(t)) |\mu'|(t) \in L^1(0,T)$. Arguing as in \cite[Theorem 1.2.5]{AGS}, we conclude that $\F(\mu(t))$ is absolutely continuous.
  \end{proof}

\begin{corollary} \label{regular corollary}
Suppose $K$ satisfies hypotheses \ref{GFcont}, \ref{GFbounds}--\ref{GFdualsobolev} and $m \geq m_0$. Then $\E_m$ is \emph{regular} and the metric slope $|\partial \E_m|$ is  a strong upper gradient for $\E_m$.
\end{corollary}
\begin{proof}
The result is an immediate consequence of Proposition \ref{LSCproposition}, Proposition \ref{subdiffalmostconvex}, and Proposition \ref{first regular proposition}, where we appeal to Proposition \ref{firstLmbound} to ensure that the estimates
\[ C(1+ \|\rho_0\|_m +\|\rho_1\|_m ) \leq C_r \Big(1 + \|\rho_0\|_m^{1+ \frac{1}{r}} + \|\rho_1\|_m^{1+ \frac{1}{r}} \Big) \leq \tilde{C}_r(1 + \E_m(\rho_1) + \E_m(\rho_2) )  \]
hold.
 \end{proof}

 \subsection{Power-Law Interaction Potentials} \label{powerlawappendix}
 In this section, we prove Proposition \ref{Kderivative} and Proposition  \ref{subdifftheorem}.

\begin{proof}[Proof of Proposition \ref{Kderivative}]
The result extends \cite[Proposition 4.6]{craig2017nonconvex}, where 
 \ref{GFbounds}, \ref{GFcty}, and  \ref{GFdualsobolev} generalize \cite[Assumption 4.1]{craig2017nonconvex}. 
Define $\bt_\alpha :=(1-\alpha) \id + \alpha \bt_{\rho_0}^{\rho_1}$ and let $\rho_\alpha := \bt_\alpha \# \rho_0$ to be the Wasserstein geodesic from $\rho_0$ to $\rho_1$. By convexity of the $L^p$-norm, defined by equation \eqref{Lpnormdef}, along Wasserstein geodesics, $\|\rho_\alpha\|_m \leq \max \{ \|\rho_0\|_m, \|\rho_1\|_m \}$ for all $\alpha \in [0,1]$.
Then,
	\beqn 
		\begin{aligned} 
\frac{d}{d \alpha} \K(\rho_\alpha) &= \lim_{h \to 0} \frac{1}{h} [\K(\rho_{\alpha+h}) - \K(\rho_\alpha)]   \\
&= \lim_{h \to 0} \frac{1}{2h} \left[ \int K*\rho_\alpha \, d\rho_{\alpha +h}  - \int K* \rho_{\alpha}\,  d \rho_\alpha \right]  \\
&\qquad\qquad +  \frac{1}{2h} \left[ \int  K* \rho_{\alpha+h} \,d\rho_{\alpha +h} - \int K*\rho_{\alpha+h}\, d\rho_{\alpha} \right]  \\
&=  \lim_{h \to 0} \frac{1}{2h} \int [(K*\rho_\alpha)\circ \bt_{\alpha +h}   -(K* \rho_{\alpha})\circ \bt_\alpha ]\,d \rho_0   \\
&\qquad \qquad +\frac{1}{2h} \int  [(K* \rho_{\alpha+h})\circ \bt_{\alpha + h} - (K*\rho_{\alpha+h})\circ \bt_{\alpha} ]\, d\rho_0   \\
&= \lim_{h \to 0} \frac{1}{2h} \int \Big( k_\alpha \circ \bt_{\alpha +h} - k_\alpha \circ \bt_\alpha \Big)\,d \rho_0 +\frac{1}{2h}  \int\Big( k_{\alpha + h}\circ \bt_{\alpha +h} - k_{\alpha + h} \circ \bt_{\alpha}\Big)\, d \rho_0 , 
		\end{aligned} \nonumber
	\eeqn
for $k_{\beta}(x) := (K * \rho_{\beta})(x)$. We consider both terms simultaneously by taking $\beta = \alpha$ or $\alpha + h$. By hypothesis \ref{GFcty},   $k_{\beta}(x)$ is continuously differentiable with respect to $x$, so
\begin{align*}
 k_{\beta} \circ \bt_{\alpha + h} - k_{\beta} \circ \bt_\alpha  &=  \int_\alpha^{\alpha +h} \frac{d}{d \gamma} k_{\beta} \circ \bt_\gamma \,d\gamma =  \int_{\alpha}^{\alpha +h} \big\la (\grad K * \rho_{\beta})\circ \bt_\gamma, \bt_{\rho_0}^{\rho_1} - \id \big\ra\, d \rho_0 d \gamma .
\end{align*}
Furthermore, since $\|\rho_{\beta} \|_m \leq \max \{ \|\rho_0\|_m, \|\rho_1\|_m\}$, \ref{GFbounds} ensures $\| \grad K * \rho_{\beta} \|_{L^2(\rho_\gamma)} <+\infty$.
Consequently, we may interchange the order of integration, and add and subtract to obtain
	\beqn \label{heightomegaconvex2}
		\begin{aligned} 
\frac{d}{d \alpha} \K(\rho_\alpha)  &= \lim_{h \to 0} \frac{1}{h} \int_{\alpha}^{\alpha + h}  \int  \big\la (\grad K * \rho_{\alpha})\circ \bt_\gamma,\bt_{\rho_0}^{\rho_1}- \id\big\ra\, d \rho_0 d \gamma   \\
& \qquad + \frac{1}{2h}  \int_\alpha^{\alpha +h}\int \big\la (\grad K * \rho_{\alpha +h})\circ \bt_\gamma - (\grad K *\rho_\alpha) \circ \bt_\gamma, \bt_{\rho_0}^{\rho_1}- \id \big\ra\, d \rho_0 d \gamma. 
		\end{aligned}
	\eeqn

In order to compute the limits on the right hand side, note that for any $\alpha,\, \tilde{\alpha},\, \beta, \, \tilde{\beta} \in [0,1]$,

	\beqn \label{heightomegaconvex3}
		\begin{aligned} 
& \int \big|(\grad K *\rho_{\tilde{\alpha}}) \circ \bt_{\tilde{\beta}} - (\grad K* \rho_{\alpha}) \circ \bt_{\beta}\big| \, \big|\bt_{\rho_0}^{\rho_1}-\id\big|\, d \rho_0  \\
&\qquad \leq  \Big[  \| (\grad K *\rho_{\tilde{\alpha}} )\circ \bt_{\tilde{\beta}} - (\grad K * \rho_{\tilde{\alpha}}) \circ \bt_\beta \|_{L^2(\rho_0)}  \\ 
	&\qquad\qquad\qquad + \| \grad K * \rho_{\tilde{\alpha}} - \grad K *\rho_{\alpha}\|_{L^2(\rho_\beta)}\, \Big] \|\bt_{\rho_0}^{\rho_1} - \id\|_{L^2(\rho_0)}  \\
 &\qquad \leq f(\rho_0,\rho_1)  \left[   \left\| \psi \big( |\bt_{\tilde{\beta}} - \bt_{\beta} |^2 \big) \right\|_{L^1(\rho_0)}^{1/2} + \psi(d_W(\rho_{\tilde{\alpha}},\rho_\alpha) )\right] d_W(\rho_0,\rho_1)  \\
    &\qquad \leq f(\rho_0,\rho_1) \left[   \psi \left(\|  \bt_{\tilde{\beta}} - \bt_{\beta}   \|_{L^2(\rho_0)}^2 \right)^{1/2} + \psi \left(d_W(\rho_{\tilde{\alpha}},\rho_\alpha) \right) \right] d_W(\rho_0,\rho_1)  \\
      &\qquad \leq f(\rho_0,\rho_1)\left[   \psi \left(|\tilde{\beta} - \beta|^2 d_W(\rho_{0},\rho_1)^2 \right)^{1/2} + \psi \Big(|\tilde{\alpha} - \alpha| d_W(\rho_{0},\rho_1) \Big) \right] d_W(\rho_0,\rho_1) .  \\
		\end{aligned}
	\eeqn
In the second inequality, we use hypotheses \ref{GFcty} and \ref{GFdualsobolev} and the fact that the 2-Wasserstein distance dominates the $(2-\epsilon)$-Wasserstein distance for all $\epsilon >0$; see inequality \ref{Wasserstein order}. In the third inequality we use Jensen's inequality for the concave function $\psi(s)$. 

When $\alpha = \tilde{\alpha}$, this estimate ensures
\[  \beta \mapsto \int  \big\la (\grad K * \rho_{\alpha})\circ \bt_\beta, \bt_{\rho_0}^{\rho_1} -\id \big\ra \, d\rho_0 \]
is continuous, so that the first term in (\ref{heightomegaconvex2}) converges to $ \int  \big\la (\grad K * \rho_{\alpha})\circ \bt_\alpha, \bt_{\rho_0}^{\mu_1} - \id \big\ra\, d \rho_0$. Likewise, when $\beta = \tilde{\beta}$, this estimate guarantees that the second term in (\ref{heightomegaconvex2}) is bounded by
\[ \lim_{h \to 0} \frac{f(\rho_0,\rho_1) }{2h} \int_\alpha^{\alpha + h} \psi \Big(h \, d_W(\rho_0,\rho_1) \Big) d_W(\rho_0,\rho_1)\, d \gamma = 0. \]
 Therefore, we conclude
\begin{align*} 
 \frac{d}{d\alpha} \K(\rho_\alpha) = \int  \big\la (\grad K * \rho_{\alpha})\circ \bt_\alpha, \bt_{\rho_0}^{\rho_1} - \id \big\ra \, d \rho_0.
 \end{align*}
 
By (\ref{heightomegaconvex3}) again, $\frac{d}{d\alpha} \K(\mu_\alpha)$ is continuous for $\alpha \in [0,1]$. Therefore,
\beqn
	\begin{aligned} \label{Taylorexpansion}
 \K(\rho_1)  &= \K(\rho_0) +  \int_0^1 \frac{d}{d\alpha} \K(\rho_\alpha) \, d \alpha \\
&= \K(\rho_0) + \left. \frac{d}{d\alpha} \K(\rho_\alpha) \right|_{\alpha =0} + \int_0^1\int  \big\la (\grad K*\rho_\alpha) \circ \bt_\alpha - \grad K*  \rho_0,\bt_{\rho_0}^{\rho_1} - \id \big\ra \, d \rho_0  d\alpha. 
	\end{aligned}
\eeqn
To prove the result, it suffices to show that the third term is $o\left(d_W(\rho_0,\rho_1) \right)$.
This follows by a final application of  inequality (\ref{heightomegaconvex3}):
\beqn
	\begin{aligned} \label{moredetailsubdiff}
&\left| \int_0^1\int  \big\la (\grad K*\rho_\alpha) \circ \bt_\alpha - \grad K*  \rho_0,\bt_{\rho_0}^{\rho_1} - \id \big\ra \, d \rho_0  d\alpha \right| \\
&\qquad \leq  f(\rho_0,\rho_1) \left[   \psi \left(  d_W(\rho_{0},\rho_1)^2 \right)^{1/2} + \psi \Big( d_W(\rho_{0},\rho_1) \Big) \right] d_W(\rho_0,\rho_1) \nonumber
	\end{aligned}
\eeqn
Finally,   by Lemma \ref{psidoublingbound}, there exists $c$, which is an increasing function of $d_W(\rho_0,\rho_1)$,  so that  $\sqrt{\psi(s^2)} \leq c \psi(s) $. Therefore, up to increasing the constant in the definition of $f(\rho_0,\rho_1)$, inequality (\ref{moredetailsubdiff}) is bounded by 
\[ f(\rho_0,\rho_1) \, \psi(d_W(\rho_0,\rho_1))  \, d_W(\rho_0,\rho_1) , \]
which completes the proof.
\end{proof}

We now turn to the proof of Proposition  \ref{subdifftheorem}.
\begin{proof}[Proof of Proposition \ref{subdifftheorem}]

Our proof follows a similar approach as \cite[Theorem 10.4.13]{AGS}, generalizing this result to nonconvex, singular interaction potentials $K$. We begin by proving (\ref{subdiffeqnEm}). Note that the implication $\impliedby$ is immediate, as the definition of the subdifferential  $\xi \in \partial \E_m(\rho)$ requires $\|\xi\|_{L^2(\rho)}< +\infty$. (See inequality (\ref{subdiffdef}).)

Suppose $|\partial \E_m|(\rho) < +\infty$. Since $\rho  \in D(\E_m)$, we have $\|\rho\|_m < +\infty$. Fix $\xi \in L^2(\rho)$ and define  $r_\alpha = (1-\alpha) \id + \alpha \xi$ and $\rho_\alpha:= (r_\alpha)_\# \rho$. Suppose  either (i) $\xi - \id \in C^\infty_c(\Rd; \Rd)$ or (ii) $\xi = 0$.
Note that both assumptions ensure that $\xi$ is differentiable almost everywhere (cf. Aleksandrov's theorem,\cite[Theorem 5.5.4]{AGS}) and there exists $\alpha_0 , C>0$ so that for all $\alpha \in [0,\alpha_0]$, $\|\rho_\alpha\|_m < C$(cf. \cite[Lemma 5.5.3]{AGS}).

Under either assumption (i) or (ii), the definition of the metric slope (\ref{metricslopedef}), Proposition \ref{Kderivative} for  $\K$, and \cite[Lemma 10.4.4]{AGS} for  $\Se_m$, ensure
	\beqn \begin{aligned} \label{ibpawesome1}
|\partial \E_m(\rho)| \,  \|\xi - \id\|_{L^2(\rho)} & \geq \limsup_{\alpha \to 0} \frac{\E_m(\rho_\alpha) - \E_m(\rho)}{d_W(\rho,\rho_\alpha)} \|\xi - \id\|_{L^2(\rho)}  \\
&= \limsup_{\alpha \to 0}  \frac{ \K(\rho_\alpha) - \K(\rho)}{\alpha}  + \frac{\Se_m(\rho_\alpha) - \Se_m(\rho)}{\alpha}  \\
&=  \int \big\la \grad K*\rho, \xi - \id \big\ra\, d  \rho   - \int \rho^m \grad \cdot (\xi -\id)\,dx .
	\end{aligned} \eeqn
Thus, by H\"older's inequality and hypothesis \ref{GFbounds},   there exists $C^\pr >0$ depending on $\rho$ so that
\begin{align} \label{ibpawesome}
|\partial \E_m(\rho)| \,  \|\xi-\id\|_{L^2(\rho)} + C^\pr \|\xi-\id\|_{L^2(\rho)} \geq - \int \rho^m \grad \cdot (\xi - \id) \,dx. 
\end{align}

First, we suppose $\xi$ satisfies assumption (i). As inequality (\ref{ibpawesome}) holds for all $f = \xi - \id \in C^\infty_c(\Rd)$, this implies $\rho^m$ is a function of bounded variation and the right hand side may be rewritten as $\int \grad \rho^m \cdot (\xi- \id)\,dx$. Applying inequality \eqref{ibpawesome} again, we obtain  $\grad \rho^m   \in L^2(\rho)$. Returning back to (\ref{ibpawesome1}), we obtain that 
\begin{align} \label{wchar}   \left\| (\grad K*\rho)  + \frac{\grad \rho^m}{\rho} \right\|_{L^2(\rho)} \leq |\partial \E_m|(\rho) . 
\end{align}
Next, suppose  $\xi = 0$. Then inequality (\ref{ibpawesome}) gives
\begin{align*}
(|\partial \E_m(\rho)| +C^\pr) M_2(\rho)^{1/2} \geq d \left( \int \rho^m \,dx \right).
\end{align*}
hence $\rho^m \in W^{1,1}(\Rd)$.

To complete the proof, it suffices to show  $(\grad K*\rho)  + \frac{\grad \rho^m}{\rho}  \in \partial \E_m(\rho)$. Remark \ref{subdiffmetricslope}  ensures  equality must hold in  (\ref{wchar}). Since the subdifferential $\partial \E_m(\rho)$ is a convex subset of $L^2(\rho)$ and the $L^2(\rho)$-norm is strictly convex, the element of $\partial \E_m(\rho)$ with minimal $L^2(\rho)$-norm is unique.

By   Proposition \ref{subdiffalmostconvex}, to show  $(\grad K*\rho)  + \frac{\grad \rho^m}{\rho}  \in \partial \E_m(\rho)$ it suffices to show that 
	\begin{multline} 
\E_m(\nu) - \E_m(\rho) \\ \geq \int \Big\la (\grad K*\rho)  + \frac{\grad \rho^m}{\rho}, \bt_\rho^\nu - \id \Big\ra \, d\rho -f(\rho,\nu) \psi(d_W(\rho,\nu))d_W(\rho,\nu) , \quad \forall \ \nu \in D(\E_m). \nonumber
\end{multline}
	Since $\rho, \nu \in D(\E_m)$,   by Proposition \ref{Kderivative}, for the subdifferential of $\K$, and by \cite[Theorem 10.4.6]{AGS}, for the subdifferential of $\Se_m$,
\begin{gather}	\K(\nu) - \K(\rho) =  \int \big\la \grad K*\rho, \bt_{\rho}^{\nu} - \id \big\ra\, d \rho -f(\rho,\nu) \psi(d_W(\rho,\nu))d_W(\rho,\nu), \label{Ksub1}\\
\Se_m(\nu) - \Se_m(\rho) \geq \int \Big\la \frac{\grad \rho^m}{\rho}, \bt_{\rho}^{\nu} - \id \Big\ra\, d \rho \label{Ssub1}
\end{gather} 
Adding together inequalities (\ref{Ksub1}) and (\ref{Ssub1}) gives the result.
\end{proof}

\subsection{Elementary Bounds}
\begin{lemma} \label{psidoublingbound}
Suppose $\psi:[0,+\infty)\to[0,+\infty)$ is a continuous, nondecreasing, concave function  with  $\psi(0) = 0$. Then for any $\bar{s}$, there exists $\bar{C}$ so that
\[ \psi(s^2) \leq \bar{C} \psi(s)^2 , \quad \text{ for all } s \in [0, \bar{s}] . \]
\end{lemma}
\begin{proof}
Since $\psi(s)$ is concave and $\psi(0) = 0$, $\psi(s)/s$ is a decreasing function and $\partial^+ \psi(0) =\lim_{s \to 0^+} \psi(s)/s$ exists. First, suppose $\partial^+ \psi(0) = 0$. Since $\psi$ is concave and nondecreasing, this implies $\psi \equiv 0$, and the result holds.

Now, suppose $\partial^+ \psi(0) >0$.
It suffices to show that $\psi(s)^2/\psi(s^2)$ is uniformly bounded below on $[0,\bar{s}]$. 
Since $\partial^+ \psi(0) >0$ and $\psi$ is increasing, $\psi(s)>0$ for all $s>0$ and $\psi(s)^2/\psi(s^2)$ is continuous and positive on $(0,\bar{s}]$, hence bounded below away from $s = 0$. Furthermore,
\[ \lim_{s \to 0} \frac{\psi(s)^2}{\psi(s^2)} = \partial^+ \psi(0) >0 . \]
Thus $\psi(s)^2/\psi(s^2)$ is bounded below on $[0, \bar{s}]$.
\end{proof}

\begin{lemma} \label{elemineq}
For $s \geq 0$ and $0 <b<1/2$, define $\psi(s,a):= s^{(1-a)}$. Then,
\[ |\psi(s,b) - \psi(s,0)| \leq b (1+ s^{2-b}) . \]
\end{lemma}

\begin{proof}
Define $\psi(s,a):= s^{(1-a)}$, so $\frac{d}{d a} \psi(s,a) = -s^{1-a} \log(s)$. By the mean value theorem, for $0 <b<1/2$, there exists $a \in [0,b]$ so that
\[ \psi(s,b) - \psi(s,0) = b \frac{d}{d a} \psi(s,a) . \]
If $s \in [0,1]$, $ \left|  \frac{d}{d a} \psi(s,a)  \right| \leq 1$, so
\begin{align} \label{sleq1claim}
 \left| \psi(s,b) - \psi(s,0)  \right| \leq b .
 \end{align}

If $s \geq 1$, we claim that 
\begin{align} \label{sgeq1claim}
\left| \psi(s,b) - \psi(s,0) \right| \leq b s^{2-b}
\end{align}
This holds since, for $s>1$,
\begin{align*}
\left| \psi(s,b) - \psi(s,0) \right| = \left| s^{1-b} - s \right| = s - s^{1-b} \leq b s^{2-b} 
& \iff s^{1+b} - s \leq bs^2  \iff s^{1+b}  \leq bs^2 + s.
\end{align*}
This is true at $s=1$, so it suffices to show the derivative with respect to $s$ of the right-hand side is larger than the derivative of the left-hand side, i.e.,
\[ (1+b) s^b \leq 2b s + 1 \]
This is also true at $s=1$, so differentiating again, it suffices to show
\[ (1+b)bs^{(b-1)} \leq 2b \]
This holds since $s \geq 1$ and $0<b^2 < b< 1/2$. Combining (\ref{sleq1claim}) and (\ref{sgeq1claim}) gives the result.
\end{proof}

\noindent {{\bf Acknowledgements}. The authors would like to thank Eric Carlen, Jos\'e Antonio Carrillo, Rupert Frank, Yao Yao, and Xu Yang  their helpful comments and suggestions. They would also like to thank the two referees for their careful reading of the paper and thoughtful insights.}

\bibliographystyle{IEEEtranS}
\def\url#1{}
\bibliography{bibl_CrTo}

\end{document}